\documentclass{amsart}
\usepackage{amssymb}
\usepackage{mathrsfs}
\usepackage{stmaryrd}
\usepackage{bbm}
\usepackage{oldgerm}
\usepackage[english]{babel}
\usepackage[T1]{fontenc}
\usepackage[latin1]{inputenc}
\usepackage[all]{xy}

\newtheorem{teo}[subsection]{Theorem}
\newtheorem{prop}[subsection]{Proposition}
\newtheorem{cor}[subsection]{Corollary}
\newtheorem{lem}[subsection]{Lemma}

\theoremstyle{definition}

\newtheorem{defi}[subsection]{Definition}
\newtheorem{rema}[subsection]{Remark}
\newtheorem{remas}[subsection]{Remarks}

\numberwithin{equation}{subsection}

\newcommand{\mR}{{\mathbb R}}
\newcommand{\mC}{{\mathbb C}}
\newcommand{\mQ}{{\mathbb Q}}

\newcommand{\mZ}{{\mathbb Z}}
\newcommand{\mF}{{\mathbb F}}
\newcommand{\mA}{{\mathbb A}}

\newcommand{\mG}{{\mathbb G}}

\newcommand{\bD}{{\bf D}}

\newcommand{\bV}{{\bf V}}
\newcommand{\bO}{{\bf O}}

\newcommand{\rt}{{\rm t}}
\newcommand{\rw}{{\rm w}}

\newcommand{\Spec}{{\rm Spec}}
\newcommand{\rig}{{\rm rig}}

\newcommand{\fil}{{\rm fil}}
\newcommand{\res}{{\rm res}}
\newcommand{\ord}{{\rm ord}}

\newcommand{\Frob}{{\rm Frob}}
\newcommand{\Supp}{{\rm Supp}}
\newcommand{\Tr}{{\rm Tr}}
\newcommand{\rN}{{\rm N}}
\newcommand{\Ind}{{\rm Ind}}
\newcommand{\Rec}{{\rm Rec}}
\newcommand{\inv}{{\rm inv}}
\newcommand{\ab}{{\rm ab}}
\newcommand{\Gr}{{\rm Gr}}
\newcommand{\gr}{{\rm gr}}
\newcommand{\id}{{\rm id}}
\newcommand{\rR}{{\rm R}}
\newcommand{\rA}{{\rm A}}
\newcommand{\rP}{{\rm P}}

\newcommand{\rL}{{\rm L}}
\newcommand{\rH}{{\rm H}}
\newcommand{\rF}{{\rm F}}
\newcommand{\rV}{{\rm V}}
\newcommand{\pr}{{\rm pr}}

\newcommand{\sw}{{\rm sw}}
\newcommand{\rk}{{\rm rk}}
\newcommand{\rsw}{{\rm rsw}}
\newcommand{\dimtot}{{\rm dimtot}}
\newcommand{\Hom}{{\rm Hom}}
\newcommand{\rW}{{\rm W}}

\newcommand{\oK}{\overline{K}}
\newcommand{\ok}{\overline{k}}
\newcommand{\os}{\overline{s}}
\newcommand{\oS}{\overline{S}}

\newcommand{\oa}{\overline{a}}
\newcommand{\oy}{\overline{y}}
\newcommand{\oh}{\overline{\hbar}}
\newcommand{\oeta}{{\overline{\eta}}}
\newcommand{\opsi}{{\overline{\psi}}}
\newcommand{\otau}{{\overline{\tau}}}

\newcommand{\ocL}{\overline{\cL}}
\newcommand{\ocH}{\overline{\cH}}

\newcommand{\odelta}{{\overline{\delta}}}

\newcommand{\omQ}{{\overline{\mQ}}}

\newcommand{\us}{\underline{s}}

\newcommand{\ut}{\underline{t}}

\newcommand{\uB}{\underline{B}}
\newcommand{\uX}{\underline{X}}
\newcommand{\uU}{\underline{U}}

\newcommand{\uK}{\underline{K}}
\newcommand{\ucG}{\underline{\cG}}
\newcommand{\ucM}{\underline{\cM}}

\newcommand{\uS}{\underline{S}}
\newcommand{\ueta}{\underline{\eta}}

\newcommand{\ta}{\widetilde{a}}
\newcommand{\tb}{\widetilde{b}}
\newcommand{\tc}{\widetilde{c}}
\newcommand{\tf}{\widetilde{f}}
\newcommand{\tz}{\widetilde{z}}
\newcommand{\tx}{\widetilde{x}}
\newcommand{\tu}{\widetilde{u}}
\newcommand{\tH}{\widetilde{H}}
\newcommand{\tU}{\widetilde{U}}

\newcommand{\tcG}{\widetilde{\cG}}
\newcommand{\tcF}{\widetilde{\cF}}

\newcommand{\hOmega}{\widehat{\Omega}}
\newcommand{\hfX}{\widehat{\fX}}

\newcommand{\cA}{{\mathscr A}}
\newcommand{\cB}{{\mathscr B}}
\newcommand{\cC}{{\mathscr C}}
\newcommand{\cD}{{\mathscr D}}
\newcommand{\cE}{{\mathscr E}}
\newcommand{\cF}{{\mathscr F}}
\newcommand{\cG}{{\mathscr G}}
\newcommand{\cJ}{{\mathscr J}}
\newcommand{\cK}{{\mathscr K}}
\newcommand{\cL}{{\mathscr L}}
\newcommand{\co}{{\mathscr O}}

\newcommand{\cH}{{\mathscr H}}
\newcommand{\cM}{{\mathscr M}}

\newcommand{\cQ}{{\mathscr Q}}

\newcommand{\cHom}{{\mathscr Hom}}

\newcommand{\fF}{{\mathfrak F}}
\newcommand{\fm}{{\mathfrak m}}

\newcommand{\fX}{{\mathfrak X}}
\newcommand{\fa}{{\mathfrak a}}
\newcommand{\fp}{{\mathfrak p}}

\newcommand{\tth}{{\tt h}}

\newcommand{\crP}{{\Check{\rP}}}
\newcommand{\crA}{{\Check{\rA}}}
\newcommand{\cf}{{\Check{f}}}
\newcommand{\cg}{{\Check{g}}}
\newcommand{\cz}{{\Check{z}}}
\newcommand{\ch}{{\Check{h}}}
\newcommand{\cx}{{\Check{x}}}

\newcommand{\cu}{{\Check{u}}}
\newcommand{\cv}{{\Check{v}}}
\newcommand{\cj}{{\Check{\jmath}}}
\newcommand{\cU}{{\Check{U}}}
\newcommand{\cT}{{\Check{T}}}

\newcommand{\ccK}{{\Check{K}}}

\newcommand{\cinfty}{{\Check{\infty}}}
\newcommand{\czero}{{\Check{0}}}
\newcommand{\cpr}{{\Check{\pr}}}
\newcommand{\ctau}{{\Check{\tau}}}

\begin{document}

\title{Local Fourier transform and epsilon factors}
\author{Ahmed Abbes and Takeshi Saito}
\address{A.A. CNRS UMR 6625, IRMAR, Université de Rennes 1,
Campus de Beaulieu, 35042 Rennes cedex, France}
\email{ahmed.abbes@univ-rennes1.fr}
\address{T.S. Department of Mathematical Sciences, 
University of Tokyo, Tokyo 153-8914, Japan}
\email{t-saito@ms.u-tokyo.ac.jp}

\begin{abstract}
Laumon introduced the local Fourier transform
for $\ell$-adic Galois representations of local fields, of equal characteristic $p$ different from $\ell$,
as a powerful tool to study the Fourier-Deligne transform of $\ell$-adic sheaves over the affine line. 
In this article, we compute explicitly the local Fourier transform of monomial representations 
satisfying a certain ramification condition, and deduce Laumon's formula relating
the $\varepsilon$-factor to the determinant of the local Fourier transform under the same condition.
\end{abstract}

\maketitle

\section{Introduction}

\subsection{}\label{intro1}
In his seminal article \cite{laumon}, Laumon introduced the local Fourier transform
for $\ell$-adic Galois representations of local fields, of equal characteristic $p$ different from $\ell$,
providing a powerful tool to study the Fourier-Deligne transform of $\ell$-adic sheaves over the affine line. 
He used it to prove that the constant of the functional equation
of the $L$-function associated to an $\ell$-adic representation of a function field 
is a product of local constants, also known as $\varepsilon$-factors. 
As a key step, he gave a cohomological interpretation of the $\varepsilon$-factor in terms 
of the determinant of the local Fourier transform. In this article, we compute explicitly
the local Fourier transform of monomial representations satisfying a certain ramification 
condition, and deduce Laumon's formula for $\varepsilon$-factors under the same condition.
Our approach, inspired by our ramification theory \cite{aml}, is local and geometric, 
while Laumon's approach is global, combining arithmetic and geometric arguments. 

\subsection{}\label{intro2}
One of the main innovations of \cite{laumon}, leading to the local Fourier transform,
is Laumon's {\em principle of stationary phase},
which has its origins in the classical theory of asymptotic integrals (cf. \cite{katz2}). 
We briefly recall the classical theory (\cite{dieudonne} IV §4).
Given two functions $\varphi\in \cC^\infty(\mR,\mR)$ and $f\in \cC^\infty_c(\mR,\mC)$, 
we are interested in studying the asymptotic behavior at $\infty$ of the integral,
depending on a real parameter $t$,  
\[
I(t)=\int f(x)e^{it \varphi(x)}dx.
\]   
If the derivative of $\varphi$ does not vanish at any point in $\Supp(f)$, 
then $I(t)$ is rapidly decreasing at $\infty$.
It follows that if $\varphi$ has only finitely many critical points in $\Supp(f)$, 
then the asymptotic behavior of $I(t)$ at $\infty$ is a finite sum of contributions, 
one from each critical point of $\varphi$ in $\Supp(f)$. If moreover all critical points of $\varphi$ 
are non-degenerate ({\em i.e.} the second derivative of $\varphi$ does not vanish at these points),
then one can give a very explicit description of $I(t)$ as $t$ tends to $\infty$. 

\subsection{}\label{intro3}
Let $k$ be a perfect field of characteristic $p$, $\rA=\Spec(k[x])$ and $\crA=\Spec(k[\cx])$
be the affine lines over $k$ (equipped with coordinates $x$ and $\cx$), 
$\ell$ be a prime number different from $p$, $\psi_0\colon \mF_p\rightarrow \omQ_\ell^\times$
be a non-trivial additive character.  We denote by $\rP$ and $\crP$ the projective lines
over $k$, completions of $\rA$ and $\crA$ respectively, and by $\infty\in \rP$ and
$\cinfty\in \crP$ the points at infinity. 
For closed points $z\in \rP$ and $\cz\in \crP$, we denote by $T_z$ and $\cT_{\cz}$
the henselizations of $\rP$ and $\crP$ at $z$ and $\cz$ respectively, 
and by $\tau_z$ and $\ctau_{\cz}$ their generic points. 
Let $\cF$ be a $\omQ_\ell$-sheaf over $\rA$. The analogue of the integral $I(t)$ is 
provided by the Fourier-Deligne transform $\fF_{\psi_0}(\cF)$, which is a complex of 
$\ell$-adic sheaves on $\crA$ (cf. \ref{dtf1}); in fact, the precise analogue of $I(t)$ 
is the sheaf $\cH^1(\fF_{\psi_0}(\cF))$ (where $t$ is replaced by $\cx$). 
The ``asymptotic behavior'' of this sheaf at $\cinfty$
is encoded in its restriction to $\ctau_{\cinfty}$, which corresponds
to an $\ell$-adic representation of the absolute Galois group of the function field 
$k(\ctau_\cinfty)$ of $\cT_\cinfty$.   
Let $U$ be a dense open subscheme of $\rA$ such that $\rA-U\subset \rA(k)$.
If $\cF$ is the extension by $0$ of a smooth $\omQ_\ell$-sheaf over $U$, 
Laumon proved that we have a canonical decomposition 
\begin{equation}\label{intro3a}
\cH^1(\fF_{\psi_0}(\cF))|\ctau_{\cinfty}\simeq \bigoplus_{z\in \rP-U} \fF_{\psi_0}^{(z,\cinfty)}(\cF|\tau_z),
\end{equation}
where the factor $\fF_{\psi_0}^{(z,\cinfty)}(\cF|\tau_z)$ is the {\em local Fourier transform} of 
$\cF|\tau_z$ at $(z,\cinfty)$. The latter transformation is a functor from $\omQ_\ell$-sheaves
over $\tau_z$ to $\omQ_\ell$-sheaves over $\ctau_\cinfty$, defined by Laumon 
using vanishing cycles (cf. \ref{results2}).

\subsection{}\label{intro5}
We assume in the sequel that $p>2$.
Let $S$ be the spectrum of a henselian discrete valuation ring, $\eta$ (resp. $s$) be its generic 
(resp. closed) point, $v\colon S\rightarrow \rA$ and $\cv\colon S\rightarrow \crP$ 
be two morphisms such that $\cv(s)=\cinfty$. We put  $z=v(s)$ and assume for simplicity 
(only in the introduction) that $z\in \rA(k)$.
We denote by $f\colon S\rightarrow T_z$ and $\cf\colon S\rightarrow \cT_\cinfty$ 
the morphisms induced by $v$ and $\cv$ respectively (cf. \ref{intro3}).
Assume that $f$ and $\cf$ are finite and étale at $\eta$.
Let $\cG$ be a $\omQ_\ell$-sheaf of rank $1$ over $\eta$. 
Our main theorem \eqref{results3} says that if $(\cG,f,\cf)$ is a {\em Legendre triple}, 
a condition that will be defined below, then we have a canonical isomorphism of sheaves over $\ctau_{\cinfty}$
\begin{equation}\label{intro5a}
\fF_{\psi_0}^{(z,\cinfty)}(f_*\cG)\simeq \cf_*\left(\cG\otimes \cL_{\psi_0}(bc)\otimes 
\cK\left(-\frac{1}{2}\frac{dc}{db}\right)\otimes \cQ\right),
\end{equation}
where the rank $1$ sheaf between brackets on the right hand side is defined as follows~: 
the pull-backs of the coordinates $x$ and $\cx$
by $f$ and $\cf$ define two functions on $\eta$, denoted respectively by $b$ and $c$. 
The sheaf $\cL_{\psi_0}(bc)$ is the Artin-Schreier sheaf over $\eta$ associated to the additive 
character $\psi_0$ and the function $bc$ (cf. \ref{results1}). The sheaf $\cK\left(-\frac{1}{2}\frac{dc}{db}\right)$
is the Kummer sheaf over $\eta$ associated to the unique non-trivial character 
$\kappa_0\colon \mu_2(k)\rightarrow \omQ_\ell^\times$
and the function $-\frac{1}{2}\frac{dc}{db}$ (cf. \ref{kummer} and \ref{results3b}).
Finally, $\cQ$ is the $\omQ_\ell$-sheaf of rank $1$ over $\Spec(k)$ corresponding to the
quadratic  Gauss sum defined by $\psi_0$ and $\kappa_0$ \eqref{kummer}.

We also prove variants of \eqref{intro5a} for the local Fourier transforms $\fF_{\psi_0}^{(\infty,\cinfty)}$
and $\fF_{\psi_0}^{(\infty,\czero)}$ \eqref{results3c}.

\subsection{}\label{intro6}
The notion of a Legendre triple relies on ramification theory. 
Let $R$ be the completion of the local ring of $S$, $K$ be its fraction field, $t$ be a uniformizer of $R$,
$\ord$ be the valuation of $K$ normalized by $\ord(t)=1$. For any $y\in K$, we put $y'=\frac{dy}{dt}$.  
Ramification theory of Artin-Schreier-Witt sheaves over $\Spec(K)$
is described in terms of the Kato filtration on the ring $\rW_{m+1}(K)$ 
of Witt vectors of length $m+1$, and the homomorphism of de Rham-Witt 
$\rF^m d\colon \rW_{m+1}(K)\rightarrow \Omega^1_K$. We refer to §\ref{calculus} for a short 
review of this theory, and to \cite{aml,kato1} for more details. 

Let $a=(a_0,\dots,a_m)$ be a non-zero element of $\rW_{m+1}(K)$, $\alpha$ be the element of $K$
defined by the equation $\rF^m d(a)=\alpha dt$ \eqref{calculus2b}, $b,c$ be non-zero elements of $K$ (for which we
will take the functions provided by \ref{intro5}). 
We say that $(a,b,c)$ is a {\em strong Legendre triple} if the following 
relations are satisfied~:
\begin{eqnarray}
&p^{m-i}\ord(a_i)\geq -n=\ord(t\alpha)& (\forall\ 0\leq i\leq m),\label{intro6b}\\
&\rF^md(a)+cdb=0,& \label{intro6a}\\
&2 \ord(tb'/b)+p \ord(tc'/c)<(p-2)n.&\label{intro6c}
\end{eqnarray} 
The inequalities in \eqref{intro6b} mean that $a$ belongs to the level $n$ of the Kato filtration 
of $\rW_{m+1}(K)$; the equality in \eqref{intro6b} 
implies that $n$ is the Swan conductor of the sheaf of rank $1$ over $\Spec(K)$ defined by $a$.
We may consider equation \eqref{intro6a} as an analogue of the Legendre transform, 
used in the method of the saddle point (\cite{dieudonne} IX §1), and   
the ramification condition \eqref{intro6c} as an analogue of the convexity condition 
required for this method (or equivalently, the non-degenerate critical points
condition in the principle of stationary phase). 

In general, $a$ and $b$ are given, and $c$ will be defined by the equations above. 
To allow more flexibility, we replace \eqref{intro6a} by the following weaker but sufficient condition~:
\begin{equation}\label{intro6d}
2\ord(\alpha+cb')\geq -n+\ord(tc'/c); 
\end{equation}
we say then that $(a,b,c)$ is a {\em Legendre triple} (cf. \ref{legendre}). 

\subsection{}\label{intro7}
We take again the notation and assumptions of \ref{intro5}. 
We say that $(\cG,b,c)$, or $(\cG,f,\cf)$, 
is a {\em Legendre triple} if we can write $\cG$ as a tensor product 
of two $\omQ_\ell$-sheaves of rank $1$ over $\eta$, $\cG\simeq \cG_\rt\otimes\cG_\rw$, 
where $\cG_\rt$ is tamely ramified at $s$ and $\cG_\rw$ is trivialized by a cyclic extension 
of order $p^{m+1}$ of $\eta$ ($m\geq 0$) and satisfies the following conditions~: there exists 
$a\in \rW_{m+1}(K)$ such that $(a,b,c)$ is a Legendre triple and the pull-back of $\cG_\rw$ 
to $\Spec(K)$ is associated to $a$ (cf. \ref{legendre1}). It follows in particular
that $\cG$ is wildly ramified.

Suppose given the pair $(\cG,f)$, we would like to compute the local Fourier transform
$\cF^{(z,\cinfty)}_{\psi_0}(f_*\cG)$. In order to apply \eqref{intro5a},  
since the morphism $\cf$ is completely determined by $c$,
the problem is to find a non-zero function $c$ over $\eta$  such that $(\cG,b,c)$ is a Legendre triple. 
It is clear that we can first choose $a$ satisfying \eqref{intro6b}, 
and then choose $c$ satisfying \eqref{intro6d}; but in general, $c$ may not satisfy \eqref{intro6c}. 
In fact, there are pairs $(\cG,b)$ such that the sheaf $\cF^{(z,\cinfty)}_{\psi_0}(f_*\cG)$ is not monomial;
therefore, \eqref{intro5a} implies that there is no $c$ such that $(\cG,b,c)$ is a Legendre triple
for such pairs $(\cG,b)$. 
On the other side, there are extreme cases for which equation \eqref{intro6c} is implied by the two others,
and hence there exists $c$ such that $(\cG,b,c)$ is a Legendre triple. 
Indeed, equations \eqref{intro6b} and \eqref{intro6d} imply that we have the following relation 
\begin{equation}\label{intro7a}
\deg(\cf)=\sw(f_*\cG)+\deg(f),
\end{equation}
where $\deg(-)$ is the degree and $\sw(-)$ is the Swan conductor  (cf. \ref{results4a}). 
If $f$ is tamely ramified, $\sw(f_*\cG)=\sw(\cG)\geq 1$ and $\sw(\cG)+\deg(f)$ is prime to $p$, 
then $\cf$ is tamely ramified, and hence equation \eqref{intro6c} is satisfied.
Notice that under the assumptions on tameness, we have $\ord(tb'/b)=\ord(tc'/c)=0$. 

A special case of formula \eqref{intro5a} was conjectured by Laumon and Malgrange (\cite{laumon} 2.6.3)
and proved by Fu \cite{fu}. It corresponds to the extreme case where $\cG_\rw$ is an Artin-Schreier sheaf 
({\em i.e.}, $m=0$) of Swan conductor $s$, and $f$ is a tamely ramified morphism of degree $r$,
such that $1\leq s<p$, $r+s$ is prime to $p$, and $k$ is algebraically closed.

\subsection{}\label{intro8}
The main idea and the key technical tool for the proof of \eqref{intro5a} come 
from our theory of ramification \cite{aml}. We denote by $\pr_1$ and $\pr_2$ the canonical 
projections of $\eta\times_k\eta$, by $j\colon \eta\times_k\eta\rightarrow S\times_kS$ 
the canonical injection, and put $b_1=\pr_1^*(b)$ and $c_2=\pr_2^*(c)$. 
The proof of \eqref{intro5a} is made in two steps.
The first and most important one is a calculus of vanishing cycles under the Legendre conditions. 
The second one is a computation of the dimension of the local Fourier transform, which 
holds in general without any restriction on the sheaf.

First, we study the complex of vanishing cycles 
of the sheaf $j_!(\pr_1^*(\cG)\otimes\cL_{\psi_0}(b_1c_2))$ relatively to the second projection
$S\times_kS\rightarrow S$. By adapting our method in \cite{aml}, 
we prove that, under the Legendre conditions, this complex can be explicitly 
described over an open subscheme of a suitable blow-up of $S\times_kS$ along a closed subscheme 
of the diagonal $S\rightarrow S\times_kS$ (cf. \ref{nbc6}). Figuratively speaking, 
we kill the ramification by blowing-up in the diagonal, which was the leitmotiv of \cite{as,aml}. 
From this, we deduce that the sheaf $\cf^*(\fF_{\psi_0}^{(z,\cinfty)}(f_*\cG))$ over $\eta$ 
has a direct factor isomorphic to  
\[
\cD=\cG\otimes \cL_{\psi_0}(bc)\otimes \cK\left(-\frac{1}{2}\frac{dc}{db}\right)\otimes \cQ,
\]
and the morphism $\cf_*(\cD)\rightarrow \fF_{\psi_0}^{(z,\cinfty)}(f_*\cG)$, 
induced by the trace homomorphism $\cf_*\cf^*\rightarrow \id$, is injective (cf. \ref{nbc3}).

Second, in order to prove that the morphism $\cf_*(\cD)\rightarrow \fF_{\psi_0}^{(z,\cinfty)}(f_*\cG)$
is an isomorphism, it is enough to show that the rank of $\cF^{(z,\cinfty)}_{\psi_0}(f_*\cG)$ 
is equal to the degree of $\cf$. 
By \eqref{intro7a}, the latter relation is a special case of a general formula proved 
by Laumon (\cite{laumon} 2.4.3)~: 
namely, for any $\omQ_\ell$-sheaf $\cF$ over $\tau_z$, we have
\begin{equation}\label{intro8a}
\rk(\cF^{(z,\cinfty)}_{\psi_0}(\cF))=\sw(\cF)+\rk(\cF).
\end{equation}
We give in the appendix \eqref{apdlft6} 
another proof of this equation using a formula of Deligne-Kato that computes 
the dimension of the nearby cycle complex 
of a sheaf on a smooth curve over a strictly henselian trait. 
Deligne considered the case where the sheaf has no vertical ramification (\cite{laumon2} 5.1.1), 
and Kato extended the formula to the general case (\cite{kato3} 6.7).  
We give in the appendix \eqref{app12} a brief review of Kato's formula for rank $1$ sheaves,
which is enough for our application, by using his refined Swan conductors. The latter 
fits perfectly in our ramification theory as proved in \cite{aml}, and hence in the 
general philosophy of this article.
 
\subsection{}\label{intro9}
Formula \eqref{intro5a} has strong relations with the theory of $\varepsilon$-factors.
First, it was suggested by explicit formulas for $\varepsilon$-factors of quasi-characters \eqref{lef7c}
(cf. also \cite{henniart}). Second, it implies Laumon's formula relating $\varepsilon$-factors
and local Fourier transforms. More precisely, 
if $k$ is finite, $\cF$ is a $\omQ_\ell$-sheaf over $\tau_0$ and  
$\cF_!$ is the extension of $\cF$ by $0$ to $T_0$, then Laumon (\cite{laumon} 3.6.2)
proved that we have
\begin{equation}\label{intro9a}
(-1)^d\det(\Rec_{\cT_\cinfty}(\cx^{-1}),\fF^{(0,\cinfty)}(\cF))=\varepsilon(T_0,\cF_!,dx),
\end{equation}
where $d$ is the dimension of 
$\fF^{(0,\cinfty)}_{\psi_0}(\cF)$, $\varepsilon(T_0,\cF_!,dx)$ is the $\varepsilon$-factor 
of the sheaf $\cF_!$ over $T_0$, and $\Rec_{\cT_\cinfty}$ is the reciprocity isomorphism of class field theory 
for the completion of the function field $k(\ctau_\cinfty)$ of $\cT_\cinfty$ (cf. \ref{lauf1}). 

Under the assumptions of \eqref{intro5}, if $(\cG,f,\cf)$ is a Legendre triple and $z=0$,
we give in \ref{lauf2} a new proof of \eqref{intro9a} for the sheaf $\cF=f_*(\cG)$. 
We deduce it from \eqref{intro5a} by using three ingredients. 
The first one is a classical explicit formula for $\varepsilon$-factors involving Gauss sums \eqref{lef7}.  
The second one is a variation on Witt's explicit reciprocity law due to Fontaine \eqref{lef9}.
The third ingredient is new \eqref{lef8}; it is an explicit formula for the Langlands $\lambda$-factor 
which appears in the induction formula for $\varepsilon$-factors.
We prove the latter by using Deligne's formula for the $\varepsilon$-factor 
of an orthogonal representation in terms of its second Stiefel-Whitney class (\cite{deligne2} 1.5), 
and Serre's formula for the second Stiefel-Whitney class of induced representations    
in terms of the Hasse-Witt invariant of quadratic forms \cite{serre2}.

\subsection{}
This article is divided into two parts and augmented by two appendices.
The first part, with a strong geometric flavor, is devoted to the proof of \eqref{intro5a}. 
Section \ref{calculus} develops the necessary tools from ramification 
theory of Artin-Schreier-Witt sheaves. It contains in particular a computation 
of Witt vectors \eqref{key1} that plays a crucial role in the sequel. 
In §\ref{results}, we review the definition of 
the local Fourier transform and state the main theorems \ref{results3} and \ref{results3c}.
Section \ref{nbc} is the heart of the article. It contains the analysis of a complex 
of vanishing cycles by blowing-up in the diagonal mentioned in \ref{intro8}. 
The proofs of the main theorems are given in §\ref{pfgm}. 
The second part, with a more arithmetic flavor, is devoted to the proof of \eqref{intro9a}.
It starts by a brief review of Stiefel-Whitney classes and a formula of Serre in §\ref{stwh}, 
followed by a short complement on refined logarithmic differents in §\ref{refdif}.
In  §\ref{lef}, we review the theory of $\varepsilon$-factors and 
develop the necessary ingredients for the proof of \eqref{intro9a}.
Finally, the proof of this formula is completed in §\ref{lauf}.
The first appendix §\ref{scsc} reviews the Deligne-Kato's formula for the dimension of 
the nearby cycle complex of a sheaf of rank $1$ on a smooth curve over a strictly henselian trait. 
In the second appendix §\ref{apdlft}, we apply this formula to compute the dimension
of the local Fourier transform.  
 
\subsection{}
C. Sabbah proved an explicit formula for the local Fourier transform 
of a formal germ of meromorphic connection of one complex variable using a blow-up technique \cite{sabbah}. 
The relation between our approaches is not clear.  
During the preparation of this article, we learned from M. Strauch that he made some expectations 
on a local principle of stationary phase, without giving precise formulas. He is motivated by 
applications to the cohomology of Lubin-Tate spaces.
We are grateful to an anonymous referee for his thorough reading of the manuscript  
and helpful comments.
The first author would like to acknowledge the hospitality of the 
Department of Mathematical Sciences at the University of Tokyo where this work was achieved.

\section*{Notation and conventions}

\subsection{}\label{not1}
In this article (except in §\ref{refdif}, §\ref{lef} and the appendix), we fix a prime number $p>2$, 
a perfect field $k$ of characteristic $p$ and an algebraic closure $\ok$ of $k$. 
For $q$ a power of $p$, we denote by $\mF_q$ the unique subfield of $\ok$
with $q$ elements. We fix also a prime number $\ell$ different from $p$, 
an algebraic closure $\omQ_\ell$ of the field $\mQ_\ell$ of $\ell$-adic numbers
and a non-trivial additive character $\psi_0\colon \mF_p\rightarrow \omQ_\ell^\times$. 
For every integer $m\geq 0$, we fix an injective homomorphism $\psi_m\colon \mZ/p^{m+1}\mZ\rightarrow \omQ_\ell^\times$
such that for any $a\in \mF_p$, we have $\psi_m(p^m a)=\psi_0(a)$, 
where $p^m a$ denotes the embedding $\mF_p\rightarrow \mZ/p^{m+1}\mZ$
induced by the multiplication by $p^m$ on $\mZ$. 

\subsection{}\label{not2}
If $X$ is a scheme and $x\in X$, we denote by $\kappa(x)$ the residue field of $X$ at $x$
and by $i_x\colon \Spec(\kappa(x))\rightarrow X$ the canonical morphism.

\subsection{}\label{not3}
For a scheme $X$, a ``$\omQ_\ell$-sheaf over $X$'' stands for a ``constructible $\omQ_\ell$-sheaf over $X$''
in the sense of (\cite{weil2} 1.1.1). We denote by $\bD^b_c(X,\omQ_\ell)$
the derived category of $\ell$-adic sheaves defined in ({\em loc. cit.}, 1.1.2 and 1.1.3).

\section{Calculus on Witt vectors}\label{calculus}

\subsection{}\label{calnt}
Let $R$ be a complete discrete valuation ring of equal characteristic $p$, 
with residue field $k$, equipped with a uniformizer $t$, $K$ be the fraction field of $R$, 
$\ord$ be the valuation of $K$ normalized by $\ord(t)=1$. 
We identify $R$ with the ring of power series $k[[t]]$. For any $x\in K$, we denote
by $x^{(i)}$ the $i$-th iterated derivative of $x$ relatively to $t$ $(i\geq 1$); we put $x^{(0)}=x$ and $x'=x^{(1)}$.

\subsection{}
The module $\Omega^1_R$ is free of rank $1$ over $R$,
and hence complete and separated. We identify it with a submodule of $\Omega^1_K$. 
For $a\in K^\times$, we put $d\log(a)=da/a \in \Omega^1_K$. We denote by $\Omega^1_R(\log)$
the sub-$R$-module of $\Omega^1_K$ generated by $\Omega^1_R$ and the elements of the form $d\log a$ 
for $a\in R-\{0\}$. Then $\Omega^1_R(\log)$ is a free $R$-module of rank $1$ generated by 
$d\log(t)$. We put $\Omega^1_k(\log)=\Omega^1_R(\log)\otimes_Rk=k\cdot d\log(t)$.
We define an increasing exhaustive filtration on $\Omega^1_K$ by setting, for $n\in \mZ$, $\fil_n\Omega^1_K=t^{-n}\Omega^1_R(\log)$. We have 
\[
\Gr_n\Omega^1_K=\fil_n\Omega^1_K/\fil_{n-1}\Omega^1_K\simeq (t^{-n}R/t^{-n+1}R)\cdot d\log(t).
\]

\subsection{}\label{calculus1}
Let $m$ be an integer $\geq 0$, $\rW_{m+1}(K)$ be the ring of Witt vectors of length $m+1$.  
Following \cite{bry,kato1}, we define an increasing exhaustive filtration on the group of Witt vectors
$\rW_{m+1}(K)$ by setting, for $n\in \mZ$, $\fil_n\rW_{m+1}(K)$ to be the subgroup of elements 
$(x_0,\dots,x_m)$ such that 
\begin{equation}\label{calculus1a}
p^{m-i}\ord(x_i)\geq -n \ \ {\rm for\ all}\ \ 0\leq i\leq m.
\end{equation}
We put 
\[
\Gr_n\rW_{m+1}(K)=\fil_n\rW_{m+1}(K)/\fil_{n-1}\rW_{m+1}(K).
\]
Let $\rV\colon \rW_{m+1}(K)\rightarrow \rW_{m+2}(K)$ be the verschiebung endomorphism. 
We have 
\[
\rV(\fil_n\rW_{m+1}(K))\subset \fil_n\rW_{m+2}(K).
\]

\subsection{}\label{calculus2}
Let $\rF\colon \rW_{\bullet+1}\Omega^1_K\rightarrow \rW_{\bullet}\Omega^1_K$
be the Frobenius endomorphism of the de Rham-Witt complex of $K$ over $k$. 
The homomorphism 
\begin{equation}\label{calculus2a}
\rF^m d\colon \rW_{m+1}(K)\rightarrow \Omega^1_K
\end{equation}
is given by the formula 
\begin{equation}\label{calculus2b}
\rF^m d(x_0,\dots,x_m)=\sum_{i=0}^mx_i^{p^{m-i}-1}dx_i.
\end{equation}
Therefore, for any integer $n$, we have 
\begin{equation}\label{calculus2c}
\rF^md(\fil_n\rW_{m+1}(K))\subset \fil_n\Omega^1_K.
\end{equation}
We deduce a canonical homomorphism 
\begin{equation}\label{calculus2d}
\gr_n(\rF^md):\Gr_n\rW_{m+1}(K)\rightarrow \Gr_n\Omega^1_K.
\end{equation}

\subsection{}\label{rappelaml}
The exact sequence 
\begin{equation}\label{deltam0}
\xymatrix{
0\ar[r]&{\mZ/p^{m+1}\mZ}\ar[r]& {\rW_{m+1}}\ar[r]^{\rF-1}& {\rW_{m+1}}\ar[r]&0}
\end{equation}
induces a surjective homomorphism
\begin{equation}\label{deltam}
\delta_{m+1}:\rW_{m+1}(K)\rightarrow \rH^1(K,\mZ/p^{m+1}\mZ).
\end{equation}
We define an increasing exhaustive filtration on $\rH^1(K,\mZ/p^{m+1}\mZ)$ by setting (for $n\in \mZ$) 
\begin{equation}
\fil_n\rH^1(K,\mZ/p^{m+1}\mZ)=\delta_{m+1}(\fil_n\rW_{m+1}(K)).
\end{equation}
We put 
\[
\Gr_n\rH^1(K,\mZ/p^{m+1}\mZ) =\fil_n\rH^1(K,\mZ/p^{m+1}\mZ)/\fil_{n-1}\rH^1(K,\mZ/p^{m+1}\mZ).
\]
By (\cite{kato1} 3.2, \cite{aml} 10.7), for any integer $n\geq 1$, there exists a unique homomorphism 
\begin{equation}
\psi_{m,n}\colon\Gr_n\rH^1(K,\mZ/p^{m+1}\mZ)\rightarrow \Gr_n\Omega^1_K
\end{equation}
making the following diagram commutative~:
\begin{equation}
\xymatrix{
{\Gr_n\rW_{m+1}(K)}\ar[rr]^-(0.5){-\gr_n(\rF^md)}\ar[d]_{\gr_n(\delta_{m+1})}&&{\Gr_n\Omega^1_K}\\
{\Gr_n\rH^1(K,\mZ/p^{m+1}\mZ)}\ar[rru]_-(0.5){\psi_{m,n}}&&}
\end{equation}

For any $\chi\in \rH^1(K,\mZ/p^{m+1}\mZ)$, the Swan conductor of $\chi$, $\sw(\chi)$,
is the smallest integer $n\geq 0$ such that $\chi\in \fil_n \rH^1(K,\mZ/p^{m+1}\mZ)$ 
(cf. \cite{bry} cor. of theo.~1, \cite{kato1}). Kato defined the {\em refined Swan conductor} of $\chi$, 
$\rsw(\chi)$, as the image of the class of $\chi$ by the homomorphism 
\[
\psi_{m,\sw(\chi)}\colon\Gr_{\sw(\chi)}\rH^1(K,\mZ/p^{m+1}\mZ)\rightarrow \Gr_{\sw(\chi)}\Omega^1_K.
\]

\subsection{}\label{key-not}
Let $k(\theta)$ be the field of rational functions in one variable $\theta$ over $k$, 
$R_L=k(\theta)[[t]]$ be the ring of power series in the variable $t$ over $k(\theta)$, 
$L$ be the fraction field of $R_L$. 
We consider $L$ as an extension of $K$ by the $k$-homomorphism $v\colon K\rightarrow L$ defined by $v(t)=t$.
Let $r$ be an integer $\geq 1$, $u\colon K\rightarrow L$ be the $k$-homomorphism defined by 
$u(t)=t(1+t^r\theta)$. In (\cite{aml} 13.4), we proved that for any integer $n$, the group homomorphism 
\begin{equation}
u-v\colon \rW_{m+1}(K)\rightarrow \rW_{m+1}(L)
\end{equation}
maps $\fil_n\rW_{m+1}(K)$ to $\fil_{n-r}\rW_{m+1}(L)$, and we expressed the induced homomorphism 
on the graded pieces
\[
\Gr_n\rW_{m+1}(K)\rightarrow \Gr_{n-r}\rW_{m+1}(L).
\]
We refine this result as follows~:

\begin{prop}\label{key1}
Let $n$ be an integer, $a=(a_0,\dots,a_m)\in \fil_n\rW_{m+1}(K)$, $\alpha$ be the element of $K$ 
such that $\rF^m d a =\alpha dt$. Then $\ord(t\alpha)\geq -n$ and we have 
\begin{equation}\label{key1a}
u(a)-a\equiv \rV^m\left(\sum_{i=1}^{p-1}\frac{t^i\alpha^{(i-1)}}{i!} (t^r\theta)^i\right)  
\mod \fil_{n-pr}\rW_{m+1}(L).
\end{equation}
\end{prop}

First, we prove some preliminary results. We define a sequence of polynomials $(n\geq 0)$ 
\[
Q_n\in \mZ[\frac{1}{p}][X_0,\dots,X_n,Y_0,\dots,Y_n]
\] 
by the inductive formula
\begin{equation}\label{key1b}
\sum_{i=0}^np^i(X_i(1+Y_i))^{p^{n-i}}=\sum_{i=0}^np^iX_i^{p^{n-i}}+\sum_{i=0}^np^i Q_i^{p^{n-i}}.
\end{equation}
Observe that for a commutative ring $A$ and elements $x=(x_0,\dots,x_m)$, $y=(y_0,\dots,y_m)$  
and $z=(z_0,\dots,z_m)$ of $\rW_{m+1}(A)$ such that $z_i=x_i(1+y_i)$ for all $0\leq i\leq m$, we have
\begin{equation}
z-x=(Q_0(x,y),Q_1(x,y),\dots,Q_m(x,y)).
\end{equation}

We denote by $\lambda(Y)$ the shifted $p$-truncated logarithm, {\em i.e.}, the polynomial of $\mZ_{(p)}[Y]$
defined by
\begin{equation}\label{key3a}
\lambda(Y)=\sum_{i=1}^{p-1}(-1)^{i+1}\frac{Y^i}{i}.
\end{equation}

\begin{lem}\label{key2}
{\rm (i)}\ The polynomials $Q_n$ belong to the ideal of $\mZ[X_0,\dots,X_n,Y_0,\dots,Y_n]$ generated 
by $(Y_0,\dots,Y_n)$. 

{\rm (ii)}\ If we attach the weight $p^i$ to the variable $X_i$ and the weight $0$ to the variable $Y_i$,
the polynomial $Q_n$ is homogeneous of weight $p^n$. 

{\rm (iii)}\ We have the following relation in $\mZ_{(p)}[X_0,\dots,X_n,Y_0,\dots,Y_n]$
\begin{equation}\label{key3}
Q_n\equiv \sum_{i=0}^{n-1} X_i^{p^{n-i}}\lambda(Y_i)+ X_nY_n \mod (p)+ (Y_0,\dots,Y_n)^p.
\end{equation}
\end{lem}

Propositions (i) and (ii) are easy. We prove (iii).
Since $Q_i$ belongs to the ideal $(Y_0,\dots,Y_i)$, we have 
\[
p^nQ_n\equiv\sum_{i=0}^n p^iX_i^{p^{n-i}}((1+Y_i)^{p^{n-i}}- 1) \mod (Y_0,\dots,Y_n)^p.
\]
So the required relation follows from the following congruence, for $i\leq n-1$ and $1\leq j\leq p-1$,
\[
\frac{1}{p^{n-i}}\left( \begin{matrix}p^{n-i}\\j\end{matrix}\right)\equiv (-1)^{j+1}\frac{1}{j}\mod p.
\]

\begin{lem}\label{key4}
{\rm (i)}\  For any $x\in K$, we have 
\begin{equation}\label{key4a}
u(x)-x\equiv \sum_{i=1}^{p-1} \frac{t^ix^{(i)}}{i!} (t^r\theta)^i\mod x t^{pr} R_L.
\end{equation}

{\rm (ii)}\ For any $x\in K^\times$, we have
\begin{equation}\label{key4b}
\ord\left(\frac{u(x)}{x}-1\right)\geq r.
\end{equation}

{\rm (iii)}\ For any $x\in K^\times$, if we put $y=\frac{x'}{x}$, we have 
\begin{equation}\label{key4c}
\lambda\left(\frac{u(x)}{x}-1\right)\equiv 
\sum_{i=1}^{p-1} \frac{t^iy^{(i-1)}}{i!}(t^r\theta)^i\mod t^{pr} R_L.
\end{equation}
\end{lem}

(i) Since $\ord(t x')\geq \ord(x)$ for any $x\in K$,
we are reduced by the Leibniz rule to proving \eqref{key4a} for $x=t$, for $x=t^{-1}$ 
and for $x\in R^\times$.
The first two cases are obvious, and the last one follows from Taylor expansion. 

(ii) It follows immediately from (i) and the fact that $\ord(t^i x^{(i)})\geq \ord(x)$ for all $i\geq 1$. 

(iii) Since both sides of \eqref{key4c} define group homomorphisms from $K^\times$ to $R_L/t^{pr}R_L$,
it is enough to prove \eqref{key4c} for $x=t$ and for $x\in R^\times$. For $x=t$, both sides 
are equal to $\lambda(t^r\theta)$. For $x\in R^\times$, we are reduced by truncation to the case where 
$x$ is a polynomial in $t$ with a non-vanishing constant term. Then after replacing $k$ by an
algebraic closure, we are further reduced to the case where $x=1- ct$ with $c\in k$.
In this case, both sides of \eqref{key4c} are equal since  
\begin{equation}\label{key4d}
\lambda\left(\frac{-ct}{1-ct}t^r\theta\right)=\sum_{i=1}^{p-1}\frac{(-1)^{i-1}}{i}
\left(\frac{-ct}{1-ct}\right)^i(t^r\theta)^i,
\end{equation}
and for $1\leq i\leq p-1$,  
\begin{equation}\label{key4e}
\left(\frac{-c}{1-ct}\right)^{(i-1)}=(-1)^{i-1}(i-1)!\left(\frac{-c}{1-ct}\right)^i.
\end{equation}

\subsection{}\label{key5}
We can now prove \ref{key1}. We set
$b=(b_0,\dots,b_m)\in L^{m+1}$, where $b_i=\frac{u(a_i)}{a_i}-1$ if $a_i\not=0$, and $b_i=0$ if $a_i=0$;
so we have $b_i\in t^rR_L$ \eqref{key4b}. It follows from \ref{key2}(ii) and \eqref{key3} that we have 
\begin{equation}
p^{m-i}\ord(Q_i(a,b))\geq -n+ p^{m-i}r,
\end{equation}
\begin{equation}
\ord\left(Q_m(a,b)-\sum_{i=0}^{m-1} a_i^{p^{m-i}}\lambda(b_i)-u(a_m)+a_m\right) \geq -n+pr.
\end{equation}
We put $c_i=\frac{a'_i}{a_i}$ if $a_i\not=0$, and $c_i=0$ if $a_i=0$. 
It follows from \eqref{key4a} and \eqref{key4c} that we have
\begin{eqnarray}
\lefteqn{\sum_{i=0}^{m-1} a_i^{p^{m-i}}\lambda(b_i)+u(a_m)-a_m}\\
&\equiv&\sum_{i=0}^{m-1} a_i^{p^{m-i}}\sum_{j=1}^{p-1}\frac{t^jc_i^{(j-1)}}{j!}(t^r\theta)^j
+\sum_{j=1}^{p-1}\frac{t^ja_m^{(j)}}{j!}(t^r\theta)^j\mod t^{pr-n}R_L\nonumber\\
&\equiv&\sum_{j=1}^{p-1} t^j\left(\sum_{i=0}^{m-1} a_i^{p^{m-i}} c_i^{(j-1)}+a_m^{(j)}\right)
\frac{(t^r\theta)^j}{j!}\mod t^{pr-n}R_L.\nonumber
\end{eqnarray}
The proposition follows since we have by definition
\begin{equation}
\alpha=\sum_{i=0}^{m-1}a_i^{p^{m-i}}c_i+a'_m.
\end{equation}

\begin{cor}\label{key6}
We keep the notation of \eqref{key1}; moreover, let $b,c$ be non-zero elements of $K$,
$\nu(b)=\ord(tb'/b)$, $\nu(c)=\ord(tc'/c)$. Assume that $\alpha+cb'=0$ and $\nu(b)+\nu(c)<(p-2)r$. Then we have
\begin{equation}\label{key6c}
u(a)-a+\rV^m(c(u(b)-b))\in \fil_{n-\nu(c)-2r}(\rW_{m+1}(L))
\end{equation}
and
\begin{equation}\label{key6d}
u(a)-a+\rV^m(c(u(b)-b))\equiv  
\rV^m(\frac{1}{2}t\alpha \frac{tc'}{c} (t^r\theta)^2) \mod \fil_{n-\nu(c)-2r-1}(\rW_{m+1}(L)).
\end{equation}
\end{cor}

Observe first that we have
\begin{equation}\label{key6a}
\alpha'+cb^{(2)}=\alpha\frac{c'}{c}, 
\end{equation}
and for any $i\geq 2$,
\begin{equation}\label{key6b}
\ord(t^{i}\alpha^{(i-1)}+t^{i}b^{(i)}c)\geq \ord(t\alpha\frac{tc'}{c}).  
\end{equation}
Indeed, equation $\alpha+cb'=0$ implies immediately \eqref{key6a} and the following equation
\[
t^i\alpha^{(i-1)}+ct^ib^{(i)}=-\sum_{j=1}^{i-1}\binom{i-1}{j}t^i b^{(i-j)}c^{(j)}.
\]
The relation $\ord(tz')\geq \ord(z)$ for any $z\in K$, implies that 
each term of the right hand side has bigger valuation than $t^2b' c'=-t^2\alpha c'/c$.

We have $\ord(t\alpha)\geq -n$ \eqref{key1}, $\ord(bc)\geq -n-\nu(b)$ and $\ord(t^2\alpha c'/c)\geq -n+\nu(c)$. 
Hence, we deduce from \eqref{key1a} and \eqref{key4a} that we have 
\begin{eqnarray}
\lefteqn{u(a)-a+\rV^m(c(u(b)-b))\equiv} \label{key6e}\\
&&\ \ \ \ \ \rV^m\left(\sum_{i=2}^{p-1}\frac{t^i}{i!}
\left(\alpha^{(i-1)}+cb^{(i)}\right) (t^r\theta)^i\right)  
\mod \fil_{n+\nu(b)-pr}\rW_{m+1}(L).\nonumber
\end{eqnarray}
The corollary follows from \eqref{key6e}, \eqref{key6a}, \eqref{key6b} and the assumptions.

\vspace{2mm}

We can replace the condition $\alpha+cb'=0$ of \ref{key6} by a weaker condition \eqref{key7b} as follows~: 

\begin{cor}\label{key7}
We keep the notation of \eqref{key1}; moreover, let $b,c$ be non-zero elements of $K$,
$\nu(b)=\ord(tb'/b)$, $\nu(c)=\ord(tc'/c)$. Assume that the following conditions are satisfied
\begin{eqnarray}
\ord(t\alpha)=-n,\label{key7a}\\
\ord(\alpha+cb')\geq -n+\nu(c)+r,\label{key7b}\\
\nu(b)+\nu(c)<(p-2)r.\label{key7c}
\end{eqnarray}
Then we have
\begin{equation}\label{key7d}
u(a)-a+\rV^m(c(u(b)-b))\in \fil_{n-\nu(c)-2r}(\rW_{m+1}(L))
\end{equation}
and
\begin{equation}\label{key7e}
u(a)-a+\rV^m(c(u(b)-b))\equiv  
\rV^m(\frac{1}{2}t\alpha \frac{tc'}{c} (t^r\theta)^2) \mod \fil_{n-\nu(c)-2r-1}(\rW_{m+1}(L)).
\end{equation}
\end{cor}

Let $c_0$ be the element of $K$ such that $\alpha+c_0b'=0$. Since $\alpha\not=0$, we have $c_0\not=0$ and
\[
1-\frac{c}{c_0}=\frac{\alpha+cb'}{\alpha}. 
\]
We deduce that $\ord(1-c/c_0)> \nu(c)+r>0$; in particular, we have $\ord(c)=\ord(c_0)$. 
The relation $\ord(tz')\geq \ord(z)$ for any $z\in K$, implies that 
\begin{equation}\label{key7f}
\ord\left(t\left(\frac{c'}{c}-\frac{c'_0}{c_0}\right)\right)=\ord\left(t\frac{c'c_0-cc'_0}{c_0^2}\right)>\nu(c)+r,
\end{equation}
and hence we have $\nu(c)=\nu(c_0)$. By \eqref{key4a}, we have 
\[
(c-c_0)(u(b)-b)=\sum_{i=1}^{p-1}(c-c_0)t^i\frac{b^{(i)}}{i!}(t^r\theta)^i \mod (c-c_0)bt^{pr} R_L.
\]
Taking in account the relations $\ord(tz')\geq \ord(z)$ for any $z\in K$, 
$\ord(b(c-c_0))>\ord(bc_0)=-n-\nu(b)$, $\ord(b'(c-c_0))\geq -n+\nu(c)+r$ which is \eqref{key7b}, and \eqref{key7c}, 
we deduce that 
\begin{equation}\label{key7g}
\ord((c-c_0)(u(b)-b))> -n+\nu(c)+2r.
\end{equation} 
The proposition follows from \eqref{key7c}, \eqref{key7f}, \eqref{key7g} and \ref{key6}.

\subsection{}\label{legendre}
Let $a$ be a non-zero element of $\rW_{m+1}(K)$, $b,c$ be non-zero elements of $K$. 
We denote by $\alpha$ the element of $K$ such that $\rF^m da=\alpha dt$, 
$n=-\ord(t\alpha)$, $\nu(b)=\ord(tb'/b)$ and $\nu(c)=\ord(tc'/c)$. 
We say that $(a,b,c)$ is a {\em Legendre triple} if the following conditions are satisfied  
\begin{eqnarray}
a\in \fil_n\rW_{m+1}(K),\label{legendrea1}\\
2\ord(\alpha+cb')\geq -n+\nu(c),\label{legendre1b}\\
2\nu(b)+p\nu(c)<(p-2) n.\label{legendre1c}
\end{eqnarray} 
Under these conditions, $n$ is finite (as $a\not=0$), and it is 
the smallest integer such that $a\in \fil_n\rW_{m+1}(K)$. Moreover, we have $n\geq 1$ and
$\ord(tb'c)=-n$. We say that $(n,\nu(b),\nu(c))$ is the {\em conductor} of the triple $(a,b,c)$.

\begin{rema}\label{legendreb}
Under the assumptions of \eqref{legendre}, if moreover $n-\nu(c)=2r$ is even, 
then the conditions of {\em loc. cit.} are equivalent to the conditions of \ref{key7}. 
\end{rema}

\subsection{}\label{legendrec}
Let $S=\Spec(R)$, $\eta=\Spec(K)$ be the generic point of $S$,  
$b,c\in K$, $\cG$ be a $\omQ_\ell$-sheaf 
of rank $1$ over $\eta$ trivialized by a cyclic extention of order $p^{m+1}$ of $\eta$ $(m\geq 0)$. 
We denote by $\chi\in \rH^1(K,\mZ/p^{m+1}\mZ)$ the class 
such that $\psi_m^{-1}\circ \chi$ is the character associated to $\cG$ \eqref{not1}.
We say that $(\cG,b,c)$ is a {\em Legendre triple} 
if there exists $a\in \rW_{m+1}(K)$ such that $\delta_{m+1}(a)=\chi$ \eqref{deltam}, 
$a,b$ and $c$ are non-zero, and $(a,b,c)$ is a Legendre triple.

\begin{defi}\label{legendre1}
Let $X$ be a smooth connected curve over $k$ (resp. the spectrum of a henselian discrete 
valuation ring of equal characteristric $p$ and residue field $k$), 
$s$ be a closed point of $X$, $U$ be the open subscheme $X-\{s\}$ of $X$, 
$x,y\in \Gamma(U,\co_X)$, $\cG$ be a smooth $\omQ_\ell$-sheaf of rank $1$ over $U$.
Let $S$ be the spectrum of the completion of the local ring of $X$ at $s$, $\eta$ be the 
generic point of $S$, $\hbar\colon S\rightarrow X$ be the canonical map.
We say that $(\cG,x,y)$ is a {\em Legendre triple} 
at $s$ if there exist $\cG_\rt$ and $\cG_\rw$ two smooth $\omQ_\ell$-sheaves of rank $1$ over $U$
satisfying the following conditions~:

{\rm (i)}\ $\cG\simeq \cG_\rt\otimes \cG_\rw$; 

{\rm (ii)}\ $\cG_\rt$ is tamely ramified at $s$; 

{\rm (iii)}\ $\cG_\rw$ is trivialized by a cyclic extension of order $p^{m+1}$ of $U$ ($m\geq 0$);

{\rm (iv)}\ $(\hbar_U^*(\cG_\rw),\hbar_U^*(x),\hbar_U^*(y))$ is a
Legendre triple in the sense of \eqref{legendrec}.  
\end{defi}

\section{Local Fourier transform}\label{results}

\subsection{}\label{results1}
Lang's isogeny $\rL$ of $\mG_{a,k}=\Spec(k[u])$, defined by $\rL^*(u)=u^p-u$, induces 
the Artin-Schreier exact sequence 
\begin{equation}
0\rightarrow \mF_p\rightarrow \mG_{a,k} \stackrel{\rL}{\rightarrow} \mG_{a,k} \rightarrow 0.
\end{equation}
The push-forward of this extension by the character $\psi^{-1}_0$ \eqref{not1} defines a $\omQ_\ell$-sheaf 
of rank one, $\cL_{\psi_0}$, on $\mG_{a,k}$. Following Deligne, if $f\colon X\rightarrow \mG_{a,k}$
is a morphism of schemes, we put $\cL_{\psi_0}(f)=f^*\cL_{\psi_0}$.

\subsection{}\label{dtf1}
Let $\rA=\Spec(k[x])$ and $\crA=\Spec(k[\cx])$ be two affine lines over $k$ (equipped with 
coordinates $x$ and $\cx$) which are in duality via the pairing $\rA\times_k\crA\rightarrow \mG_{a,k}$
defined by $(x,\cx)\mapsto u=x\cx$.
We denote by $\rP$ and $\crP$ the projective lines over $k$, completions of $\rA$ and $\crA$
respectively, by $\infty\in \rP(k)$ and $\cinfty\in \crP(k)$ the points at infinity, 
by $j\colon \rA\rightarrow \rP$ and $\cj\colon \crA\rightarrow \crP$ the canonical 
injections and by $\pr$ and $\cpr$ the canonical projections of $\rA\times_k\crA$.
We have the $\omQ_\ell$-sheaf $\cL_{\psi_0}(x\cx)$ on $\rA\times_k\crA$;
we put $\ocL_{\psi_0}(x\cx)=(j\times \cj)_!\cL_{\psi_0}(x\cx)$ on $\rP\times_k\crP$.
For a complex $K$ of $\bD^b_c(\rA,\omQ_\ell)$,  the {\em Fourier transform} of $K$ 
is the complex $\fF_{\psi_0}(K)$ of $\bD^b_c(\crA,\omQ_\ell)$ defined by 
\begin{equation}
\fF_{\psi_0}(K)=\rR\cpr_!(\pr^*K\otimes \cL_{\psi_0}(x\cx)).
\end{equation}

In the sequel, we will omit the subscript $\psi_0$ from the notation $\fF_{\psi_0}$
and $\cL_{\psi_0}$ when there is no risque of confusion.  

\subsection{}\label{kummer}
The Kummer covering of order $2$ is the exact sequence 
\begin{equation}\label{kummera}
1\rightarrow \mu_2(k)\rightarrow \mG_{m,k}\stackrel{[2]}{\rightarrow}  \mG_{m,k}\rightarrow 1,
\end{equation}
where $[2]$ is the square power map. We denote by $\cK$ the $\omQ_\ell$-sheaf of rank $1$ on $\mG_{m,k}$ obtained by 
push-forward of this extension by the unique non-trivial character $\mu_2(k)\rightarrow \omQ_\ell^\times$. 
For a morphism $f\colon X\rightarrow \mG_{m,k}$, we put $\cK(f)=f^*\cK$. 

Consider the open subschemes $U=\rA-\{0\}$ of $\rA$ and $\cU=\crA-\{\czero\}$ of $\crA$, equipped 
with the isomorphisms $x\colon U\rightarrow \mG_{m,k}$ and $\cx\colon \cU\rightarrow \mG_{m,k}$. 
Let $\cQ$ be the $\omQ_\ell$-sheaf on $\Spec(k)$ defined by 
\begin{equation}\label{kummerb}
\cQ=\rH^1_c(U_{\ok},\cK(x)\otimes\cL_{\psi_0}(x)).
\end{equation}
Then $\cQ$ has rank $1$ and the $\rH^i_c(U_{\ok},\cK(x)\otimes\cL_{\psi_0}(x))$, for $i\not=1$, 
vanish. Moreover, we have canonical isomorphisms
\begin{eqnarray}
\fF(j_*\cK(x))[1]&\simeq &\cj_*\cK(\cx)\otimes \cQ, \label{kummerc}\\
\cj^*\rR \cpr_{!}(\cL_{\psi_0}(x^2\cx))[1]&\simeq &\cK(\cx)\otimes \cQ.\label{kummerd}
\end{eqnarray}
Indeed, the first assertion and the isomorphism \eqref{kummerc} are proved in (\cite{laumon} 1.4.3.1).
Consider the morphism $\pi\colon \rA\rightarrow \rA$ defined by $\pi(x)=x^2$. 
By the projection formula, we have a canonical isomorphism
\[
\rR \cpr_{!}(\cL_{\psi_0}(x^2\cx))\simeq \rR\cpr_{!}((\pi\times 1)^*\cL_{\psi_0}(x\cx))\simeq 
\rR\cpr_{!}((\pi\times 1)_*(\omQ_\ell)\otimes \cL_{\psi_0}(x\cx)).
\]
Since we have $\pi_*(\omQ_\ell)\simeq j_*\cK\oplus \omQ_\ell$, the isomorphism \eqref{kummerd} follows
from \eqref{kummerc} and (\cite{laumon} 1.2.2.2).

\subsection{}\label{results2}
Let $z$ be a closed point of $\rP$, $\cz$ be a closed point of $\crP$, 
$T$ and $\cT$ be the henselizations of $\rP$ and $\crP$ 
at $z$ and $\cz$ respectively, $\tau$ and $\ctau$ be the generic points of $T$ and $\cT$ respectively, 
$h\colon T\rightarrow \rP$ and $\ch\colon \cT\rightarrow \crP$ be the canonical morphisms,
$\pr$ and $\cpr$ be the canonical projections of $T\times_{k}\cT$. 
We denote also by $x$ and $\cx$ the pull-backs of the coordinates $x$ and $\cx$ of $\rA$ and $\crA$
over $\tau$ and $\ctau$ respectively, 
and (abusively) by $\ocL_{\psi_0}(x\cx)$ the sheaf $(h\times \ch)^*\ocL_{\psi_0}(x\cx)$ over $T\times_k\cT$. 
Let $\cF$ be a $\omQ_\ell$-sheaf over $\tau$, $\cF_!$ be its extension by zero to $T$. 
By (\cite{laumon} 2.3.2.1 and 2.3.3.1), the complex of vanishing cycles 
$\Phi(\pr^*(\cF_!)\otimes \ocL_{\psi_0}(x\cx))$
in $\bD_c^b(T\times_k\ctau,\omQ_\ell)$, relatively to the projection 
$\cpr\colon T\times_k\cT\rightarrow \cT$, 
is supported on $z\times_k\ctau$ and has cohomology only in degree $1$. 
Following Laumon, we call {\em local Fourier transform} of $\cF$ at $(z,\cz)$, 
and denoted by $\fF^{(z,\cz)}(\cF)$, the $\omQ_\ell$-sheaf over $z\times_k\ctau$ defined by
\begin{equation}
\fF^{(z,\cz)}(\cF)=(i_z\times 1)^*(\Phi^1(\pr^*(\cF_!)\otimes \ocL_{\psi_0}(x\cx))).
\end{equation} 
In fact, $\fF^{(z,\cz)}(\cF)$ vanishes if $(z,\cz)\in \rA\times\crA$ (\cite{laumon} 2.3.2.1 and 2.3.3.1).
Observe that if $z$ or $\cz$ is $k$-rational (which is the case if $(z,\cz)\notin \rA\times\crA$), 
then $z\times_k\cT$ is connected; more precisely, if $z$ is 
$k$-rational, then $z\times_k\cT=\cT$, and if $\cz$ is $k$-rational, then $z\times_k\cT$ is a finite
étale covering of $\cT$. 

\subsection{}\label{lftbc}
We keep the notation of \eqref{results2}, moreover,
let $k'$ be a finite extension of $k$, $u\colon k(z)\rightarrow k'$ and $\cu\colon k(\cz)\rightarrow k'$
be two $k$-homomorphisms, where $k(z)$ and $k(\cz)$ are the residue fields 
of $z$ and $\cz$ respectively. The pairs $(z,u)$ and $(\cz,\cu)$
define rational points $z'\in \rP_{k'}(k')$ and $\cz'\in \crP_{k'}(k')$. 
Let $T'$ and $\cT'$ be the henselizations of $\rP_{k'}$ and $\crP_{k'}$ at $z'$ and $\cz'$ respectively, 
$\tau'$ and $\ctau'$ be the generic points of $T'$ and $\cT'$ respectively, 
$f\colon T'\rightarrow T$ and $\cf\colon \cT'\rightarrow \cT$ be the canonical morphisms.
The canonical morphism $\cT'\rightarrow k'\otimes_k\cT$
induces a morphism $\tf\colon \cT'\rightarrow z\times_k\cT$. 
For any $\omQ_\ell$-sheaf  $\cF'$ over $\tau'$, we can consider 
the local Fourier transform $\fF^{(z',\cz')}(\cF')$, of $\cF'$ at $(z',\cz')$,
which is $\omQ_\ell$-sheaf over $\ctau'=z'\times_{k'}\ctau'$.

\begin{prop}\label{bext}
We keep the notation of \eqref{lftbc}.

{\rm (i)}\ For any $\omQ_\ell$-sheaf $\cF$ over $\tau$,
we have a canonical functorial isomorphism over $\ctau'$~:
\begin{equation}\label{bexta}
\tf^*(\fF^{(z,\cz)}(\cF))\simeq \fF^{(z',\cz')}(f^*\cF).
\end{equation}

{\rm (ii)}\ Assume that $\cz$ is $k$-rational. Then for any $\omQ_\ell$-sheaf $\cF'$ over $\tau'$,
we have a canonical functorial isomorphism over $z\times_k\ctau$~:
\begin{equation}\label{bext1}
\tf_*(\fF^{(z',\cz')}(\cF'))\simeq \fF^{(z,\cz)}(f_*\cF').
\end{equation}
\end{prop}

Consider the following commutative diagram with cartesian squares
\[
\xymatrix{
{T'\times_{k'}\cT'}\ar[r]^\gamma\ar[rdd]&{T'\times_{k}\cT'}\ar[r]\ar[dd]&
{T'\times_k\cT}\ar[r]^-(0.5){\pr_1}\ar[d]&{T'}\ar[d]^{f}\\
&&{T\times_k\cT}\ar[r]^-(0.5){\pr}\ar[d]^-(0.5){\cpr}&{T}\\
&{\cT'}\ar[r]^{\cf}&{\cT}&}
\]
Observe that $\gamma$ is an open and a closed immersion, and hence is étale. 
Let $\cG$ be a $\omQ_\ell$-sheaf over $T'\times_k\cT$, $\cG'$ be its pull-back over $T'\times_{k'}\cT'$. 
We consider the complexes of vanishing cycles $\Phi(\cG)$ in $\bD^b_c(T'\times_k\ctau,\omQ_\ell)$
and $\Phi(\cG')$ in $\bD^b_c(T'\times_{k'}\ctau',\omQ_\ell)$, relatively to the second projections
$T'\times_k\cT\rightarrow \cT$ and $T'\times_{k'}\cT'\rightarrow \cT'$ respectively.
We denote by 
\begin{equation}
p\colon T'\times_{k'}\ctau'\rightarrow T'\times_k\ctau
\end{equation}
the canonical morphism.
By (\cite{sga4.5} [Th. finitude] 3.7) and (\cite{sga7-2} XIII 2.1.7.2), 
and since $\gamma$ is étale, we have a canonical isomorphism
\begin{equation}\label{bexta3}
p^*(\Phi(\cG))\stackrel{\sim}{\rightarrow} \Phi(\cG').
\end{equation}

(i) We consider the sheaf $\cH=\pr^*(\cF_!)\otimes\ocL_{\psi_0}(x\cx)$ over $T\times_k\cT$,
and its complex of vanishing cycles $\Phi(\cH)$ 
in $\bD^b_c(T\times_k\ctau,\omQ_\ell)$, relatively to the projection $\cpr$. 
We take for $\cG$ the inverse image of $\cH$ 
over $T'\times_k\cT$. By (\cite{sga7-2} XIII 2.1.7.2), since $f$ is étale, we have a canonical isomorphism 
\begin{equation}\label{bexta4}
(f\times 1)^*(\Phi(\cH))\stackrel{\sim}{\rightarrow}
\Phi(\cG).
\end{equation}
Then the proposition follows from \eqref{bexta3} and \eqref{bexta4}.

(ii) We consider the sheaf $\cH=\pr^*(f_*\cF'_!)\otimes\ocL_{\psi_0}(x\cx)$ over $T\times_k\cT$,
and its complex of vanishing cycles $\Phi(\cH)$ 
in $\bD^b_c(T\times_k\ctau,\omQ_\ell)$, relatively to the projection $\cpr$. 
We take $\cG=\pr_1^*(\cF'_!)\otimes\ocL_{\psi_0}(x\cx)$ over $T'\times_k\cT$, 
where we denoted (abusively) by $\ocL_{\psi_0}(x\cx)$ the pull-back of the sheaf $\ocL_{\psi_0}(x\cx)$
over $T'\times_k\cT$. 
By (\cite{sga7-2} XIII 2.1.7.1) and the projection formula, 
since $f$ is finite, we have a canonical isomorphism 
\begin{equation}\label{bext4}
(f\times 1)_*(\Phi(\cG))\stackrel{\sim}{\rightarrow}
\Phi(\cH).
\end{equation}
On the other hand, the canonical morphism $\cT'\rightarrow k'\otimes_k\cT$ is an isomorphism
by assumption. Therefore, $p$ is an isomorphism, and we deduce from \eqref{bexta3} a functorial isomorphism 
\begin{equation}\label{bext3}
\Phi(\cG)\stackrel{\sim}{\rightarrow} p_*(\Phi(\cG')).
\end{equation}
The proposition follows from \eqref{bext4} and \eqref{bext3}.

\begin{teo}\label{results3}
Let $S$ be the spectrum of a henselian discrete valuation ring of equal characteristic $p$,
with perfect residue field, $s$ (resp. $\eta$) be the closed (resp. generic) point of $S$,
$\cG$ be a $\omQ_\ell$-sheaf of rank $1$ over $\eta$,
$v\colon S\rightarrow \rA$ and $\cv\colon S\rightarrow \crP$ be two non-constant morphisms
(with the notation of \ref{dtf1}). 
We put $z=v(s)$, $\cz=\cv(s)$, $b$ and $c$ the functions on $\eta$ deduced by pull-back 
from the coordinates $x$ and $\cx$ of $\rA$ and $\crA$ respectively. 
We take again the notation of \eqref{results2}
relatively to $z$ and $\cz$, and denote by $f\colon S\rightarrow T$ and 
$\cf\colon S\rightarrow \cT$ the morphisms induced by $v$ and $\cv$ respectively,
by $q\colon T\rightarrow z$ the canonical morphism, and by 
$\tf\colon S\rightarrow z\times_k\cT$ the morphism $(q\circ f, \cf)$.
We assume that $(\cG,b,c)$ is a {\em Legendre triple} at $s$ \eqref{legendre1},
and $f$ and $\cf$ are finite and étale at $\eta$. Then $\cz=\cinfty$ and we have a canonical isomorphism
\begin{equation}\label{results3a}
\fF^{(z,\cz)}(f_*\cG)\stackrel{\sim}{\rightarrow} \tf_*\left(\cG\otimes \cL_{\psi_0}(bc)\otimes 
\cK\left(-\frac{1}{2}\frac{dc}{db}\right)\otimes \cQ\right).
\end{equation}
\end{teo}

The following diagram summarizes the geometric picture of theorem \ref{results3}~:
\begin{equation}
\xymatrix{
{S\times_kS}\ar[rr]\ar[dd]&&S\ar[d]_f\ar@/^1pc/[rdd]^v&\\
&{T\times_k\cT}\ar[r]\ar[d]\ar[rd]&{T}\ar[rd]^-(0.4)h&\\
S\ar@/_1pc/[rrd]_{\cv}\ar[r]^{\cf}&{\cT}\ar[rd]_-(0.4){\ch}&{\rP\times_k\crP}\ar[r]\ar[d]&{\rP}\\
&&{\crP}&}
\end{equation}

The proof of theorem \ref{results3} is given in §\ref{pfgm}. For $(z,\cz)=(\infty,\czero)$
or $(\infty,\cinfty)$, the result is also valid under an extra condition \eqref{results3c}.

\begin{rema}\label{results3b}
We keep the notation of \eqref{results3}. It is clear that $\Omega^1_{T/k}$
is a free $\co_T$-module of rank~$1$. So the $\co_S$-module $\Omega^1_{S/k}$ is of finite
type, and hence free of rank $1$ because its completion along the closed point of $S$ is free of rank $1$.
In particular, $\frac{dc}{db}$ is a well defined function over $\eta$. 
\end{rema}

\begin{teo}\label{results3c}
Let $S$ be the spectrum of a henselian discrete valuation ring of equal characteristic $p$,
with perfect residue field, $s$ (resp. $\eta$) be the closed (resp. generic) point of $S$,
$\cG$ be a $\omQ_\ell$-sheaf of rank $1$ over $\eta$,
$v\colon S\rightarrow \rP$ and $\cv\colon S\rightarrow \crP$ be two non-constant morphisms
(with the notation of \ref{dtf1}). 
We put $z=v(s)$, $\cz=\cv(s)$, $b$ and $c$ the functions on $\eta$ deduced by pull-back 
from the coordinates $x$ and $\cx$ of $\rA$ and $\crA$ respectively. 
We take again the notation of \eqref{results2}
relatively to $z$ and $\cz$, and denote by $f\colon S\rightarrow T$ and 
$\cf\colon S\rightarrow \cT$ the morphisms induced by $v$ and $\cv$ respectively.
We assume that the following conditions are satisfied~:

{\rm (i)}\ $(\cG,b,c)$ is a {\em Legendre triple} at $s$ \eqref{legendre1}; 

{\rm (ii)}\ $f$ and $\cf$ are finite and étale at $\eta$; 

{\rm (iii)}\ $z=\infty$.

\noindent Moreover, we assume that one of the following conditions is satisfied~:

{\rm (iv)}\ $\cz=\czero$ and all the slopes of $f_*(\cG)$ are $<1$;

{\rm (iv')}\ $\cz=\cinfty$ and all the slopes of $f_*(\cG)$ are $>1$.

Then we have a canonical isomorphism
\begin{equation}\label{results3cc}
\fF^{(z,\cz)}(f_*\cG)\stackrel{\sim}{\rightarrow} \cf_*\left(\cG\otimes \cL_{\psi_0}(bc)\otimes 
\cK\left(-\frac{1}{2}\frac{dc}{db}\right)\otimes \cQ\right).
\end{equation}
\end{teo}

The proof of theorem \ref{results3c} is given in §\ref{pfgm}. 

\begin{rema}\label{results3d}
We keep the notation of \eqref{results3c}, and put $\sw(f_*\cG)$ and $\rk(f_*\cG)$
the Swan conductor and the rank of $f_*(\cG)$. If conditions (i), (ii) and (iii) are satisfied, 
if $f$ is tamely ramified and if $\sw(f_*\cG)\not=\rk(f_*\cG)$, 
then one of the conditions (iv) or (iv') is satisfied (cf. \eqref{psnum1a} and \cite{katz} 1.13).  
\end{rema}

\section{Nearby cycles and blow-up of the diagonal}\label{nbc}

\subsection{}\label{nbc1}
We keep the notation of \eqref{dtf1}. 
Let $X=\Spec(B)$ be a smooth connected affine curve over $k$, $s\in X(k)$, $U$ be the open 
subscheme $X-\{s\}$ of $X$,
$t\in B$ be a local parameter at $s$ which is invertible on $U$. 
Let $g\colon X\rightarrow \rP$ and $\cg\colon X\rightarrow \crP$ be two 
non-constant $k$-morphisms, $z=g(s)$, $\cz=\cg(s)$. 
We assume that $g$ and $\cg$ are generically étale and that $g(U)\subset \rA$ and $\cg(U)\subset \crA$.
The coordinates $x$ and $\cx$ of $A$ and $\crA$ define two sections in $\Gamma(U,\co_X)=B_t$,
that we denote also by $x$ and $\cx$.

We take again the notation of \eqref{results2} relatively to the points $z$ and $\cz$.
Let $S=\Spec(B^\tth)$ be the henselization of $X$ at $s$, 
$\hbar\colon S\rightarrow X$ be the canonical morphism,
$f\colon S\rightarrow T$ and $\cf\colon S\rightarrow \cT$ 
be the morphisms induced by $g$ and $\cg$ respectively.
We denote (abusively) by $t$ the uniformizer of $B^\tth$ image of $t\in B$
and by $s$ the closed point of $S$.
Let $\eta$ be the generic point of $S$, $\oeta$ (resp. $\os$) be a geometric point of $S$ 
above $\eta$ (resp. $s$). We denote by $b$ and $c$ the images of respectively 
$x$ and $\cx$ by the canonical homomorphism $B_t\rightarrow B^\tth_t$. 
Recall \eqref{results3b} that $\frac{dc}{db}$  is a well defined function over $\eta$.

We denote by $\pr_1$ and $\pr_2$ the canonical projections of $X\times_kS$ or $X\times_k\cT$~:
\[
\xymatrix{
{X\times_kS}\ar[rr]\ar[dd]_{\pr_2}\ar[dr]_-(0.5){\pr_1}&&{X\times_k\cT}\ar[ld]^-(0.5){\pr_1}\ar[dd]^{\pr_2}\\
&X&\\
S\ar[rr]^{\cf}&&{\cT}}
\]
We consider the sheaf $\ocL_{\psi_0}(x \cx)$ over $\rP\times_k\crP$, and denote also by $\ocL_{\psi_0}(x \cx)$ 
its pull-back by $g\times \ch$ over $X\times_k\cT$, and by $\ocL_{\psi_0}(x c)$ its pull-back 
by $g\times (\ch \circ \cf)$ over $X\times_kS$. 
These notation are coherent with our conventions, and do not lead to any ambiguity.

\subsection{}\label{nbc2}
Let $\cG$ be a smooth $\omQ_\ell$-sheaf of rank $1$ over $U$, $\cG_!$ be the 
extension by $0$ of $\cG$ to $X$, $\cG_\eta=\hbar_U^*(\cG)$.  
The purpose of this section is to study the complex of vanishing cycles 
\[
\Phi(\pr_1^*(\cG_!)\otimes\ocL_{\psi_0}(x \cx))
\]
in $\bD_c^b(X\times_k\ctau,\omQ_\ell)$, relatively to the projection $\pr_2\colon X\times_k\cT\rightarrow \cT$.
By (\cite{sga4.5} [Th. finitude] 3.7), $(1\times \cf)^*(\Phi(\pr_1^*(\cG_!)\otimes\ocL_{\psi_0}(x \cx)))$ is 
canonically isomorphic to the complex of vanishing cycles 
\[
\Phi(\pr_1^*(\cG_!)\otimes\ocL_{\psi_0}(x c))
\]
in $\bD_c^b(X\times_k\eta,\omQ_\ell)$, relatively to the projection $\pr_2\colon X\times_kS\rightarrow S$.

\begin{prop}\label{nbc3}
Assume that $(\cG,x,\cx)$ is a Legendre triple at $s$ \eqref{legendre1}. 

{\rm (i)}\ The complex $\Phi(\pr_1^*(\cG_!)\otimes\ocL_{\psi_0}(x c))$ is supported at 
$s\times_k\eta$ and has cohomology only in degree~$1$,
and the complex $\Phi(\pr_1^*(\cG_!)\otimes\ocL_{\psi_0}(x \cx))$ is supported at 
$s\times_k\ctau$ and has cohomology only in degree~$1$.

{\rm (ii)}\ The sheaf $(i_s\times 1)^*(\Phi^1(\pr_1^*(\cG_!)\otimes\ocL_{\psi_0}(x c)))$ over $\eta$
has a direct factor canonically isomorphic~to 
\begin{equation}\label{nbc3a}
\cD=\cG_\eta\otimes \cL_{\psi_0}(bc)\otimes \cK\left(-\frac{1}{2}\frac{dc}{db}\right)\otimes \cQ,
\end{equation}
where the sheaves $\cK$ and $\cQ$ are defined in \eqref{kummer}. 

{\rm (iii)}\ The morphism 
\begin{equation}\label{nbc3.5a}
\cf_*(\cD)\rightarrow (i_s\times 1)^*(\Phi^1(\pr_1^*(\cG_!)\otimes\ocL_{\psi_0}(x \cx)))
\end{equation}
induced by the trace morphism $\cf_*\cf^*\rightarrow \id$, is injective.
\end{prop}

Proposition (i) is due to Laumon (\cite{laumon} 2.3.2.1 and 2.3.3.1). We prove the first
statement (which implies the second). 
By the $t$-exactness of the functor $\Phi$ (\cite{bbd} 4.4.2 and \cite{illusie} 4.2),
it is enough to prove that the complex $\Phi(\pr_1^*(\cG_!)\otimes\ocL_{\psi_0}(x c))$ 
is supported on $s\times_k\eta$. Let $\ord$ be the valuation of $B^\tth_t$ 
normalized by $\ord(t)=1$. If $\ord(c)\geq 0$, the sheaf $\pr_1^*(\cG_!)\otimes\ocL_{\psi_0}(x c)$ 
is smooth over $U\times_kS$, and the assertion follows from (\cite{sga7-2} XIII 2.1.5). 
If $\ord(c)<0$, the assertion is a consequence of (\cite{laumon} 1.3.1.2). 

Propositions \ref{nbc3}(ii)-(iii) will be proved in \ref{nbc9} and \ref{nbc10}.

\subsection{}\label{nbc4}
Let $X\times^{\log}_{k}S$ be the logarithmic product of $X$ and $S$ over $k$, 
{\em i.e.}, the open subscheme of the blow-up of $X\times_{k}S$ along the closed point $s\times s$, 
obtained by removing the strict transforms of the axes $X\times_{k}s$ and 
$s\times_{k}S$ (\cite{aml} 4.3). The parameter $t$ identifies $X\times^{\log}_{k}S$ 
with the affine scheme defined by the $k$-algebra 
\begin{equation}\label{nbc4a}
B\otimes_{k}B^\tth[w,w^{-1}]/(t\otimes 1-1\otimes t\cdot w).
\end{equation}
We denote also by $\pr_1$ and $\pr_2$ the canonical projections from $X\times^{\log}_{k}S$ to $X$ and $S$
respectively. Then $\pr_2$ is smooth, and we have a canonical isomorphism 
$\pr_2^{-1}(\eta) \simeq  U\times_k\eta$. 
The strict transform of the graph $S\rightarrow X\times_{k}S$ of $\hbar$ defines a closed embedding 
\begin{equation}\label{nbc4b}
\delta\colon S\rightarrow X\times^{\log}_{k}S,
\end{equation} 
whose ideal in the ring \eqref{nbc4a} is generated by $w-1$. 
We put $Y=\pr_2^{-1}(s)$ and $e=\delta(s)\in Y(k)$. 
Then $Y$ is canonically isomorphic to the multiplicative group 
$\mG_{m,k}=\Spec(k[w,w^{-1}])$ and $e$ is the neutral element $1$.

\subsection{} We denote by $\cH$ the sheaf over $U\times_k\eta$ defined by
\begin{equation}\label{nbc2a}
\cH=\cHom(\pr_2^*(\cG_\eta),\pr_1^*(\cG)),
\end{equation}
and consider the complex of nearby cycles 
\[
\Psi_\eta(\cH\otimes\cL_{\psi_0}(c(x-b)))
\]
in $\bD_c^b(Y\times_k\eta,\omQ_\ell)$, relatively to the projection $\pr_2\colon X\times_k^{\log}S\rightarrow S$.
The following proposition is a refinement of \ref{nbc3}(i)-(ii). It will not be used in this article, 
and will be proved in \ref{pfnbc3bis}. 

\begin{prop}\label{nbc3bis}
Assume that $(\cG,x,\cx)$ is a Legendre triple at $s$ \eqref{legendre1} and $\cz\in \{\czero,\cinfty\}$. 

{\rm (i)}\ The complex $\Psi_\eta(\cH\otimes\cL_{\psi_0}(c(x-b)))$ is supported on $\Sigma\times_k\eta$,
where $\Sigma$ is a finite subgroup-scheme of $Y=\mG_{m,k}$, and has cohomology only in degree $1$.

{\rm (ii)}\ The sheaf $(i_e\times 1)^*(\Psi_\eta^1(\cH\otimes\cL_{\psi_0}(c(x-b))))$ over $\eta$
has a direct factor canonically isomorphic to $\cK\left(-\frac{1}{2}\frac{dc}{db}\right)\otimes \cQ$, 
where the sheaves $\cK$ and $\cQ$ are defined in \eqref{kummer}.
\end{prop}

\begin{rema}
Condition \eqref{legendre1c} contained in the definition of a Legendre triple is not
necessary for proposition \ref{nbc3bis}(i).
\end{rema}

\subsection{}
We assume in the sequel of this section that $(\cG,x,\cx)$ is a Legendre triple at $s$. 
Let $R$ be the completion of $B^\tth$, $K$ be the fraction field of $R$.
We identify $R$ with the ring of power series $k[[t]]$.
We denote also by $b$ and $c$ the images of $b$ and $c$ in $K$.  
By definition \eqref{legendre1}, there exist $\cG_\rt$ and $\cG_\rw$ 
two smooth $\omQ_\ell$-sheaves of rank $1$ over $U$ and $a\in \rW_{m+1}(K)$ ($m\geq 0$), 
satisfying the following conditions~:

{\rm (i)}\ $\cG\simeq \cG_\rt\otimes \cG_\rw$; 

{\rm (ii)}\ $\cG_\rt$ is tamely ramified at $s$; 

{\rm (iii)}\ $\cG_\rw$ is trivialized by a cyclic extension of order $p^{m+1}$ of $U$.

{\rm (iv)}\ If we put $\chi=\delta_{m+1}(a)\in \rH^1(K,\mZ/p^{m+1}\mZ)$, then $\psi_m^{-1}\circ \chi$ is 
the character of the pull-back of $\cG_\rw$ to $\Spec(K)$ \eqref{not1}.  

{\rm (v)}\ $(a,b,c)$ is a Legendre triple. Let $(n,\nu(b),\nu(c))$ be its conductor \eqref{legendre}.  

\vspace{2mm}

We consider the sheaves $\cH_\rt$ and $\cH_\rw$ over $U\times_k\eta$ defined by 
\begin{eqnarray}
\cH_\rt&=&\cHom(\pr_2^*(\hbar_U^*(\cG_\rt)),\pr_1^*(\cG_\rt)),\label{nbd1}\\
\cH_\rw&=&\cHom(\pr_2^*(\hbar_U^*(\cG_\rw)),\pr_1^*(\cG_\rw)),\label{nbd2}
\end{eqnarray}
so we have $\cH\simeq \cH_\rt\otimes \cH_\rw$ \eqref{nbc2a}.
We take again the notation of §\ref{calculus}, and define $\alpha\in K$ by the equation $\rF^md(a)=\alpha dt$. 
So we have $a\in \fil_n\rW_{m+1}(K)$, $\ord(t\alpha)=-n$, $2\ord(\alpha+cb')\geq -n+\nu(c)$
and $2\nu(b)+p\nu(c)<(p-2)n$. In particular, we have $n-\nu(c)\geq 1$. 
We denote by $\gamma$ the non-zero element of $k$ defined by 
\begin{equation} \label{nbc5f}
\gamma=\frac{1}{2}\left(t^{n+1}\alpha\frac{t^{1-\nu(c)}c'}{c}\right) \mod t R.
\end{equation}

\begin{prop}\label{nbd3}
The sheaf $\cH_\rt$ over $U\times_k\eta$ extends to a smooth sheaf over $X\times_k^{\log}S$
whose pull-back by $\delta$ over $S$ is constant. 
\end{prop}

We denote by $\kappa_0$ the generic point of $Y$ \eqref{nbc4}, 
by $R_{L_0}$ the completion of the local ring of $X\times^{\log}_{k}S$ at $\kappa_0$
(which is a discrete valuation ring), by $L_0$ the fraction field of $R_{L_0}$, 
by $I_K^\rt$ (resp. $I_{L_0}^\rt$) the tame inertia group of $K$ (resp. $L_0$). 
The restriction of $\cH_\rt$ to $\Spec(L_0)$ is tamely ramified. 
Since the projections $\pr_1$ and $\pr_2$ of $X\times_k^{\log}S$ are smooth, they induce the 
same isomorphism $I_{L_0}^\rt\stackrel{\sim}{\rightarrow} I_{K}^\rt$. 
We deduce that the representation of $I_{L_0}^\rt$ defined by the sheaf $\cH_\rt$ 
is trivial, and hence the restriction of $\cH_\rt$ to $\Spec(L_0)$ is unramified.
The first assertion follows by the Zariski-Nagata's purity theorem (SGA 2 X 3.4). 
The second assertion is a consequence of the first one and the fact that the restriction 
of $\cH_{\rt}$ to $\delta(\eta)$ is trivial.

\subsection{}\label{nbc5}
Let $r$ be an integer $\geq 1$, $S_r$ be the closed subscheme of $S$ defined by $t^r$, 
$(X\times^{\log}_kS)_{[r]}$ be the blow-up of $X\times_k^{\log}S$ along $\delta(S_r)$,
$(X\times^{\log}_{k}S)_{(r)}$ be the dilatation of $X\times_k^{\log}S$ along $\delta$
of thickening $r$, that is, the open subscheme of $(X\times^{\log}_kS)_{[r]}$ obtained 
by removing the strict transform of $Y$, or equivalently, 
the maximal open subscheme of $(X\times^{\log}_kS)_{[r]}$ where the exceptional divisor is defined by $\pr_2^*(t^r)$
(\cite{aml} 3.1). We denote by $\Theta_r$ the exceptional divisor on $(X\times^{\log}_{k}S)_{(r)}$, by 
\begin{equation}\label{nbc5a}
\delta_{(r)}\colon S\rightarrow (X\times^{\log}_{k}S)_{(r)}
\end{equation}
the unique lifting of $\delta$ ({\em i.e.} the strict transform of $\delta$),
and abusively by $\pr_1$ and $\pr_2$ the projections from $(X\times^{\log}_{k}S)_{(r)}$ 
to $X$ and $S$ respectively. Then $\pr_2$ is smooth; the commutative diagram
\begin{equation}\label{nbv5b}
\xymatrix{
{\Theta_r}\ar[r]\ar[d]&{(X\times^{\log}_{k}S)_{(r)}}\ar[d]^{\pr_2}&{U\times_{k}\eta}\ar[l]\ar[d]\\
{S_r}\ar[r]&{S}&{\eta}\ar[l]}
\end{equation}
has Cartisian squares; and $\Theta_r$ is canonically isomorphic to the vector bundle 
$\bV(\Omega^1_{X/k}(\log s)\otimes_X\co_{S_r}(S_r))$ over $S_r$ (\cite{aml} 4.6).

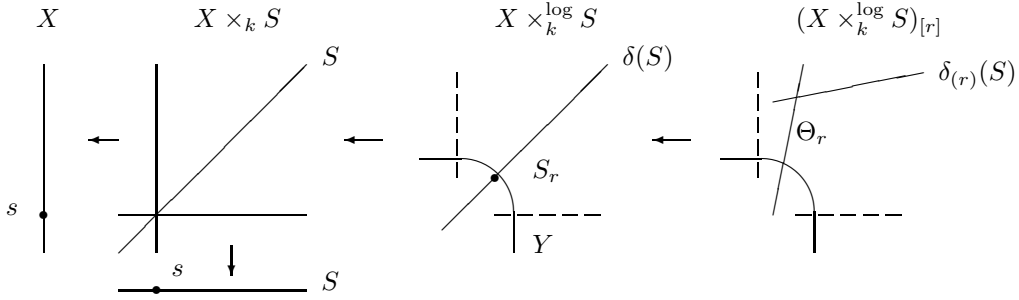
\begin{center}
\setlength{\unitlength}{1mm}
\begin{picture}(140,50)
\put(5,5){\line(0,1){25}}
\put(5,10){\circle*{1}}
\put(4,35){$X$}
\put(0,10){$s$}
\put(15,20){\vector(-1,0){4}}
\put(15,0){\line(1,0){25}}
\put(20,0){\circle*{1}}
\put(30,6){\vector(0,-1){4}}
\put(22,2){$s$}
\put(42,0){$S$}
\put(15,5){\line(1,1){25}}
\put(20,5){\line(0,1){25}}
\put(15,10){\line(1,0){25}}
\put(42,30){$S$}
\put(25,35){$X\times_k S$}
\put(50,20){\vector(-1,0){5}}
\put(58,8){\line(1,1){22}}
\put(55,5){\oval(25,25)[tr]}
\put(65,15){\circle*{1}}
\multiput(60,15)(0,3){5}{\line(0,1){2}}
\multiput(65,10)(3,0){5}{\line(1,0){2}}
\put(70,5){$Y$}
\put(70,15){$S_r$}
\put(82,30){$\delta(S)$}
\put(65,35){$X\times_k^{\log} S$}
\put(91,20){\vector(-1,0){5}}
\put(102,10){\line(1,5){4}}
\put(102,25){\line(5,1){20}}
\put(95,5){\oval(25,25)[tr]}
\multiput(100,15)(0,3){5}{\line(0,1){2}}
\multiput(105,10)(3,0){5}{\line(1,0){2}}
\put(124,28){$\delta_{(r)}(S)$}
\put(105,20){$\Theta_r$}
\put(105,35){$(X\times_k^{\log} S)_{[r]}$}
\end{picture}
\end{center}

\vspace{5mm}

It follows from \eqref{nbc4a} that $(X\times_k^{\log}S)_{(r)}$ is the affine scheme of ring 
\begin{equation}\label{nbv5e}
B\otimes_kB^\tth[\theta,(1+1\otimes t^r \cdot \theta)^{-1}]/(t\otimes 1-1\otimes t(1+1\otimes t^r \cdot \theta)),
\end{equation}
$\Theta_r\otimes_{R}k$ is the affine line $\rA=\Spec(k[\theta])$ over $k$ (with coordinate $\theta$), 
and $\delta_{(r)}$ is defined by the equation $\theta$.
Let $\kappa$ be the generic point of $\Theta_r$,  
$R_L$ be the completion of the local ring of $(X\times^{\log}_{k}S)_{(r)}$ at $\kappa$
(which is a discrete valuation ring), $L$ be the fraction field of $R_L$,
\begin{eqnarray}
u\colon K&\rightarrow L \label{nbv5c}\\
v\colon K&\rightarrow L \label{nbv5d}
\end{eqnarray}
be the homomorphisms induced respectively by $\pr_1$ and $\pr_2$. 
We consider $L$ as an extension of $K$ by $v$.
By \eqref{nbv5e}, we can identify $R_L$ with the ring $k(\theta)[[t]]$. Then the $k$-homomorphisms
$u$ and $v$ are defined by $u(t)=t(1+t^r\theta)$ and $v(t)=t$.

\begin{prop}\label{nbc6} 
Assume that $n-\nu(c)=2r$ is even.
Let $\Gamma_\eta$ be the closed subscheme of $U\times_k\eta$ inverse image 
by the morphism $\cg\times_k\eta$ of the section $\eta\rightarrow \crA\times_k\eta$ defined by $c$
(which is also the closed subscheme of $U\times_k\eta$ defined by the equation $\cx-c$),
$\Gamma_{(r)}$ be the schematic closure of $\Gamma_\eta$ in $(X\times_k^{\log}S)_{(r)}$. Then~:

{\rm (i)}\ The scheme $\Gamma_{(r)}$ is quasi-finite over $S$,
and $\delta_{(r)}(S)$ is the finite part of $\Gamma_{(r)}$.

{\rm (ii)}\ The sheaf $\cH\otimes\cL_{\psi_0}(c(x-b))$ over $U\times_k\eta$ extends to a smooth sheaf 
over $(X\times^{\log}_kS)_{(r)}$ whose restriction to $\Theta_r\otimes_Rk=\rA$ 
is canonically isomorphic to the sheaf $\cL_{\psi_0}(\gamma\theta^2)$.
\end{prop}

(i) Recall first that any quasi-finite separated scheme $Z$ over $S$ can de decomposed canonically 
into a sum $Z^{f}\amalg Z^g$, where $Z^{f}$ is finite over $S$ (called the finite part of $Z$) 
and the special fiber of $Z^g$ is empty. 
It follows from \eqref{key4a} and the inequality $\nu(c)<r(p-1)$ that we have 
\begin{equation}\label{nbc6b}
\frac{u(c)-c}{t^{r+\nu(c)}c} \equiv \frac{t^{1-\nu(c)}c'}{c}  \theta \mod tR_L.
\end{equation}
We deduce that the function $t^{-r-\nu(c)}(\cx/c-1)$ on $U\times_k\eta$ extends to a regular 
function on a neighborhood of $\kappa$ in $(X\times_k^{\log}S)_{(r)}$. 
Since $(X\times_k^{\log}S)_{(r)}$ is a smooth curve over $S$ with an integral special fiber,
$t^{-r-\nu(c)}(\cx/c-1)$ extends to a regular function on $(X\times_k^{\log}S)_{(r)}$.
The latter belongs to the ideal of $\Gamma_{(r)}$; so it defines a closed subscheme $\Gamma'_{(r)}$ of  
$(X\times_k^{\log}S)_{(r)}$ containing $\Gamma_{(r)}$. 
Moreover, it follows from \eqref{nbc6b} that the special fiber of $\Gamma'_{(r)}$ is 
the origin of $\Theta_r\otimes_Rk=\Spec(k[\theta])$. Therefore, the three closed subschemes 
$\delta_{(r)}(S)\subset \Gamma_{(r)}\subset \Gamma'_{(r)}$ of $(X\times_k^{\log}S)_{(r)}$
have the same special fiber, which implies the proposition. 

(ii) It follows from \ref{nbd3} that $\cH_\rt$ extends to a smooth sheaf on $(X\times_k^{\log}S)_{(r)}$ 
whose restriction to $\Theta_r\otimes_Rk$ is constant. 
The pull-back of $\cH_\rw\otimes\cL_{\psi_0}(c(x-b))$ to $\Spec(L)$ 
corresponds to the Witt vector 
\[
u(a)-a+\rV^m(c(u(b)-b))\in \rW_{m+1}(L).
\]
By \ref{key7}, the latter belongs to $\rW_{m+1}(R_L)$, and its residue class modulo $t$ 
is equal to $\rV^m(\gamma \theta^2)\in \rW_{m+1}(k(\theta))$. Therefore, the pull-back of
$\cH_\rw\otimes\cL_{\psi_0}(c(x-b))$ to $\Spec(L)$ is canonically isomorphic to $\cL_{\psi_0}(\gamma \theta^2)$. 
We deduce by the Zariski-Nagata's purity theorem (SGA 2 X 3.4)
that $\cH_\rw\otimes\cL_{\psi_0}(c(x-b))$ extends to a smooth sheaf on $(X\times_k^{\log}S)_{(r)}$ 
whose restriction to $\Theta_r\otimes_Rk=\rA$ is canonically isomorphic to $\cL_{\psi_0}(\gamma \theta^2)$. 

\begin{rema}\label{nbc6e}
We keep the notation and assumptions of \eqref{nbc6}. Recall that we have a commutative diagram 
\[
\xymatrix{
{S}\ar[r]^{\hbar}\ar[d]_{\cf}&{X}\ar[d]^{\cg}\\
{\cT}\ar[r]^{\ch}&{\crP}}
\]
We consider $\eta$ as a scheme over $\crA$ by the composed morphism $\cg\circ \hbar=\ch\circ \cf$. 
Then the closed subscheme $\Gamma_\eta$ of $U\times_k\eta$ is canonically 
isomorphic to $U\times_{\crA}\eta$. The morphism $\delta$ \eqref{nbc4b} induces a closed emdedding 
$\eta\rightarrow U\times_k\eta$, which determines a connected component of $U\times_{\crA}\eta$.
Proposition \ref{nbc6}(i) says that only this connected component extends to a closed subscheme 
of $(X\times_k^{\log}S)_{(r)}$ which is finite and flat over $S$, namely $\delta_{(r)}(S)$; 
the other connected components of $U\times_{\crA}\eta$ are closed in $(X\times_k^{\log}S)_{(r)}$. 
\end{rema}

\begin{lem}\label{nbc8}
Let $\cK$ and $\cQ$ de the sheaves defined in \eqref{kummer}, 
$\gamma$ be the element of $k$ defined in \eqref{nbc5f}, 
$\cQ_\gamma$ be the $\omQ_\ell$-sheaf over $\Spec(k)$ defined~by 
\begin{equation}\label{nbc8a}
\cQ_\gamma=\rH^1_c(\rA_{\ok},\cL_{\psi_0}(\gamma \theta^2)).
\end{equation}

{\rm (i)}\ The sheaf $\cQ_\gamma$ has rank $1$, the $\rH^i_c(\rA_{\ok},\cL_{\psi_0}(\gamma\theta^2))$, 
for $i\not=1$, vanish, and the canonical morphism 
\[
\rH^1_c(\rA_{\ok},\cL_{\psi_0}(\gamma \theta^2))\rightarrow \rH^1(\rA_{\ok},\cL_{\psi_0}(\gamma \theta^2))
\] 
is an isomorphism. 

{\rm (ii)}\ If $n-\nu(c)=2r$ is even, the sheaf $\cK(-\frac{1}{2}\frac{dc}{db})\otimes \cQ$
over $\eta$ is unramified and isomorphic to the geometrically constant sheaf $\cQ_\gamma$.

{\rm (iii)}\ If $n-\nu(c)=2r+1$ is odd, the sheaf $\cK(-\frac{1}{2}\frac{dc}{db})\otimes \cQ$ over $\eta$
is tamely ramified, and its restriction to the quadratic extension 
$\ueta=\eta[\ut]/(\ut^2-t)$ is unramified and isomorphic to the geometrically constant sheaf $\cQ_\gamma$. 
\end{lem}

Observe first that it is enough to prove the lemma after replacing $\eta$ by $\Spec(K)$. 

(i) Since the sheaf $\cL_{\psi_0}(\gamma\theta^2)$ is smooth on $\rA$ and its Swan conductor at 
$\infty$ is $2$, the assertion follows from the Grothendieck-Ogg-Shafarevich formula
and (\cite{sga4.5} [Sommes trig.] 1.19 and 1.19.1). 

(ii) Condition $2\ord(\alpha+cb')\geq -n+\nu(c)$ implies that 
\[
\frac{1}{2}\left(t^{n+1}\alpha\frac{t^{1-\nu(c)}c'}{c}\right)\equiv-\frac{1}{2}(t^{2+2r}b'c') \mod t^{1+r+\nu(c)}R;
\]
so $-\frac{c'}{2b'}\gamma^{-1}$ is a square in $K$. Therefore, the sheaf $\cK(-\frac{c'}{2b'})$
over $\Spec(K)$ is unramified and isomorphic to the geometrically constant sheaf $\cK(\gamma)$. 
The last assertion follows from \eqref{kummerd}. 

(iii) Condition $2\ord(\alpha+cb')\geq -n+\nu(c)$ implies that 
\[
\frac{1}{2}\left(t^{n+1}\alpha\frac{t^{1-\nu(c)}c'}{c}\right)\equiv-\frac{1}{2}(t^{3+2r}b'c') \mod t^{2+r+\nu(c)}R;
\]
so $-\frac{tc'}{2b'}\gamma^{-1}$ is a square in $K$. The proposition follows as in (ii).

\subsection{}\label{sigma}
Let $\sigma$ be a $k$-automorphism of $S$. We put
$\delta^{(\sigma)}\colon S\rightarrow X\times_k^{\log}S$ the strict transform of the graph of $\hbar\circ \sigma$. 
The automorphism $1\times \sigma$ 
of $X\times_kS$ lifts uniquely to $X\times_k^{\log}S$ and the following diagram
\[
\xymatrix{
S\ar[rr]^-(0.5){\delta^{(\sigma)}}\ar[d]_\sigma&&{X\times_k^{\log}S}\ar[d]^{1\times \sigma}\\
S\ar[rr]^-(0.5){\delta}&&{X\times_k^{\log}S}}
\]
is commutative. Hence, $1\times \sigma$ 
induces an isomorphism between the dilatations of $X\times_k^{\log}S$
along $\delta$ and $\delta^{(\sigma)}$ with the same thickening.

Let $\sigma$ be a $k$-automorphism of $X$ such that $\sigma(x)=x$. We denote also by $\sigma$
the $k$-automorphism of $S$ induced by $\sigma$.   
The automorphism $\sigma\times 1$ of $X\times_kS$ lifts uniquely to $X\times_k^{\log}S$ and the following diagram
\[
\xymatrix{
S\ar[r]^-(0.5){\delta}\ar[rd]_-(0.5){\delta^{(\sigma)}}&{X\times_k^{\log}S}\ar[d]^{\sigma\times 1}\\
&{X\times_k^{\log}S}}
\]
is commutative. Hence, $\sigma\times 1$ 
induces an isomorphism between the dilatations of $X\times_k^{\log}S$
along $\delta$ and $\delta^{(\sigma)}$ with the same thickening.

\subsection{}\label{nbc9}
We can now prove proposition \ref{nbc3}(ii).  Observe first that we have
\[
(i_s\times 1)^*(\Psi(\pr_1^*(\cG_!)\otimes\ocL_{\psi_0}(x c)))=
(i_s\times 1)^*(\Phi(\pr_1^*(\cG_!)\otimes\ocL_{\psi_0}(x c))).
\]
Let $\ocH$ be the extension by $0$ of $\cH$ \eqref{nbc2a} to $X\times_k\eta$,
\[
\cM=\ocH\otimes\ocL_{\psi_0}(c(x-b)).
\]
We identify $\cM$ with $\pr_2^*(\cG_\eta^\vee\otimes \cL_{\psi_0}(-bc))\otimes \pr_1^*(\cG_!)\otimes\ocL_{\psi_0}(x c)$,
where $\cG_\eta^\vee$ is the dual $\omQ_\ell$-sheaf of $\cG_\eta$ over $\eta$. 
By (\cite{illusie} 4.7) (applied with $Y=S$),  we have a canonical isomorphism
\begin{equation}\label{nbc9a}
\Psi_\eta(\cM)\stackrel{\sim}{\rightarrow}
\pr_2^*(\cG_\eta^\vee\otimes \cL_{\psi_0}(-bc))\otimes \Psi_\eta(\pr_1^*(\cG_!)\otimes\ocL_{\psi_0}(x c))
\end{equation}
in $\bD^b_c(X\times_k\eta,\omQ_\ell)$. 
So we are reduced to proving that the sheaf $(i_s\times 1)^*(\Psi^1_\eta(\cM))$
over $\eta$ has a direct factor isomorphic to $\cK(-\frac{1}{2}\frac{db}{dc})\otimes\cQ$.

Assume first that $n-\nu(c)=2r$ is even. It follows from \ref{nbc6}(ii) and (\cite{sga7-2} 2.1.7.1) 
that we have canonical morphisms of representations of $\pi_1(\eta,\oeta)$ 
\begin{equation}\label{nbc9b}
\rH^1_c(\rA^1_{\ok},\cL_{\psi_0}(\gamma \theta^2))\stackrel{u}{\rightarrow}
\Psi^1_\eta(\cM)_{(\os,\oeta)}\stackrel{v}{\rightarrow}
\rH^1(\rA^1_{\ok},\cL_{\psi_0}(\gamma \theta^2)),
\end{equation}
where the source of $u$ and the target of $v$ 
are considered as unramified representations of $\pi_1(\eta,\oeta)$.
Moreover, $v\circ u$ is the canonical morphism; so it is an isomorphism by \ref{nbc8}(i).
Hence the required assertion follows by \ref{nbc8}(ii).  

Assume next that $n-\nu(c)$ is odd. 
We put $\uB=B[\ut]/\ut^2-t$, $\uX=\Spec(\uB)$, $\us\in \uX(k)$ the unique point of $\uX$ above $s$, 
and denote by an underline the objects deduced from objects over $X$ by the base change $\uX\rightarrow X$.
So $\uS$ is the spectrum of the henselization of the local ring of $\uX$ at $\us$ 
and $\ueta$ is the generic point of $\uS$. 
Let $\rho\colon \ueta\rightarrow \eta$ be the canonical morphism, 
$G$ be the Galois group of $\ueta$ over $\eta$,  that we identify with the group of $X$-automorphisms of $\uX$. 
We consider the pull-back $\ucM$ of $\cM$ over $\uX\times_k\ueta$, and the complex 
of nearby cycles $\Psi_{\ueta}(\ucM)$ in $\bD^b_c(\uX\times_k\ueta,\omQ_\ell)$,
relatively to the second projection $\uX\times_k\uS\rightarrow \uS$. 
It follows from (\cite{sga7-2} XIII 2.1.7.1), (\cite{sga4.5} [Th. finitude] 3.7) and 
the Hochschild-Serre spectral sequence applied to the following diagram
\[
\xymatrix{{\uX\times_k\uS}\ar[r]\ar[rd]&{X\times_k\uS}\ar[r]\ar[d]&{X\times_kS}\ar[d]\\
&{\uS}\ar[r]&S}
\]
that we have an isomorphism 
\begin{equation}\label{nbc9c}
(i_s\times 1)^*(\Psi_\eta(\cM))\simeq (\rho_*((i_{\us}\times 1)^*(\Psi_{\ueta}(\ucM))))^{G\times G},
\end{equation}
where the group $G\times G$ acts on $\Psi_{\ueta}(\ucM)$ via its action on $\uX\times_k\uS$. 

We put $\uK=K[\ut]/(\ut^2-t)$. The image of $(a,b,c)$ in $\rW_{m+1}(\uK)\times \uK\times \uK$ 
is a Legendre triple of conductor $(2n,2\nu(b),2\nu(c))$ \eqref{legendre}. 
Hence, we can apply \ref{nbc6} over $(\uX\times_k^{\log}\uS)_{(r)}$ with $r=n-\nu(c)$. 
We deduce, as in the even case, that we have canonical morphisms of sheaves over $\ueta$ 
\begin{equation}\label{nbc9d}
\cK(-\frac{1}{2}\frac{db}{dc})\otimes\cQ\stackrel{u}{\rightarrow} (i_{\us}\times 1)^*(\Psi^1_{\ueta}(\ucM)) 
\stackrel{v}{\rightarrow} \cK(-\frac{1}{2}\frac{db}{dc})\otimes\cQ
\end{equation}
such that $v\circ u$ is the identity. We let $G\times G$ acts on the sheaf 
$\cK(-\frac{1}{2}\frac{db}{dc})\otimes\cQ$ over $\ueta$ through the action of 
the second factor $G$ on $\ueta$. It follows from \ref{sigma}, \eqref{nbv5e} 
and \ref{nbc6}(ii) that $u$ and $v$ are $(G\times G)$-equivariant.
The required assertion follows by using \eqref{nbc9c}.

\subsection{}\label{equiv}
Let $\oS$ be the integral closure of $S$ in $\oeta$, $\os$ be its closed point. 
We denote by $X\times_k^{\log}\oS$ the base change of $X\times_k^{\log}S$ by the morphism 
$\oS\rightarrow S$, which is also the logarithmic product of $X$ and $\oS$ over $k$ (\cite{aml} 4.3),
by $\oh\colon \oS\rightarrow X$ the morphism induced by $\hbar\colon S\rightarrow X$ and by
\begin{equation}\label{equiva}
\odelta\colon \oS\rightarrow  X\times_k^{\log}\oS
\end{equation}
the base change of $\delta$ \eqref{nbc4b}, which is also the strict transform of the graph of $\oh$. 
For any $k$-automorphism $\sigma$ of $\oS$, we denote by $\odelta^{(\sigma)}\colon \oS\rightarrow  X\times_k^{\log}\oS$
the strict transform of the graph of $\oh\circ \sigma$.  
Let $r$ be an integer $\geq 1$, $\oS_r$ be the closed subscheme of $\oS$ defined by $t^r$.
We denote by $\cJ_{\sigma}$ the ideal of $\odelta^{(\sigma)}(\oS_r)$ in $X\times_k^{\log}\oS$,
by $(X\times^{\log}_k\oS)_{[r]}^{(\sigma)}$ the blow-up of $X\times_k^{\log}\oS$ along $\cJ_{\sigma}$,
and by $(X\times^{\log}_k\oS)_{(r)}^{(\sigma)}$ the dilatation of $X\times^{\log}_k\oS$ along $\odelta^{(\sigma)}$
of thickening $r$, {\em i.e.}, the open subscheme of $(X\times^{\log}_k\oS)_{[r]}^{(\sigma)}$ where the 
exceptional divisor is generated by $\pr_2^*(t^r)$.

We consider $\eta$ as an étale covering of $\ctau$ by $\cf$.   
Observe that for any $\sigma\in \pi_1(\ctau,\oeta)$, the morphism $\odelta^{(\sigma)}$
depends only on the class of $\sigma$ in $\pi_1(\eta,\oeta)\backslash\pi_1(\ctau,\oeta)$
(the quotient of $\pi_1(\ctau,\oeta)$ by the subgroup $\pi_1(\eta,\oeta)$ acting by translation on the left). 
Moreover, the natural action of $\pi_1(\ctau,\oeta)$ on $X\times_k\oS$ 
lifts to $X\times_k^{\log}\oS$, and for any $\sigma,\sigma'\in \pi_1(\ctau,\oeta)$,
we have a commutative diagram
\begin{equation}\label{equivb}
\xymatrix{
\oS\ar[rr]^-(0.5){\odelta^{(\sigma'\sigma)}}\ar[d]_{\sigma}&&{X\times_k^{\log}\oS}\ar[d]^{1\times \sigma}\\
\oS\ar[rr]^-(0.5){\odelta^{(\sigma')}}&&{X\times_k^{\log}\oS}}
\end{equation}
In particular, $1\times \sigma$ induces an isomorphism from $(X\times^{\log}_k\oS)_{[r]}^{(\sigma'\sigma)}$
to $(X\times^{\log}_k\oS)_{[r]}^{(\sigma')}$, that transforms $(X\times^{\log}_k\oS)_{(r)}^{(\sigma'\sigma)}$
into $(X\times^{\log}_k\oS)_{(r)}^{(\sigma')}$.

\begin{lem}\label{equiv1}
Let $r$ be an integer $\geq 1$, $\fX$ be the blow-up of $X\times_k^{\log}\oS$ along 
the ideal $\cJ=\prod_{\sigma\in \pi_1(\eta,\oeta)\backslash\pi_1(\ctau,\oeta)}\cJ_{\sigma}$
(the quotient of $\pi_1(\ctau,\oeta)$ by the subgroup $\pi_1(\eta,\oeta)$ acting by translation on the left).  

{\rm (i)}\ The action of $\pi_1(\ctau,\oeta)$ on $X\times_k^{\log}\oS$ lifts uniquely to $\fX$.
For every $\sigma\in \pi_1(\ctau,\oeta)$, there exists 
a unique morphism $\varphi_{\sigma}\colon \fX\rightarrow (X\times^{\log}_k\oS)_{[r]}^{(\sigma)}$ 
over $X\times_k^{\log}\oS$. 

{\rm (ii)}\ For every $\sigma\in \pi_1(\ctau,\oeta)$, 
$\varphi_\sigma$ induces an isomorphism above $(X\times^{\log}_k\oS)_{(r)}^{(\sigma)}$. 
We put $\fX_{(\sigma)}=\varphi_{\sigma}^{-1}((X\times^{\log}_k\oS)_{(r)}^{(\sigma)})$.

{\rm (iii)}\ For every $\sigma,\sigma'\in \pi_1(\ctau,\oeta)$, 
we have $\sigma(\fX_{(\sigma'\sigma)})=\fX_{(\sigma')}$. 

{\rm (iv)}\ $\fX_{\os}$ is connected. 

{\rm (v)}\ Assume that $n-\nu(c)=2r$. 
Then for every $\sigma,\sigma'\in \pi_1(\ctau,\oeta)$ such that $\sigma'\not\in \pi_1(\eta,\oeta)\sigma$,
we have $\fX_{(\sigma)}\cap \fX_{(\sigma')} \cap \fX_{\os}=\emptyset$.
\end{lem}

(i) Since the action of $\pi_1(\ctau,\oeta)$ on $X\times_k^{\log}\oS$ preserves the
ideal $\cJ$, it lifts to an action on $\fX$. 
Since the ideal $\cJ\co_\fX$ is invertible, each ideal $\cJ_{\sigma}\co_\fX$ is invertible
(see Bourbaki Alg. Comm. chap.~II §5.6 théo.~4). Hence, 
for every $\sigma\in \pi_1(\ctau,\oeta)$, the canonical morphism
$\fX\rightarrow  X\times_k^{\log}\oS$ lifts uniquely to a morphism 
$\varphi_\sigma\colon \fX\rightarrow (X\times^{\log}_k\oS)_{[r]}^{(\sigma)}$.

(ii) By the universal property of blow-ups, it is enough to prove that for any 
$\sigma,\sigma'\in \pi_1(\ctau,\oeta)$, the inverse image of the ideal $\cJ_{\sigma'}$
over $(X\times^{\log}_k\oS)_{(r)}^{(\sigma)}$ is invertible. 
We denote by $C$ the ring of the affine scheme $X\times_k^{\log}\oS$ and 
(abusively) by $t$ the function $\pr_2^*(t)\in C$. 
The embedding $\delta$ \eqref{nbc4b} is defined by the equation $w-1$ of the ring \eqref{nbc4a}. 
Let $W$ be the image of $w-1$ in $C$. For every $\sigma\in \pi_1(\ctau,\oeta)$,  
we put $W^{(\sigma)}=(1\times\sigma)^*(W)\in C$. It follows from 
\eqref{equivb} that the closed embedding $\odelta^{(\sigma)}$ is defined by the equation
$W^{(\sigma)}$. Hence, we have $\cJ_\sigma=(W^{(\sigma)}, t^r)$.
For every $\sigma,\sigma'\in \pi_1(\ctau,\oeta)$, we have by \eqref{nbc4a}  
\[
W^{(\sigma')}=\frac{\sigma^*(t)}{\sigma'^*(t)}(W^{(\sigma)}+1)-1.
\]
By construction of $(X\times^{\log}_k\oS)_{(r)}^{(\sigma)}$, 
$t$ is not a zero divisor and $W^{(\sigma)}$ is a multiple of $t^r$ there. 
Therefore, the inverse image of $(W^{(\sigma')},t^r)$ over $(X\times^{\log}_k\oS)_{(r)}^{(\sigma)}$ 
is an invertible ideal, equal to the inverse image of the invertible ideal 
$(t^r,\frac{\sigma^*(t)}{\sigma'^*(t)}-1)$ of $\co_{\oS}$. 

(iii) It follows from \eqref{equivb}.  

(iv) Let  $\pi\colon \fX\rightarrow  X\times_k^{\log}\oS$ be the canonical morphism.
Since $\pi$ is proper and surjective, and since the special fiber of $X\times_k^{\log}\oS$ is isomorphic to $\mG_{m,\os}$ 
\eqref{nbc4}, it is enough to prove that all fibers of $\pi$ are connected,
or equivalently that the canonical morphism $\co_{X\times_k^{\log}\oS} \rightarrow \pi_*(\co_{\fX})$
is an isomorphism. We know that  $\pi_*(\co_{\fX})$ is a coherent $(\co_{X\times_k^{\log}\oS})$-algebra.
Since $X\times_k^{\log}\oS=\Spec(\cA)$ is affine, we are reduced to showing that $\cB=\Gamma(\fX,\co_\fX)$
is isomorphic to $\cA$. On the one hand, $\cA$ is a normal domain because $X\times_k^{\log}\oS$
is smooth over $\oS$. On the other hand, $\cB$ is a domain with the same fraction field as $\cA$
because $\pi$ is a blow-up. Since $\cB$ is finite over $\cA$, we conclude that $\cB\simeq \cA$.

(v) We provide two proofs. The first one uses rigid geometry. 
We keep the notation of (iv). Recall that $\hfX^\rig$ is an annulus. Each open $\fX_{(\sigma)}$ of $\fX$ defines 
a closed subdisc $D_\sigma$ of $\hfX^\rig$. The action of $\pi_1(\ctau,\oeta)$ on $\fX$ induces an action on $\hfX^\rig$. 
For every $\sigma,\sigma' \in \pi_1(\ctau,\oeta)$, we have $\sigma(D_{\sigma'\sigma})=D_{\sigma'}$. 
By \ref{nbc6}(i), if $\sigma\not\in \pi_1(\eta,\oeta)$, then $D_\sigma$ is not contained in $D_{\id}$ (cf. \ref{nbc6e}).
Hence, for $\sigma,\sigma'\in \pi_1(\ctau,\oeta)$ such that $\sigma'\not\in \pi_1(\eta,\oeta)\sigma$, 
the disks $D_\sigma$ and $D_{\sigma'}$ are disjoint, which implies the proposition.

The second proof. For every $\sigma\in \pi_1(\ctau,\oeta)$,
we know, by (ii), (iii) and \ref{nbc5}, that $\fX_{(\sigma)} \cap \fX_{\os}$ is an affine line over $\os$. 
Let $\sigma,\sigma' \in \pi_1(\ctau,\oeta)$ be such that $\sigma'\not\in \pi_1(\eta,\oeta)\sigma$. 
If $\fX_{(\sigma)}\cap \fX_{(\sigma')} \cap \fX_{\os}$ is not empty,
it is dense open in both $\fX_{(\sigma)} \cap \fX_{\os}$  and $\fX_{(\sigma')} \cap \fX_{\os}$.
So the projective completions of these two affine lines are equal~; we denote it by $\rP$. 
We claim that the strict transform of $\odelta^{(\sigma')}(\oS)$ in $\fX$ is not contained in 
$\fX_{(\sigma)}$, while it is clearly contained in $\fX_{(\sigma')}$.
Indeed, we are reduced by (iii) and \eqref{equivb} to the case where $\sigma=\id$ and $\sigma'\not\in \pi_1(\eta,\oeta)$;
then the claim follows from (ii) and \ref{nbc6}(i) (cf. \ref{nbc6e}).
We conclude that $\fX_{(\sigma)} \cap \fX_{\os}$  and $\fX_{(\sigma')} \cap \fX_{\os}$ are different.
Then $\rP=(\fX_{(\sigma)} \cap \fX_{\os})\cup (\fX_{(\sigma')} \cap \fX_{\os})$; 
in particular, $P$ is open in $\fX_{\os}$. Since $P$ is projective over $\os$, it is also closed in $\fX_{\os}$.  
Therefore $\rP$ is a connected component of $\fX_{\os}$. 
Hence, $\fX_{\os}$ is not connected as it contains also the strict transform 
of the special fiber $Y\times_k\os$ of $X\times_k^{\log}\oS$ \eqref{nbc4}.  
We get a contradiction with (iv).

\subsection{}\label{nbc10}
We can now prove proposition \ref{nbc3}(iii). Assume first that $n-\nu(c)=2r$ is even. 
We constructed in \ref{nbc9} canonical morphisms of representations of $\pi_1(\eta,\oeta)$
\[
\cD_\oeta\stackrel{u}{\rightarrow}(\Phi^1(\pr_1^*(\cG_!)\otimes\ocL_{\psi_0}(x \cx)))_{(\os,\oeta)}
\stackrel{v}{\rightarrow} \cD_\oeta,
\]
such that $v\circ u$ is the identity. As a $\omQ_\ell$-vector space, 
$\cD_\oeta$ corresponds to the contribution of the nearby cycle complex of 
$\pr_1^*(\cG_!)\otimes\ocL_{\psi_0}(x c)$ over the dilatation $(X\times_k^{\log}S)_{(r)}$.
We consider the sheaf $\pr_1^*(\cG_!)\otimes\ocL_{\psi_0}(x c)$ over $X\times_k\oeta$, equiped with 
the action of $\pi_1(\ctau,\oeta)$ by transport of structure, and the blow-up $\fX$ of 
$X\times_k^{\log}\oS$ defined in \ref{equiv1}. Then $\cD_\oeta$ corresponds also 
to the contribution of the nearby cycle complex of $\pr_1^*(\cG_!)\otimes\ocL_{\psi_0}(x c)$ over 
the open subscheme $\fX_{(\id)}$ of $\fX$ (cf. \ref{equiv1}(ii)). 
By \ref{equiv1}(iii), for every $\sigma\in \pi_1(\ctau,\oeta)$, 
$\sigma\circ u$ is the contribution of the nearby cycle complex of $\pr_1^*(\cG_!)\otimes\ocL_{\psi_0}(x c)$ over 
$\fX_{(\sigma)}$. If $\sigma\not\in \pi_1(\eta,\oeta)$,
then we have $\fX_{(\sigma)}\cap \fX_{(\id)} \cap \fX_{\os}=\emptyset$ by \ref{equiv1}(v),
and hence $v\circ \sigma\circ u=0$, which implies the required assertion.
 
Assume next that $n-\nu(c)$ is odd, and consider as in \ref{nbc9} the base change $\uX\rightarrow X$ obtained
by taking a square root of $t$. We keep the same notation, moreover, we put 
$\ucG$ the pull-back of $\cG$ to $\uU$ and $\ucG_!$ the extension by zero of $\ucG$ to $\uX$.
We consider the complex of vanishing cycles $\Phi(\pr_1^*(\ucG_!)\otimes\ocL_{\psi_0}(x \cx))$ 
in $\bD^b_c(\uX\times_k\ctau,\omQ_\ell)$,  
relatively to the projection $\pr_2\colon \uX\times_k\cT\rightarrow \cT$. 
It follows from the even case that the morphism of sheaves over $\ctau$ 
\[
\cf_*(\rho_*(\rho^*\cD))\rightarrow (i_{\us}\times_k1)^*\Phi^1(\pr_1^*(\ucG_!)\otimes\ocL_{\psi_0}(x \cx))
\]
induced by the trace morphism $\cf_*\rho_*\rho^* \cf^*\rightarrow \id$, is injective.
Moreover, by the definition of $\cD$ in this case \eqref{nbc9}, we have a commutative diagram
\[
\xymatrix{
{\rho_*(\rho^*\cD)}\ar[r]&{\cf^*((i_{\us}\times_k1)^*\Phi^1(\pr_1^*(\ucG_!)\otimes\ocL_{\psi_0}(x \cx)))}\\
{\cD}\ar[r]\ar[u]&{\cf^*((i_s\times_k1)^*\Phi^1(\pr_1^*(\cG_!)\otimes\ocL_{\psi_0}(x \cx)))}\ar[u]}
\]
where the vertical arrows are induced by the adjunction map $\id\rightarrow \rho_*\rho^*$
and (\cite{sga7-2} XIII 2.1.7.1). The required assertion follows since $\cD\rightarrow \rho_*(\rho^*\cD)$
is injective.

\subsection{}\label{pfnbc3bis}
We can now prove proposition \ref{nbc3bis}, that will not be used in the sequel of this article. 

(i) By the $t$-exactness of the functor $\Psi$ (\cite{bbd} 4.4.2 and \cite{illusie} 4.2), 
it is enough to prove the first statement, which amounts to proving that 
$\Psi_\eta(\pr_1^*(\cG)\otimes\cL_{\psi_0}(cx))$ is supported on a finite subgroup of $Y\otimes_k\ok=\mG_{m,\ok}$. 
We may assume $k$ algebraically closed. By \ref{nbd3}, we may reduce to the case where $\cG=\cG_\rw$. 
Since $\cG$ is trivialized by a cyclic extension of degree $p^{m+1}$ of $U$, 
we may further reduce to the case where $\cG$ is a locally constant 
sheaf of $\Lambda$-modules of rank $1$ over $U$, and $\Lambda$ is a finite field of characteristic $\ell$.
We fix injective homomorphisms $\opsi_i\colon\mZ/p^{i+1}\mZ\rightarrow \Lambda^\times$ 
$(i\geq 0)$ such that for any $z\in \mF_p$, we have $\opsi_i(p^iz)=\opsi_0(z)$.
For every point $y\in Y(k)$, we denote by $H$ the henselization of $X\times_k^{\log}S$ at $y$,
that we consider as an $S$-scheme by the morphism induced by $\pr_2$. We put
\begin{equation}\label{pfnbc3bisa}
\rho(y)=\varphi_{s}(H,H_\eta,\pr_1^*(\cG)\otimes\cL_{\opsi_0}(xc)|H_\eta)
\end{equation}
the invariant defined in \eqref{app11b}. Since the pull-back of $\pr_1^*(\cG)\otimes\cL_{\opsi_0}(xc)$ 
to $H_\oeta$ is not constant, we have 
\begin{equation}\label{pfnbc3bisb}
\Psi^0_\eta(\pr_1^*(\cG)\otimes\cL_{\opsi_0}(xc))_y=0.
\end{equation}
Then by \ref{app12}, it is enough to prove that the support 
of the function $\rho(y)$ is a finite subgroup of $Y=\mG_{m,k}$.

Let $\kappa_0$ be the generic point of $Y$ \eqref{nbc4}, 
$R_{L_0}$ be the completion of the local ring of $X\times^{\log}_{k}S$ at $\kappa_0$
(which is a discrete valuation ring), $L_0$ be the fraction field of $R_{L_0}$, 
\begin{eqnarray*}
u\colon K&\rightarrow&L_0,\\
v\colon K&\rightarrow&L_0
\end{eqnarray*}
be the homomorphisms induced respectively by $\pr_1$ and $\pr_2$. We consider $L_0$ as an extension
of $K$ by $v$. By \eqref{nbc4a}, we can identify $R_{L_0}$ with the ring $k(w)[[t]]$. 
Then the $k$-homomorphisms $u$ and $v$ are defined by $u(t)=tw$ and $v(t)=t$.
We have $u(c)/c\equiv w^{\ord(c)} \mod tR_{L_0}$ and $u(b')/b'\equiv w^{\ord(b')} \mod tR_{L_0}$.
Since $\ord(t^{n+1}cb')=0$ \eqref{legendre} and $\ord(c)\not=0$ by assumption, 
then $t^{n+1}(c-u(c))u(b')$ is a unit of $R_{L_0}$. We denote by $P$ the reduction 
of $t^{n+1}(c-u(c))u(b')$ in $k(w)$, and by $\lambda$ the reduction of $t^{n+1}cb'$ in $k$. 
Then we have $P=\lambda (1-w^{\ord(c)}) w^{\ord(b')}$. 
It is enough to prove that for any $y\in Y(k)$, we have
\begin{equation}\label{pfnbc3bisd}
\rho(y)=-\ord_y(P).
\end{equation}

The pull-back of $\pr_1^*(\cG)\otimes\cL_{\opsi_0}(cx)$ to $\Spec(L_0)$ corresponds to the Witt
vector  $u(a)+\rV^m(cu(b))\in \rW_{m+1}(L_0)$. Since $\ord(bc)=-n-\nu(b)$ and $\ord(\alpha+cb')\geq -n$,
we have $u(a)+\rV^m(cu(b))\in \fil_{n+\nu(b)}\rW_{m+1}(L_0)$ and 
\begin{equation}\label{pfnbc3bisc}
\rF^md(u(a)+\rV^m(cu(b)))\equiv(c-u(c))d(u(b))+u(b)dc\mod t^{-n+1}\Omega^1_{R_{L_0}}(\log).
\end{equation}

Let $(S',\eta',s')$ be a finite covering of $(S,\eta,s)$ such that  $(H,H_\eta,\pr_1^*(\cG)\otimes\cL_{\opsi_0}(cx)|H_\eta)_{S'}$ is stable \eqref{app8}, 
$R'$ be the completion of the local ring of $S'$, $K'$ be the fraction field of $R'$,
$R'_{L'_0}=R_{L_0}\otimes_RR'$, $L'_0$ be the fraction field of $R'_{L'_0}$. 
After replacing $(S',\eta',s')$ by a finite covering (\ref{app8}(i)), 
we may assume that $\eta'$ is inseparable over $\eta$ and the image of $dt$ by the canonical morphism 
$\Omega^1_{\eta/k}\rightarrow \Omega^1_{\eta'/k}$ vanishes. 
Since we have $d(u(b))=u(b')(wdt+tdw)$ and $dc=c'dt$ in $\Omega^1_{L_0}$, 
then the following relation holds in $\Omega^1_{L'_0}$ 
\begin{equation}
(c-u(c))d(u(b))+u(b)dc=t(c-u(c))u(b') dw.
\end{equation}
It follows that the Swan conductor of the pull-back of 
$\pr_1^*(\cG)\otimes\cL_{\opsi_0}(cx)$ to $\Spec(L'_0)$ is $[K':K] n$, and its refined Swan is 
the residue class of $t(c-u(c))u(b') dw$ in $t^{-n}\Omega^1_{R'_{L'_0}}(\log)\otimes_{R'}k$, 
which is equal to $t^{-n}P dw$. Equation \eqref{pfnbc3bisd} is proved.

(ii) The proof is similar to that of proposition \ref{nbc3}(ii), given in \ref{nbc9}.

\section{Proofs of theorems \ref{results3} and \ref{results3c}}\label{pfgm}

\subsection{}\label{pfgm-not}
We observe first that, by \ref{bext}(ii), we may reduce theorems \ref{results3} and \ref{results3c} to the case 
where the residue field of $S$ at $s$ is $k$ (in particular, we have $\tf=\cf$ in \ref{results3}). 
Hence, in this section, we denote by $S$ the spectrum of a henselian discrete valuation ring of equal 
characteristic $p$, with residue field $k$, by $s$ (resp. $\eta$) the closed 
(resp. generic) point of $S$, by $\cG$ a $\omQ_\ell$-sheaf of rank $1$ over $\eta$, and by
$v\colon S\rightarrow \rP$ and $\cv\colon S\rightarrow \crP$ two non-constant morphisms
(with the notation of \ref{dtf1}). 
We put $z=v(s)$, $\cz=\cv(s)$, $b$ and $c$ the functions on $\eta$ deduced by pull-back 
from the coordinates $x$ and $\cx$ of $\rA$ and $\crA$ respectively. 
We take again the notation of \eqref{results2}
relatively to $z$ and $\cz$, and denote by $f\colon S\rightarrow T$ and 
$\cf\colon S\rightarrow \cT$ the morphisms induced by $v$ and $\cv$ respectively.
We assume that $(\cG,b,c)$ is a Legendre triple at $s$ \eqref{legendre1}, 
and $f$ and $\cf$ are finite and étale at $\eta$. 

We denote by $R$ the completion of the ring of $S$, by $K$ the fraction field of $R$,
and also by $b$ and $c$ the images of $b$ and $c$ in $K$.  
By definition \eqref{legendre1}, there exist $\cG_\rt$ and $\cG_\rw$ 
two $\omQ_\ell$-sheaves of rank $1$ over $\eta$, satisfying the following conditions~:

{\rm (i)}\ $\cG\simeq \cG_\rt\otimes \cG_\rw$; 

{\rm (ii)}\ $\cG_\rt$ is tamely ramified; 

{\rm (iii)}\ $\cG_\rw$ is trivialized by a cyclic extension of order $p^{m+1}$ of $\eta$ $(m\geq 0)$.

{\rm (iv)}\ The pull-back of $(\cG_\rw,b,c)$ over $\Spec(K)$ is a 
Legendre triple in the sense of \eqref{legendrec}.

\begin{prop}\label{results4}
We keep the assumptions of \eqref{pfgm-not}, and assume moreover that one of the following conditions
is satisfied~:

{\rm (i)}\ $z\in \rA$;

{\rm (ii)}\ $(z,\cz)=(\infty,\czero)$ and all the slopes of $f_*(\cG)$ are $<1$;

{\rm (iii)}\ $(z,\cz)=(\infty,\cinfty)$ and all the slopes of $f_*(\cG)$ are $>1$.
 
\noindent Then the rank of $\fF^{(z,\cz)}(f_*\cG)$ is equal to the degree of $\cf$, and in case {\rm (i)}, 
we have $\cz=\cinfty$.
\end{prop}

Let $t$ be a uniformizer, $\ord$ be the valuation of $K$ normalized by $\ord(t)=1$. 
We put $b'=\frac{db}{dt}$. 
Since $(\cG,b,c)$ is a Legendre triple, we have by \ref{rappelaml} and \ref{legendre},
\begin{equation}\label{psnum1b}
-\ord(c)=\sw(\cG)+\ord(tb'/b)+\ord(b).
\end{equation}
If $z\in \rA(k)$, then we have $\cz=\cinfty$ since $\sw(\cG)\geq 1$. 
In this case, we may replace $b$ by $b-x(z)$ in the equation above.
In general, we deduce by (\cite{serre1} VI §2) that we have 
\begin{equation}\label{psnum1a}
-\ord(c)=\left\{ \begin{array}{clcr}
\sw(f_*\cG)+\rk(f_*\cG)& {\rm if}\ z\in \rA(k),\\
\sw(f_*\cG)-\rk(f_*\cG)& {\rm if}\ z=\infty.
\end{array}
\right.
\end{equation} 
Hence, we have
\begin{equation}\label{results4a}
\deg(\cf)=\left\{\begin{array}{clcr}
\sw(f_*\cG)+\rk(f_*\cG)& {\rm if}\ (z,\cz)\in \rA\times \cinfty,\\
\sw(f_*\cG)-\rk(f_*\cG)& {\rm if}\ (z,\cz)=(\infty,\cinfty),\\
\rk(f_*\cG)-\sw(f_*\cG)& {\rm if}\ (z,\cz)=(\infty,\czero).
\end{array}
\right.
\end{equation} 
The proposition follows from \eqref{results4a} and (\cite{laumon} 2.4.3), for which we give 
a new proof in \ref{apdlft6}.

\subsection{}\label{pfgm3}
We identify $S$ with the henselization of the affine line $\mA^1_k=\Spec(k[u])$ at 
the origin $0$, and put $\mG_{m,k}=\Spec(k[u,u^{-1}])$ and $\hslash\colon S\rightarrow \mA^1_{k}$
the canonical morphism. By Kummer theory, $\cG_\rt$ is the pull-back of 
a smooth $\omQ_\ell$-sheaf of rank $1$, $\cM$ over $\mG_{m,k}$, 
tamely ramified at $0$ and $\infty$ (\cite{laumon} 2.2.2.1).
On the other hand, there exists a connected affine elementary étale neighborhood $(X,s)\rightarrow (\mA^1_k,0)$ 
satifying the following properties. Let $B=\Gamma(X,\co_X)$, $U$ be the open $X-\{s\}$
of $X$, $\hbar\colon S\rightarrow X$ be the unique morphism lifting $\hslash$. 
Then,

{\rm (a)}\ $\cG_\rw$ is the pull-back of a smooth $\omQ_\ell$-sheaf of rank $1$, $\tcG_\rw$ over $U$, 
trivialized by a cyclic extension of order $p^{m+1}$ of $U$ (by using Artin-Schreier-Witt theory). 

{\rm (b)}\ There exist $\tb,\tc\in \Gamma(U,\co_X)$ such that $\hbar_U^*(\tb)=b$ and $\hbar_U^*(\tc)=c$.  

{\rm (c)}\  There exists $t\in B$, which is a parameter at $s$ and invertible on $U$.

We denote by $\tcG_\rt$ be the pull-back of $\cM$ to $U$, by $\tcG=\tcG_\rt\otimes\tcG_\rw$,
by $g\colon X\rightarrow \rP$ and $\cg\colon X\rightarrow \crP$ the $k$-morphisms such that 
$g^*_\rA(x)=\tb$ and $\cg_{\crA}^*(\cx)=\tc$.  
\begin{equation}
\xymatrix{
{\cT}\ar[d]_{\ch}&S\ar[d]^{\hbar}\ar[l]_{\cf}\ar[r]^f&T\ar[d]^h\\
{\crP}&X\ar[l]_{\cg}\ar[r]^g&{\rP}}
\end{equation}
By construction, $(\tcG,\tb,\tc)$ is a Legendre triple at $s$, and we can apply \ref{nbc3}.

\subsection{}
With the notation of \ref{nbc1} and \ref{nbc2}, we have a canonical isomorphism over $\ctau$ 
\begin{equation}\label{pfgm4a}
(i_s\times 1)^*(\Phi^1(\pr_1^*(\tcG_!)\otimes\ocL_{\psi_0}(x\cx)))\simeq \fF^{(z,\cz)}(f_*\cG).
\end{equation}
It is a consequence of the functorial properties 
of the complex of nearby cycles and the fact that $\hbar\colon S\rightarrow X$ 
is universally locally acyclic and $f\colon S\rightarrow T$ is finite.
Then propositions \ref{nbc3} and \ref{results4} 
imply theorems \ref{results3} and \ref{results3c}.  

\section{Review of Stiefel-Whitney classes}\label{stwh}

\subsection{}\label{stwh1}
In this section, $K$ denotes a field of characteristic $\not=2$, 
$\oK$ a separable closure and $G_K$ the Galois group of $\oK$ over $K$. 
We denote by $1_K$ the trivial representation of $G_K$. 
By Kummer theory, $\rH^1(K,\mZ/2\mZ)$ is identified with the group $K^\times/K^{\times 2}$. 
For $a\in K^\times/K^{\times 2}$ (or in $K^\times$), we denote by $\{a\}$ the associated
element of $\rH^1(K,\mZ/2\mZ)$, and by $\kappa_a \colon G_K\rightarrow \{\pm1\}$ its image 
by the isomorphism $\rH^1(K,\mZ/2\mZ)=\Hom(G_K,\{\pm1\})$
(which is also the character induced by the quadratic extension $K(\sqrt{a})$ of $K$).
For $a,b\in K^\times/K^{\times 2}$, we denote by $\{ a ,b \}$ the 
cup-product $\{a\}\cup\{b\}$ in $\rH^2(K,\mZ/2\mZ)$.

\subsection{}\label{stwh2}
For a non-degenerate quadratic form $\cQ=\cQ(X_1,\dots,X_n)$ of rank $n$ over $K$,  
we denote by $w_m(\cQ)\in \rH^m(K,\mZ/2\mZ)$ $(m\geq 0)$ its $m$-th Stiefel-Whitney class and by 
\[
w(\cQ)=1+w_1(\cQ)+\dots \in \rH^*(K,\mZ/2\mZ)=\prod_{m} \rH^m(K,\mZ/2\mZ)
\]
its total Stiefel-Whitney class (\cite{serre2} 1.2). Recall that, if $\cQ\sim a_1X_1^2+\dots+a_nX_n^2$,
where $a_i\in K^\times$, then we have $w(\cQ)=\prod_i(1+\{a_i\})$. If $d\in K^\times/K^{\times 2}$
is the discriminant of $\cQ$, we have $w_1(\cQ)=\{d\}$.

\subsection{}\label{stwh3}
Let $V$ be a finite dimensional complex vector space, equipped with a non-degenerate quadratic form $\cQ$,
$\rho\colon G_K\rightarrow \bO(V,\cQ)$ be a continuous orthogonal representation
of $G_K$ ({\em i.e.}, the kernel of $\rho$ is open). Deligne (\cite{deligne2} 1.3 and 5.1)
associated to $(V,\rho)$ Stiefel-Whitney classes $w_m(V)\in \rH^m(K,\mZ/2\mZ)$ $(m\geq 0)$.
The class $w_1(V)$ is identified via the isomorphism $\rH^1(K,\mZ/2\mZ)=\Hom(G_K,\{\pm1\})$
with the character $\det(V)\colon G_K\rightarrow \{\pm1\}$. 
The total Stiefel-Whitney class 
\[
w(V)=1+w_1(V)+\dots \in \rH^*(K,\mZ/2\mZ)
\]
satisfies the following properties~:

(i)\ If $V$ is an orthogonal direct sum of two subrepresentations $V'$ and $V''$, then 
$w(V)=w(V')w(V'')$. 

(ii)\ If $W$ is a totally isotropic invariant subspace of $V$, 
the quadratic form $\cQ$ on $V$ induces a quadratic form on $W^\perp/W$,
and a duality between $W$ and $V/W^\perp$, and hence a quadratic form on $W\oplus V/W^\perp$; 
then we have
\begin{equation}\label{stwh3a}
w(V)=w(W\oplus V/W^\perp)w(W^\perp/W)=(1+\{-1\})^{\dim W}w(W^\perp/W).
\end{equation}

\subsection{}\label{stwh40}
Let $L$ be a finite separable extension of $K$ contained in $\oK$, $G_L$ be the 
Galois group of $\oK$ over $L$. The discriminant of $L$ over $K$, $d_{L/K}\in K^\times/K^{\times 2}$,
is by definition the discriminant of the quadratic form $\Tr_{L/K}(x^2)$, for $x\in L$.
For $\alpha\in L^\times$, the quadratic form $\Tr_{L/K}(\alpha x^2)$, for $x\in L$,
has discriminant $d_{L/K}\rN_{L/K}(\alpha)$, where $\rN_{L/K}(\alpha)$ is the norm of $\alpha$.
We denote by $w(L, \Tr_{L/K}(\alpha x^2))$ its total Stiefel-Whitney class. 

For a complex character $\chi$ of $G_L^\ab$ (or of $G_L$), we denote by 
$\rN_{L/K}(\chi)$ the composition of $\chi$ with the transfer $G_K^\ab\rightarrow G_L^\ab$.   
For any finite dimensional complex representation $V$ of $G_L$, we have
\begin{equation}\label{stwh41}
\det(\Ind_{G_L}^{G_K}V)=\kappa_{d_{L/K}}^{\dim(V)}\cdot\rN_{L/K}(\det V).
\end{equation}
This follows from  (\cite{deligne1} 1.2)  and the fact that $\det(\Ind_{G_L}^{G_K}1_L)=\kappa_{d_{L/K}}$
(\cite{serre2} 1.4).

\begin{prop}\label{stwh4}
Let $L=K(t)$ be a finite separable extension of $K$ contained in $\oK$, 
of degree $n$, generated by an element $t\in L$, $G_L$ be the Galois group of $\oK$
over $L$, $f(X)\in K[X]$ be the minimal polynomial of $t$. We put $D=f'(t)\in L^\times$
and $\kappa_D$ the associated quadratic character \eqref{stwh1}. 
Then we have 
\begin{eqnarray}
d_{L/K}&=&(-1)^{\binom{n}{2}}\rN_{L/K}(D) \in K^\times/K^{\times 2}, \label{stwh4a}\\
w_2(\Ind_{G_L}^{G_K}\kappa_{D})&=&\binom{n}{4}\{-1,-1\}+\{d_{L/K},2\}.\label{stwh4b}
\end{eqnarray}
\end{prop}

Recall (\cite{serre1} III lem.~2) that we have 
\begin{equation}\label{stwh4c}
\Tr_{L/K}(D^{-1}t^i)=\left\{
\begin{array}{clcr}
0&{\rm if}\ 0\leq i\leq n-2,\\
1&{\rm if}\ i=n-1.
\end{array}
\right.
\end{equation}
Therefore, the discriminant of the quadratic form $\Tr_{L/K}(D^{-1}x^2)$ over $K$ is 
$(-1)^{\binom{n}{2}}\in K^\times/K^{\times 2}$, which implies equation \eqref{stwh4a}. 
By (\cite{serre2} §4 theorem~1' and §1 1.5), we have 
\begin{equation}\label{stwh4d}
w_2(\Ind_{G_L}^{G_K}\kappa_D)=w_2(L,\Tr_{L/K}(D^{-1}x^2))+\{d_{L/K},2\}.
\end{equation}
Let $m$ be the largest integer such that $2m\leq n$. We denote by $W$ the sub-$K$-vector space 
of $L$ generated by $1,t,\dots,t^{m-1}$, and by $W^\perp$ the orthognal subspace 
relatively to the quadratic form $\Tr_{L/K}(D^{-1}x^2)$. By \eqref{stwh4c}, $W$ is totally isotropic; 
we have $W^\perp/W=K t^{m}$ if $n=2m+1$, and  $W=W^\perp$ otherwise. 
We deduce by \eqref{stwh3a} and \eqref{stwh4c} that we have
\begin{equation}\label{stwh4e}
w(L,\Tr_{L/K}(D^{-1}x^2))=(1+\{-1\})^m.
\end{equation} 
Equation \eqref{stwh4b} follows from \eqref{stwh4d}, 
\eqref{stwh4e} and the fact that $\binom m 2\equiv \binom n 4 \mod 2$.

\section{Refined logarithmic different}\label{refdif}

This short section is independent of the rest of the article, and does not use conventions \eqref{not1}. 

\subsection{}\label{stwh5}
Let $K$ be a complete discrete valuation field, with residue field $k$,
$L$ be a finite separable extension of $K$. 
We denote by $\co_K$ (resp. $\co_L$) the valuation ring of $K$ (resp. $L$),
by $\fm_K$ (resp. $\fm_L$) the maximal ideal of $\co_K$ (resp. $\co_L$), 
and by $k_L$ the residue field of $\co_L$. 
Recall that the different $\cD_{L/K}$ of $L$ over $K$ is the ideal of $\co_L$ 
such that the inverse $\cD_{L/K}^{-1}$ is the maximal fractional ideal $\fa$ of $L$ 
satisfying the condition $\Tr_{L/K}(\fa)\subset \co_K$. Following Kato (\cite{kato2} 2.1),
we call {\em logarithmic different} of $L$ over $K$, and denote by $\cD_{L/K}^{\log}$, 
the fractional ideal of $L$ defined by
\begin{equation}\label{stwh5a}
\cD_{L/K}^{\log}=\fm_K^{-1}\fm_L\cD_{L/K}.
\end{equation}
In fact, the ideal $(\cD_{L/K}^{\log})^{-1}$ is 
the maximal fraction ideal $\fa$ of $L$ such that $\Tr_{L/K}(\fm_L\fa)\subset \fm_K$,
and also the minimal fraction ideal $\fa$ of $L$ such that $\Tr_{L/K}(\fa)\supset \co_K$. 
We call {\em refined logarithmic different} of $L$ over $K$, 
a generator $\delta$ of the $\co_L$-module $\cD_{L/K}^{\log}$ such that, for any $a\in \co_L$, 
we have 
\begin{equation}\label{stwh5b}
\Tr_{L/K}(\delta^{-1} a)=\Tr_{k_L/k}(\oa) \mod \fm_K.
\end{equation}
Observe that $\delta$ is unique in $L^\times/1+\fm_L$. 

\subsection{}\label{stwh5c}
Let $M$ be a finite separable extension of $L$, $\delta_{L/K}$ (resp. $\delta_{M/L}$) 
be a refined logarithmic different of $L$ over $K$ (resp. of $M$ over $L$). Then 
$\delta_{M/K}=\delta_{M/L}\delta_{L/K}$ is a refined logarithmic different of $M$ over $K$.

\subsection{}\label{stwh5d}
Assume that $L$ is totally ramified over $K$ of degree $n$.
Let $t$ be a uniformizer of $L$, 
$f(X)\in \co_K[X]$ be the minimal polynomial of $t$. Then it follows from \eqref{stwh4c} 
and (\cite{serre1} III §6) that $\delta=t^{1-n}f'(t)$ is a refined logarithmic different of $L$ over $K$.

\subsection{}\label{stwh5e}
Assume that $k$ is perfect, $K$ has characteristic $p$ and $L$ is totally ramified over $K$.
Let $x$ (resp. $t$) be a uniformizer of $K$ (resp. $L$).
Then $\delta=\frac{d\log(x)}{d\log(t)}$ is a refined logarithmic different of $L$ over $K$.
Observe first that the class of $\frac{d\log(x)}{d\log(t)}$ in $L^\times/1+\fm_L$
does not depend on the choice of $x$. Let $f(X)\in \co_K[X]$ be the minimal polynomial of $t$,
$n$ be the degree of $L$ over $K$.
Since $L$ is totally ramified over $K$, $f$ is an Eisenstein polynomial, and we may assume 
that $x=-f(0)$. Therefore, we have $\frac{dx}{dt}\in f'(t)(1+\fm_L)$ and 
$tx^{-1}\in t^{1-n}(1+\fm_L)$, and the assertion follows from \ref{stwh5d}.

\section{Local epsilon factors}\label{lef}

\subsection{}\label{lef1}
In this section, $K$ denotes a complete discrete valuation field (of equal or unequal characteristics), 
with finite residue field $k$ of order $q=p^f$, 
$\oK$ a separable closure of $K$, $W_K$ the Weil group of $\oK$ over $K$
and $I$ the inertia subgroup of $W_K$. 
We denote by $\co_K$ the valuation ring of $K$, by $\fm_K$ the maximal ideal of $\co_K$,
by $\ord$ the valuation of $K$ normalized by $\ord(K^\times)=\mZ$, by 
$\co_{\oK}$ the integral closure of $\co_K$ in $\oK$ and by $\ok$ the residue field of $\co_{\oK}$. 
In the sequel, a representation of $W_K$ stands for a pair $(V,\rho)$,
where $V$ is a finite dimensional $\omQ_\ell$-vector space and $\rho$ is a continuous 
homomorphism $W_K\rightarrow {\rm GL}(V)$ ({\em i.e.} an open subgroup of $I$ acts trivially).

The quotient group $W_K/I$ is canonically isomorphic to $\mZ$, 
generated by $\Frob$, the geometric Frobenius of $k$ (i.e. the inverse of the automorphism $x\mapsto x^q$ of $\ok$).
Class field theory provides an isomorphism 
\begin{equation}\label{lef1a}
\Rec_K\colon K^\times\stackrel{\sim}{\rightarrow} W^\ab_K,
\end{equation} 
that we normalize by mapping uniformizers of $K$ to liftings of $\Frob$ (\cite{deligne1} 2.3).
We use $\Rec_K$ to identify the isomorphism classes of representations of dimension $1$
of $W_K$ with quasi-characters (i.e. continuous homomorphisms) $K^\times \rightarrow \omQ_\ell^\times$.    
If the characteristic of $K$ is $\not=2$, we denote the Hilbert symbol over $K$ by 
\begin{equation}\label{lef1b}
(\ ,\ )_K\colon K^\times/K^{\times 2}\times K^\times/K^{\times 2}\rightarrow \{\pm1\}.
\end{equation}

\subsection{}\label{lef2}
We fix a non-trivial additive character $\psi\colon K\rightarrow \omQ_\ell^\times$,  
and the Haar measure $dx$ on the additive group of $K$ such that $\int_{\co_K}dx=1$.
We call conductor of $\psi$, and denote by $\ord(\psi)$, 
the biggest integer $n$ such that $\psi|\fm_K^{-n}=1$. 
Let $\chi$ be a quasi-character of $K^\times$.  
The conductor of $\chi$, denoted by $a(\chi)$, is $0$ if $\chi$ is unramified, 
and the smallest integer $m$ such that $\chi(1+\fm_K^m)=1$ if $\chi$ is ramified. 
The Swan conductor of $\chi$, denoted by $\sw(\chi)$, is $0$ is $\chi$ is unramified,
and $a(\chi)-1$ if $\chi$ is ramified.

\subsection{}\label{lef3}
Deligne and Langlands attached to every representation $V$ of $W_K$ a {\em local $\varepsilon$-factor} 
$\varepsilon(V,\psi)\in \omQ_\ell^\times$, characterized by the following conditions~:

(i) For any exact sequence of representations 
$0\rightarrow V'\rightarrow V\rightarrow V''\rightarrow 0$, we have 
\begin{equation}\label{lef3a}
\varepsilon(V,\psi)=\varepsilon(V',\psi)\varepsilon(V'',\psi).
\end{equation}
In particular, $\varepsilon(V,\psi)$ depends only on the class of $V$ in the Grothendieck 
group of representations of $W_K$, and we can define $\varepsilon(V,\psi)$ 
when $V$ is a virtual representation of $W_K$. 

(ii) For every finite extension $L$ of $K$ contained in $\oK$, there exists a constant 
\begin{equation}\label{lef3b}
\lambda(L/K,\psi)\in \omQ_\ell^\times
\end{equation}
such that, for any representation $V_L$ of $W_L$ and $V_K$ the induced representation of $W_K$,
we have 
\begin{equation}\label{lef3c}
\varepsilon(V_K,\psi)= \lambda(L/K,\psi)^{\dim(V_L)}\varepsilon(V_L,\psi\circ \Tr_{L/K}).
\end{equation}

(iii) If $V$ has dimension $1$, corresponding to a quasi-character $\chi\colon K^\times\rightarrow \omQ_\ell^\times$, 
then $\varepsilon(V,\psi)$ is the constant $\varepsilon(\chi,\psi)$ of the local functional equation of Tate
(\cite{deligne1} §3). Recall that if $\chi$ is unramified and $\ord(\psi)=0$, 
then $\varepsilon(\chi,\psi)=1$; and if $\chi$ is ramified, then
\begin{equation}\label{lef3d}
\varepsilon(\chi,\psi)=\int_{K^\times}\chi^{-1}(x)\psi(x) dx.
\end{equation}

We omitted the Haar measure $dx$ from the notation $\varepsilon(V,\psi,dx)$, as it has been fixed 
in \eqref{lef2}. Following Deligne (\cite{deligne1} §5), we put, for a representation $V$ of $W_K$,  
\begin{equation}\label{lef3e}
\varepsilon_0(V,\psi)=\det(-\Frob,V^I) \varepsilon(V,\psi).
\end{equation}
The function $\varepsilon_0$ satisfies clearly properties (i) and (ii) with the same constant \eqref{lef3b}. 

\begin{remas}\label{lef4}
{\rm (i)}\ For any $a\in K^\times$, we have 
\begin{equation}\label{lef4a}
\varepsilon(V,\psi(ax))=\det(V)(a) \cdot q^{\ord(a)\dim(V)} \cdot \varepsilon(V,\psi);
\end{equation}
and similarly for $\varepsilon_0$. 

{\rm (ii)}\ If $L$ is a finite, separable and unramified extension of $K$ contained in $\oK$
and  $\ord(\psi)=0$, then $\lambda(L/K,\psi)=1$. This follows from (\cite{deligne1} 5.5.3)
and the fact that $\ord(\psi\circ \Tr_{L/K})=0$.
\end{remas}

\subsection{}\label{lef5}
We fix a non-trivial additive character $\psi_k\colon k\rightarrow \omQ_\ell^\times$. 
For a character $\chi\colon k^\times\rightarrow \omQ_\ell^\times$, we denote by $\tau(\chi,\psi_k)$
the Gauss sum
\begin{equation}\label{lef5a}
\tau(\chi,\psi_k)=-\sum_{x\in k^\times}\chi^{-1}(x)\psi_k(x).
\end{equation}
We have $\tau(1,\psi_k)=1$. If the characteristic $p$ of $k$ is odd, 
we denote by $\kappa_0\colon k^\times\rightarrow \{\pm1\}$ the unique character of order $2$, 
and by $G_{\psi_k}$ the quadratic Gauss sum associated to $\psi_k$, defined by
\begin{equation}\label{lef5b}
G_{\psi_k}=\sum_{x\in k}\psi_k(x^2).
\end{equation}
Then we have $\tau(\kappa_0,\psi_k)=-G_{\psi_k}$ and, by (\cite{sga4.5} [Sommes trig.] 4.4),
\begin{equation}\label{lef5c}
q=\kappa_0(-1)G_{\psi_k}^2.
\end{equation}

\subsection{}\label{lef6}
We call a {\em $\psi_k$-gauge} of $\psi$ an element $\beta\in K^\times$ such that,
for any $a\in \co_K$, with residue class $\oa$ in $k$, we have 
\begin{equation}\label{lef6a}
\psi(\beta^{-1}a)=\psi_k(\oa).
\end{equation}
Such an element $\beta$ exists, is unique in $K^\times/1+\fm_K$, and we have $\ord(\beta)=\ord(\psi)+1$.

\begin{prop}\label{lef7}
Let $\chi$ be a quasi-character of $K^\times$,
$\beta\in K^\times$ be a $\psi_k$-gauge of $\psi$, $\pi$ be a uniformizer of $K$.

{\rm (i)}\  Assume $\chi$ is at most tamely ramified (i.e., $a(\chi)\leq 1$), 
and let $\chi_k\colon k^\times\rightarrow \omQ_\ell^\times$
be the character defined by $\chi$. Then we have 
\begin{equation}\label{lef7a}
\varepsilon_0(\chi,\psi)=-\chi(\beta) \ q^{\ord(\psi)} \ \tau(\chi_k,\psi_k).
\end{equation}

{\rm (ii)}\ Assume $p\not=2$ and $\chi$ is wildly ramified (i.e., $a(\chi)\geq 2$).
We put $n=\sw(\chi)$ and let $r$ be the smallest integer such that $2r\geq n$; so we have 
$n=2r$ or $n=2r-1$. 
Let $c$ be an element of $K^\times$ such that, for any $x\in \fm_K^r$, we have 
\begin{equation}\label{lef7b}
\chi(1+x+\frac{x^2}{2})=\psi(cx).
\end{equation}
Then $\ord(\beta c)=-n$ and $c$ is unique in $K^\times/1+\fm_K^{n-r+1}$.  
We have 
\begin{equation}\label{lef7c}
\varepsilon_0(\chi,\psi)=\chi^{-1}(c)\psi(c)q^{-\ord(c)}\kappa_0(-1)^{\binom{-n}{2}}G_{\psi_k}^{-n-1}
\times \left\{
\begin{array}{clcr}
1&{\rm if} \ n\ {\rm is \ odd},\\
(-2\beta c,\pi)_K &{\rm if} \ n\ {\rm is \ even}.
\end{array}
\right.
\end{equation}
\end{prop}

(i) Assume that $\chi$ is unramified. By \eqref{lef4a}, we may assume that $\beta$ is a uniformizer
of $K$ and $\ord(\psi)=0$. Then we have $\varepsilon(\chi,\psi)=1$, and both sides of \eqref{lef7a} 
are equal to $-\chi(\beta)$.

Assume that $\chi$ is tamely ramified. By \eqref{lef4a}, we may assume that $\beta=1$ and $\ord(\psi)=-1$.
Then it follows from \eqref{lef3d} that we have $\varepsilon(\chi,\psi)=\varepsilon_0(\chi,\psi)
=-q^{-1}\tau(\chi_k,\psi_k)$. 

(ii) Since $3r\geq a(\chi)$, for any $x,y\in \fm^r_K$, we have 
\[
\chi(1+x+\frac{x^2}{2})\chi(1+y+\frac{y^2}{2})=\chi(1+x+y+\frac{(x+y)^2}{2}),
\]
which implies easily the existence and the uniqueness of $c$; the valuation of $c$ is clear. 

Let $m$ be the smallest integer such that $2m>n=\sw(\chi)$; so $m=r$ if $n=2r-1$,
and $m=r+1$ if $n=2r$. In both cases, we have $n+1=m+r$. 
For any $x\in \fm_K^m$, we have $\chi(1+x)=\psi(cx)$. We compute the integral 
\[
\varepsilon_0(\chi,\psi)=\varepsilon(\chi,\psi)=\int_{K^\times}\chi^{-1}(x)\psi(x)dx
\]
by splitting it according the classes $K^\times/1+\fm_K^m$. Only the classes contained in $c(1+\fm_K^r)$
remain~: 
\[
\varepsilon_0(\chi,\psi)=\int_{c(1+\fm^r_K)}\chi^{-1}(x)\psi(x)dx=q^{-\ord(c)}\chi^{-1}(c)\psi(c)
\int_{\fm^r_K}\chi^{-1}(1+x)\psi(cx)dx.
\]
For any $x\in \fm^r_K$, we have 
\[
\chi^{-1}(1+x)\psi(cx)=\chi^{-1}(1+x)\chi(1+x+\frac{x^2}{2})=\chi(1+\frac{x^2}{2})
=\psi(\frac{cx^2}{2}).
\]
We deduce that 
\begin{equation}\label{lef7d}
\varepsilon_0(\chi,\psi)=\chi^{-1}(c)\psi(c)q^{-\ord(c)}\int_{\fm^r_K}\psi(\frac{cx^2}{2})dx.
\end{equation}

If $n=2r-1$ is odd, then $\ord(c)+2r=-\ord(\psi)$ and $\int_{\fm^r_K}\psi(\frac{cx^2}{2})dx=q^{-r}$;
so equation \eqref{lef7c} follows by \eqref{lef5c}. 

Assume that $n=2r$ is even, so $\ord(c)+2r=-\ord(\beta)$. By \eqref{lef6a}, we have 
\[
\int_{\fm^r_K}\psi(\frac{cx^2}{2})dx=q^{-r-1}\sum_{x\in \fm^r_K/\fm^{r+1}_K}\psi_k(\overline{\frac{\beta c x^2}{2}})
=q^{-r-1} G_{\psi_k} (2\beta c,\pi)_K.
\]
So equation \eqref{lef7c} follows by \eqref{lef5c} and the relation $(-1,\pi)_K=\kappa_0(-1)$.

\begin{prop}\label{lef8}
Assume $p\not=2$ and $\ord(\psi)=0$. Let $\beta\in K^\times$ be a $\psi_k$-gauge of $\psi$ \eqref{lef6}, 
$L$ be a finite, separable, totally ramified extension of $K$ of degree $n$, 
$\pi_L$ be a uniformizer of $L$, 
$\cD_{L/K}$ be the different of $L$ over $K$, $\delta$ be a refined logarithmic 
different of $L$ over $K$ \eqref{stwh5}, $m=\ord_L(\cD_{L/K})$, 
where $\ord_L$ is the valuation of $L$ normalized by $\ord_L(\pi_L)=1$. 
Then we have 
\begin{equation}\label{lef8a}
\lambda(L/K,\psi)=\kappa_0(-1)^{\binom{m+1}{2}}G_{\psi_k}^{-m}
\times \left\{
\begin{array}{clcr}
1&{\rm if}\ m\ {\rm is\ even},\\
(-2\beta \delta,\pi_L)_L &{\rm if}\ m\ {\rm is\ odd}.
\end{array}
\right.
\end{equation}  
\end{prop} 

Let $f(X)\in \co_K[X]$ be the minimal polynomial of $\pi_L$, $D=f'(\pi_L)$ (which is a generator of $\cD_{L/K}$),
$\kappa_D\colon W_L\rightarrow \{\pm1\}$ be the character defined by the class of $D$ 
in $L^\times/L^{\times 2}$ \eqref{stwh1}, $V_K=\Ind_{W_L}^{W_K}\kappa_D$, $\psi_L=\psi\circ \Tr_{L/K}$.
We have 
\begin{equation}\label{lef8b}
\lambda(L/K,\psi)=\frac{\varepsilon(V_K,\psi)}{\varepsilon(\kappa_D,\psi_L)}.
\end{equation}
It is clear that $V_K$ is an orthogonal representation of $W_K$. 
By \eqref{stwh41} and \eqref{stwh4a}, the determinant of $V_K$ is the unramified character 
$\kappa_{-1}^{\binom{n}{2}}$. Therefore, by (\cite{serre3} theo.~1), the Artin conductor 
$a(V_K)$ of $V_K$ is even. Let $r$ be the smallest integer such that $2r\geq m=\ord_L(D)$. 
Since $a(V_K)=m+a(\kappa_D)$, we have $a(V_K)=2r$; moreover, $\kappa_D$ is unramified 
if and only if $m=2r$ is even. 

We identify $\rH^2(K,\mZ/2\mZ)$ with $\{\pm1\}$ by the isomorphism $\inv_K$, and the Hilbert symbol
$(\ ,\ )_K$ with the pairing $\{ \ , \ \}$ induced by the cup-product \eqref{stwh1}.  
By (\cite{deligne2} 1.5), since $\det(V_K)$ is unramified and $\ord(\psi)=0$, we have 
\begin{equation}\label{lef8c}
\varepsilon(V_K,\psi)=w_2(V_K) q^r,
\end{equation}
where $w_2(V_K)\in \{\pm1\}$ is the second Stiefel-Whitney class of $V_K$ \eqref{stwh3}. 
Since $(-1,-1)_K=(-1,2)_K=1$, we deduce from \ref{stwh4} that we have
\begin{equation}\label{lef8d}
w_2(V_K)=(d_{L/K},2)_K=(D,2)_L=\left\{
\begin{array}{clcr}
1& {\rm if} \ m\ {\rm is\ even},\\
(2,\pi_L)_L&{\rm if} \ m\ {\rm is\ odd}.
\end{array}
\right.
\end{equation}
To prove \eqref{lef8d}, we expressed the Hilbert symbol in terms of the tame symbol. 
By \eqref{lef5c}, we have 
\begin{equation}
q^r=\kappa_0(-1)^rG_{\psi_k}^{2r}=\kappa_0(-1)^{\binom{m+1}{2}} \times \left\{
\begin{array}{clcr}
G_{\psi_k}^m& {\rm if} \ m=2r,\\
G_{\psi_k}^{m+1}&{\rm if} \ m=2r-1.
\end{array}
\right.
\end{equation}

We put $\psi'_L(x)=\psi_L(D^{-1}x)$. Then we have $\ord_L(\psi'_L)=0$
and $\beta \pi_L^{1-n}$ is a $\psi_k$-gauge of $\psi'_L$.
Indeed, since $D \pi_L^{1-n}$ is a refined logarithmic different of $L$ over $K$ \eqref{stwh5d}, 
for any $a\in \co_L$, with residue class $\oa$ in $\co_L/\pi_L\co_L$, we have by \eqref{stwh5b}
\[
\psi'_L(\beta^{-1} \pi_{L}^{n-1}a)=\psi(\beta^{-1}\Tr_{L/K}(D^{-1} \pi_{L}^{n-1}a))=\psi_k(\oa).
\]
Since $\kappa_D(D)=(D,D)_L=(D,-1)_L=\kappa_0(-1)^m$, we have by \eqref{lef4a} 
\begin{equation}
\varepsilon(\kappa_D,\psi_L)=\kappa_D(D) q^m\varepsilon(\kappa_D,\psi'_L)=G_{\psi_k}^{2m}
\varepsilon(\kappa_D,\psi'_L). 
\end{equation}
If $m=2r$ is even, then $\kappa_D$ is unramified and we have $\varepsilon(\kappa_D,\psi'_L)=1$,
which implies equation \eqref{lef8a} in this case. Assume that $m=2r-1$ is odd, so $\kappa_D$ is tamely ramified;
in particular, the character $k^\times \rightarrow \{\pm 1\}$ defined by $\kappa_D$
is non-trivial, and hence is equal to $\kappa_0$. By \ref{lef7}(i), we have 
\begin{equation}
\varepsilon(\kappa_D,\psi'_L)=
\varepsilon_0(\kappa_D,\psi'_L)=-\kappa_D(\beta \pi_{L}^{1-n})\tau(\kappa_0,\psi_k)=\kappa_D(\beta \pi_{L}^{1-n})G_{\psi_k}.
\end{equation}
Since $\beta \pi_L^{1-n}$ is a uniformizer of $L$, we have  $\ord_L(D\beta \pi_L^{1-n})=2r$ and 
\[
\kappa_D(\beta \pi_L^{1-n})=(D,\beta \pi_L^{1-n})_L=(-D\beta \pi_L^{1-n},\beta \pi_L^{1-n})_L=
(-\beta D \pi_L^{1-n},\varpi_L)_L,
\]  
where $\varpi_L$ is any uniformizer of $L$. Moreover, a refined logarithmic different of $L$ over $K$ being unique 
in $L^\times/1+\pi_L\co_L$, we have $(D \pi_L^{1-n},\varpi_L)_L=(\delta,\varpi_L)_L$,
which proves equation \eqref{lef8a} in this case.

\subsection{}\label{lef9p}
Assume $K$ is of equal characteristic $p$. Recall that we fixed a non-trivial additive 
character $\psi_k\colon k\rightarrow \omQ_\ell^\times$.
We denote by $\res\colon \Omega^1_{K}\rightarrow k$ the residue homomorphism
and by $\ord\colon \Omega^1_K-\{0\}\rightarrow \mZ$ the valuation defined
by $\ord(xdy)=\ord(x)$, if $x,y\in K^\times$ and $\ord(y)=1$.  
For a non-zero element $\omega$ of  $\Omega^1_{K}$, we denote by 
$\psi_\omega\colon K\rightarrow \omQ_\ell^\times$ 
the non-trivial additive character defined, for any $a\in K$, by 
\begin{equation}\label{lef9pa}
\psi_\omega(a)=\psi_k(\res(a\omega)).
\end{equation} 
Let $x$ be a uniformizer of $K$, $\beta$ be the element of $K^\times$ such that $\omega=\beta x^{-1}dx$. 
Then $\beta$ is a $\psi_k$-gauge of $\psi_\omega$ and we have $\ord(\psi_\omega)=\ord(\omega)=\ord(\beta)-1$.

\begin{cor}\label{lef10}
Assume $K$ is of equal characteristic $p\not=2$. Let $L$ be a finite, separable, totally ramified extension 
of $K$ of degree $n$, $x$ be a uniformizer of $K$, $t$ be a uniformizer of $L$, 
$d_{L/K}\in K^\times/K^{\times 2}$ be the discriminant of $L$ over $K$ \eqref{stwh40}.
We put $x'=\frac{dx}{dt}$ and $m=\ord_L(x')$, where $\ord_L$ is the valuation of $L$ normalized by $\ord_L(t)=1$. 
Then we have 
\begin{eqnarray}
d_{L/K}&=&(-1)^{\binom{n}{2}}\rN_{L/K}(t^nx'/x), \label{lef10a}\\
\lambda(L/K,\psi_{dx})&=&\kappa_0(-1)^{\binom{m+1}{2}}G_{\psi_k}^{-m}\times 
\left\{
\begin{array}{clcr}
1&{\rm if}\ m \ {\rm is \ even},\\
(2x',t)_L&{\rm if}\ m \ {\rm is \ odd}.
\end{array}
\right.\label{lef10b}
\end{eqnarray}
\end{cor}

Let $f(X)\in \co_K[X]$ be the minimal polynomial of $t$.
Since $t^{1-n}f'(t)$ and $tx'/x$ are refined logarithmic differents of $L$ over $K$ \eqref{stwh5}, 
the quotient $f'(t)(t^nx'/x)^{-1}$ belongs to $1+t\co_L$. Hence, equation \eqref{lef10a} follows from \eqref{stwh4a}. 
On the other hand, $x$ is a $\psi_k$-gauge of $\psi_{dx}$ and we have $\ord(\psi_{dx})=0$. 
Then equation \eqref{lef10b} follows from \ref{lef8}.

\begin{prop}\label{lef9}
Assume $K$ is of equal characteristic $p$. Let $\chi\colon K^\times \rightarrow \omQ_\ell^\times$ be a 
wildly ramified quasi-character of Swan conductor $n=\sw(\chi)\geq 1$, 
$c\in K^\times$, $\omega$ be a non-zero element of $\Omega^1_K$, 
$\psi_m\colon \mZ/p^{m+1}\mZ \rightarrow \omQ_\ell^\times$ $(m\geq 0)$   
be injective homomorphisms. We assume the following conditions satisfied~:

{\rm (i)}\ There exist a character $\chi_\rw\colon K^\times \rightarrow \mZ/p^{m+1}\mZ$ $(m\geq 0)$ 
and a tamely ramified quasi-character $\chi_\rt\colon K^\times\rightarrow \omQ_\ell^\times$
such that $\chi=\chi_\rt\cdot (\psi_m^{-1}\circ \chi_\rw)$. 
We denote by $\gamma\in \rH^1(K,\mZ/p^{m+1}\mZ)$ the cohomology class corresponding to $\chi_\rw$
by the reciprocity isomorphism \eqref{lef1a}.

{\rm (ii)}\ There exists $a\in \fil_n\rW_{m+1}(K)$ such that $\delta_{m+1}(a)=\gamma$ \eqref{deltam} and
\begin{equation}\label{lef9a}
2\ \ord(\rF^m da+c \omega)\geq -n,
\end{equation}
where $\rF^md$ is the homomorphism defined in \eqref{calculus2a}. 

{\rm (iii)}\ $\psi_k=\psi_0\circ \Tr_{k/\mF_p}$ and $\psi_m(p^ma)=\psi_0(a)$ for any $a\in \mF_p$,
where $p^m a$ denotes the embedding $\mF_p\rightarrow \mZ/p^{m+1}\mZ$
induced by the multiplication by $p^m$ on $\mZ$. 

Let $\psi_\omega\colon K\rightarrow \omQ_\ell^\times$ 
be the additive character defined in \eqref{lef9pa}, 
$r$ be the smallest integer such that $2r\geq n$. 
Then, for any $x\in \fm_K^r$, we have 
\begin{equation}\label{lef9b}
\chi(1+x+\frac{x^2}{2})=\psi_\omega(cx).
\end{equation}
\end{prop}

We will deduce the proposition from Witt's explicit reciprocity law 
according to Fontaine (\cite{fontaine} 2.4.3). 
Let $W=\rW(k)$, $W_{m+1}=\rW_{m+1}(k)$, $\co_\cE$ be the $p$-adic completion of the ring $W((t))$
of Laurent power series over $W$ in the variable $t$ (which is an absolutely unramified,
complete, discrete valuation ring), $\cE$ be the fraction field of $\co_\cE$. 
We identify the residue field of $\co_\cE$
with $K$ by mapping the residue class of $t$ to a uniformizer of $K$. 
We denote by $\hOmega^1_{\co_\cE/W}$ the module of continuous differential forms of $\co_\cE$
over $W$ and by $\res_t\colon \hOmega^1_{\co_\cE/W}\rightarrow W$ the residue homomorphism 
(cf. \cite{fontaine} 2.2).
For $z\in \rW_{m+1}(K)$ and $u\in K^\times$, we denote by $[z,u)_m$ the element 
of $\mZ/p^{m+1}\mZ\subset \rW_{m+1}(\oK)$ defined by 
\[
[z,u)_m=g_u(\xi)-\xi,
\]
where $\xi$ is an element of $\rW_{m+1}(\oK)$ such that 
$\rF(\xi)-\xi=z$, $\rF$ is the Frobenius homomorphism, 
and $g_u\in G_K^\ab$ is the image of $u$ by the reciprocity isomorphism \eqref{lef1a}. 
If we put $\co_{\cE,m}=\co_{\cE}/p^{m+1} \co_{\cE}$, we have a homomorphism 
\[
w_m\colon \rW_{m+1}(K)\rightarrow \co_{\cE,m}
\]
defined for an element $z=(z_0,z_1,\dots,z_m)$ of $\rW_{m+1}(K)$, by 
$w_m(z)=\sum_{0\leq j\leq m}p^j(\tz_j)^{p^{m-j}}$, where $\tz_j$ is any lifting 
of $z_j$ in $\co_{\cE,m}$. Then, if $z\in \rW_{m+1}(K)$ and if $\tu$
is a unit of $\co_{\cE,m}$ lifting an element $u$ of $K^\times$, we have 
\begin{equation}\label{lef9d}
[z,u)_m=-\Tr_m(\res_m(w_m(z)d\log \tu)),
\end{equation}
where $\Tr_m$ (resp. $\res_m$) is the reduction modulo $p^{m+1}$
of the trace homomorphism of $W$ over $\mZ_p$ (resp. $\res_t$). 
Notice that the minus sign in \eqref{lef9d} does not appear in (\cite{fontaine} 2.4.3)
because the reciprocity map used there is the inverse of \eqref{lef1a}
(the reciprocity map used by Fontaine sends uniformizers to arithmetic Frobenius).   

Multiplication by $p^m$ on $\hOmega^1_{\co_\cE/W}$ induces a homomorphism 
\[
\Omega^1_K=\hOmega^1_{\co_\cE/W}\otimes_{\co_\cE}K\rightarrow 
\hOmega^1_{\co_\cE/W}\otimes_{\co_\cE}\co_{\cE,m},
\]
that we abusively denote by a multiplication by $p^m$. For any $z\in \rW_{m+1}(K)$, we have 
\begin{equation}
dw_m(z)=p^m\rF^md(z).
\end{equation}

We can now prove the proposition. Since $r\geq 1$, we may assume $\chi_\rt=1$. 
We put $a=(a_0,\dots,a_m)\in \rW_{m+1}(K)$ and $\nu_i=\ord(a_i)$.
Let $\ta_i$ be a lifting of $a_i$ in $t^{\nu_i}W_{m+1}[[t]]$ $(0\leq i\leq m)$, 
$x\in \fm^r_K$, $\tx$ be a lifting of $x$ in $t^{r}W_{m+1}[[t]]$. 
It follows from the choice of the $\ta_i$ and the fact that $3r\geq n+1$, that 
we have $w_m(a)\in t^{-n}W_{m+1}[[t]]$, and $\res_m(w_m(a)\tx^id\tx)=0$  for $i\geq 2$. 
Therefore, we have
\begin{eqnarray*}
[a,1+x+\frac{x^2}{2})_m&=&-\Tr_m(\res_m( w_m(a)\frac{1+\tx}{1+\tx+\frac{\tx^2}{2}}d\tx))\\
&=&-\Tr_m(\res_m(w_m(a)d\tx))\\
&=&\Tr_m(\res_m(\tx d(w_m(a))))\\
&=&p^m\Tr_{k/\mF_p}(\res_{K}(x\rF^md(a))),
\end{eqnarray*}
where on the right hand side, $p^m$ denotes the embedding $\mZ/p\mZ\rightarrow \mZ/p^{m+1}\mZ$
induced by the multiplication by $p^m$ on $\mZ$. 
We conclude by \eqref{lef9a} that for any $x\in \fm^r_K$, we have
\begin{eqnarray*}
\chi(1+x+\frac{x^2}{2})&=&\psi^{-1}_m([a,1+x+\frac{x^2}{2})_m)\\
&=&\psi_0(-\Tr_{k/\mF_p}(\res_{K}(x\rF^md(a)))\\
&=&\psi_0(\Tr_{k/\mF_p}(\res_{K}(xcw)))\\
&=&\psi_\omega(cx).
\end{eqnarray*}

\section{Laumon's formula for local epsilon factors}\label{lauf}  

\subsection{}\label{lauf1}
Let $T$ be the spectrum of a henselian discrete valuation field of equal characteristic $p$, 
with finite residue field $k$ of order $q=p^f$, $\tau$ (resp. $\otau$) be the generic point 
(resp. a geometric generic point) of $T$, $G=\pi_1(\tau,\otau)$. We denote by 
$K$ the completion of the function field $k(\tau)$ of $T$ and by 
\begin{equation}\label{lauf1a}
\Rec_T\colon K^\times\rightarrow G^\ab
\end{equation}
the reciprocity homomorphism, normalized as in \eqref{lef1a}. 

Recall that we fixed a non-trivial additive character $\psi_0\colon \mF_p\rightarrow \omQ_\ell^\times$ \eqref{not1}.
Let $\psi_k\colon k\rightarrow \omQ_\ell^\times$ be the additive character $\psi_0\circ \Tr_{k/\mF_p}$. 
For a complex $C$ of $\bD^b_c(T,\omQ_\ell)$ and a non-zero meromorphic differential form 
$\omega$ on $T$ (i.e. $\omega\in \Omega^1_{k(\tau)}-\{0\}$), Laumon 
attached a local $\varepsilon$-factor $\varepsilon(T,C,\omega)\in \omQ_\ell^\times$ (\cite{laumon} 3.1.5.4). 
For any $\omQ_\ell$-sheaf $\cF$ over $\tau$, we have (with the notation of \ref{lef3})
\begin{eqnarray}
\varepsilon(T,j_*\cF,\omega)&=&\varepsilon(\cF_\otau,\psi_\omega),\label{lauf1c}\\
\varepsilon(T,j_!\cF,\omega)&=&\varepsilon_0(\cF_\otau,\psi_\omega),\label{lauf1d}
\end{eqnarray}
where $j\colon \tau\rightarrow T$ is the canonical injection, $\psi_\omega$ is the 
additive character defined in \eqref{lef9pa}. 

In the situation of \eqref{results2}, we use the notation above for $T$ and $\cT$.
We equip with a $\vee$ the objects relative to $\cT$.

\begin{teo}[\cite{laumon} 3.6.2]\label{lauf2}
The assumptions are those of \eqref{results3}, moreover, we assume that 
$k$ is finite, $z=0$ and $\cz=\cinfty$. We denote by $\cG_!$ the extension by $0$ of $\cG$ to $S$, 
and by $d$ the dimension of $\fF^{(0,\cinfty)}(f_*(\cG))$ over $\omQ_\ell$. Then we have 
\begin{equation}\label{lauf2a}
(-1)^d\det(\Rec_{\cT}(\cx^{-1}),\fF^{(0,\cinfty)}(f_*(\cG)))=\varepsilon(T,f_*(\cG_!),dx).
\end{equation}
\end{teo} 

Let $k'$ be the residue field of $S$ at $s$. First, we reduce the theorem to the case where $k'=k$. 
We denote by $0'\in \rP_{k'}(k')$ and $\cinfty'\in \crP_{k'}(k')$ the points induced by $0\in \rP(k)$ and 
$\cinfty\in \crP(k)$, by $v'\colon S\rightarrow \rP_{k'}$
and $\cv'\colon S\rightarrow \crP_{k'}$ the morphisms induced by $v$ and $\cv$, 
by $T'$ and $\cT'$ the henselizations of $\rP_{k'}$ and $\crP_{k'}$ at $0'$ and $\infty'$ respectively, 
by $f'\colon S\rightarrow T'$ and $\cf'\colon S\rightarrow \cT'$ the morphisms induced by $v'$ and $\cv'$
(or by $f$ and $\cf$) and by $w\colon \cT'\rightarrow \cT$ the canonical morphism. 
We have $T'=T\otimes_kk'$ and $\cT'=\cT\otimes_kk'$. By \ref{bext}(ii), we have 
\begin{equation}\label{rstep1}
w_*(\fF^{(0',\cinfty')}(f'_*(\cG)))\simeq \fF^{(0,\infty)}(f_*(\cG)).
\end{equation}
Let $d'$ be the dimension of $\fF^{(0',\cinfty')}(f'_*(\cG))$ over $\omQ_\ell$.
Since $\cT'$ is a finite étale covering of $\cT$ of degree $d/d'$, we deduce by (\cite{deligne1} 1.2)
that we have  
\begin{eqnarray}\label{rstep2}
\lefteqn{\det(\Rec_{\cT}(\cx^{-1}),\fF^{(0,\infty)}(f_*(\cG)))}\\
&=&\det(\Rec_{\cT}(\cx^{-1}),w_*(1))^{d'}\cdot \det(\Rec_{\cT'}(\cx^{-1}),\fF^{(0',\infty')}(f'_*(\cG))) \nonumber\\
&=&(-1)^{d-d'} \det(\Rec_{\cT'}(\cx^{-1}),\fF^{(0',\infty')}(f'_*(\cG))). \nonumber
\end{eqnarray}
Since $T'$ is a finite étale covering of $T$, we have, by \ref{lef4}(ii) and \ref{lef9p},
\begin{equation}\label{rstep3}
\varepsilon(T',f'_*(\cG_!),dx)=\varepsilon(T,f_*(\cG_!),dx).
\end{equation}

Equations \eqref{rstep2} and \eqref{rstep3} show that we may assume $k'=k$. 
We denote by $L$ the completion of the function field $k(\eta)$ of $S$, 
by $\co_L$ its valuation ring, by $t$ a uniformizer of $k(\eta)$, 
and by $\ord_L$ the valuation of $L$ normalized by $\ord_L(t)=1$.
For any $y\in L$, we put $y'=\frac{dy}{dt}$; if $y\in k(\eta)$, 
then $y'$ is well defined in $k(\eta)$ \eqref{results3b}.
We consider $L$ as a finite, separable, totally ramified extension of both $K$ and $\ccK$, 
the completions of the function fields of $T$ and $\cT$ respectively. 
Since $x$ is a uniformizer of $K$ and $\cx^{-1}$ is a uniformizer of $\ccK$, 
we have $[L:K]=\ord_L(b)$, $[L:\ccK]=-\ord_L(c)$, 
and $tb'/b$ (resp. $-tc'/c$) is a refined logarithmic different of $L$ over $K$ 
(resp. of $L$ over $\ccK$) \eqref{stwh5e}. We put $m=\ord_L(b')$, $i=-\ord(c)$, 
and denote by $d_{L/\ccK}$ the discriminant of $L$ over $\ccK$ \eqref{stwh40}. 

We denote by $\chi\colon L^\times\rightarrow \omQ_\ell^\times$ the quasi-character defined 
by the sheaf $\cG$ over $\eta$, by $(\ ,\ )_L$ the Hilbert symbol over $L$, 
by $\kappa_0\colon k^\times\rightarrow \{\pm1\}$ the unique character of order $2$,
by $\psi_k\colon k\rightarrow \omQ_\ell^\times$ the additive character $\psi_0\circ \Tr_{k/\mF_p}$,
and by $G_{\psi_k}$ the quadratic Gauss sum associated to $\psi_k$ \eqref{lef5b}.
Observe that we have the following equality of additive characters of $L$
\begin{equation}
\psi_{db}=\psi_{dx}\circ \Tr_{L/K}.
\end{equation}

On the one hand, by \ref{results3}, we have a canonical isomorphism of sheaves over $\ctau$ 
\[
\fF^{(0,\cinfty)}(f_*(\cG))\simeq \cf_*(\cG\otimes \cL_{\psi_0}(bc)\otimes\cK(-\frac{1}{2}b'c')\otimes \cQ).
\]
We deduce by \eqref{stwh41} that we have 
\begin{eqnarray}\label{lauf2b}
\lefteqn{\det(-\Rec_{\cT}(\cx),\fF^{(0,\cinfty)}(f_*(\cG)))=}\\
&&(-1)^{i} (\cx,d_{L/\ccK})_{\ccK}\det(\Rec_{S}(c),\cG\otimes \cL_{\psi_0}(bc)\otimes\cK(-\frac{1}{2}b'c')\otimes \cQ).
\nonumber
\end{eqnarray}
Moreover, we have 
\begin{eqnarray}
(\cx,d_{L/\ccK})_{\ccK}&=&(\cx,(-1)^{\binom{i}{2}}\rN_{L/\ccK}(-t^i c'/c))_{\ccK}\label{lauf2g}\\
&=&\kappa_0(-1)^{\binom{i}{2}}(c,t^i c')_L,\nonumber\\
\det( \Rec_S(c),\cG)&=&\chi(c),\label{lauf2c}\\
\det( \Rec_S(c),\cL_{\psi_0}(bc))&=&\psi^{-1}_0(-\Tr_{k/\mF_p}(\res_L (bc\frac{dc}{c})))\label{lauf2d}\\
&=&(\psi_k(\res_L (cdb)))^{-1}=(\psi_{db}(c))^{-1},\nonumber\\
\det( \Rec_S(c),\cK(-\frac{1}{2}b'c'))&=&(c,-2b'c')_L,\label{lauf2e}\\
\det( \Rec_S(c),\cQ)&=&(-1)^i G_{\psi_k}^{-i}.\label{lauf2f}
\end{eqnarray}
Equation \eqref{lauf2e} is obvious from the definitions, 
equation \eqref{lauf2g} follows from \eqref{lef10a}, equation 
\eqref{lauf2f} follows from (\cite{laumon} 1.4.3.2(ii)), 
and equation \eqref{lauf2d} is a consequence of (\cite{serre1} XIV §5 prop.~15)~:
the power $-1$ above $\psi_0$ is due to the convention in \eqref{results1},
and the minus sign before $\Tr_{k/\mF_p}$ is due to the 
convention in the definition of the reciprocity law \eqref{lef1a}.
We deduce that we have 
\begin{equation}\label{lauf2bb}
\det(-\Rec_{\cT}(\cx),\fF^{(0,\cinfty)}(f_*(\cG)))=
\kappa_0(-1)^{\binom{i}{2}}\chi(c)(\psi_{db}(c))^{-1}G_{\psi_k}^{-i}(c,-2t^ib')_L.
\end{equation}

On the other hand, we have 
\begin{equation}\label{lauf3a}
\varepsilon(T,f_*(\cG_!),dx)=\lambda(L/K,\psi_{dx})\varepsilon(S,\cG_!,db)
=\lambda(L/K,\psi_{dx})\varepsilon_0(\chi,\psi_{db}).
\end{equation}
By assumption, $(\cG,b,c)$ is a Legendre triple \eqref{legendre1}.
So $\chi$ is wildly ramified. We put $n=\sw(\chi)$ and $r$ the smallest integer such that $2r\geq n$. 
By \ref{lef9}, for any $y\in t^r\co_L$, we have 
\[
\chi(1+y+\frac{y^2}{2})=\psi_{db}(cy).
\]
For any $y\in \co_L$, with residue class $\oy$ in $k$, we have by \eqref{stwh5b}
\[
\psi_{db}((tb')^{-1}y)=\psi_{dx}(x^{-1}\Tr_{L/K}((tb'/b)^{-1}y))=\psi_k(\oy).
\]
We deduce from \ref{lef7}(ii) that we have $1+m-i=\ord(tb'c)=-n$ and
\begin{equation}\label{lauf3b}
\varepsilon_0(\chi,\psi_{db})=\chi^{-1}(c)\psi_{db}(c)q^{i}\kappa_0(-1)^{\binom{-n}{2}}G_{\psi_k}^{-n-1}
\times\left\{
\begin{array}{clcr}
1&{\rm if}\ n\ {\rm is\ odd}\\ 
(2b'c,t)_L&{\rm if}\ n\ {\rm is\ even}
\end{array}
\right.
\end{equation}
By \eqref{lef10b}, we have 
\begin{equation}\label{lauf3c}
\lambda(L/K,\psi_{dx})=\kappa_0(-1)^{\binom{m+1}{2}}G_{\psi_k}^{-m}\times 
\left\{
\begin{array}{clcr}
1&{\rm if}\ m \ {\rm is \ even},\\
(2b',t)_L&{\rm if}\ m \ {\rm is \ odd}.
\end{array}
\right.
\end{equation}

To conclude the proof, it remains to check that the product of the right hand sides of equations
\eqref{lauf2bb}, \eqref{lauf3b} and \eqref{lauf3c} is equal to $1$. 
Since we have $q=\kappa(-1)G_{\psi_k}^2$ \eqref{lef5c}, $1+m=i-n$ and hence 
$\binom{m+1}{2}=\binom{i}{2}+\binom{-n}{2}-in$, we are reduced to checking that 
\begin{equation}
1=\kappa_0(-1)^{i(n+1)}(c,-2t^ib')_L\times 
\left\{
\begin{array}{clcr}
1&{\rm if} \ m\ {\rm is \ even\ and}\ n \ {\rm is \ odd},\\
(c,t)_L&{\rm if} \ m\ {\rm is \ odd\ and}\ n \ {\rm is \ even},\\
(2b'c,t)_L&{\rm if} \ m\ {\rm and}\ n \ {\rm are \ even},\\
(2b',t)_L&{\rm if} \ m\ {\rm and}\ n \ {\rm are \ odd}.
\end{array}
\right.
\end{equation}
Since we have $\kappa_0(-1)=(-1,t)=(t,t)$, we are further reduced to 
checking that 
\begin{equation}
1=\left\{
\begin{array}{clcr}
(c,-2b')_L&{\rm if} \ m\ {\rm is \ even\ and}\ n \ {\rm is \ odd},\\
(c,-2b't)_L&{\rm if} \ m\ {\rm is \ odd\ and}\ n \ {\rm is \ even},\\
(-2b',tc)_L&{\rm if} \ m\ {\rm and}\ n \ {\rm are \ even},\\
(-2tb',tc)_L&{\rm if} \ m\ {\rm and}\ n \ {\rm are \ odd}.
\end{array}
\right.
\end{equation}
In each case, the valuations of both terms of the Hilbert symbol are even, which proves the required result.

\appendix

\section{Semi-continuity of the Swan conductor}\label{scsc}

\subsection{}\label{app1}
In this section, $(S,\eta,s)$ denotes an excellent henselian trait, of equal characteristic $p>0$, 
with algebraically closed residue field $k$, {\em i.e.}  
$S=\Spec(V)$ is the spectrum of an excellent henselian discrete valuation ring, of equal characteristic $p>0$, 
$\eta$ and $s$ are the generic and the closed points of $S$. 
We denote by $K$ the fraction field of $V$.
We fix a geometric generic point $\oeta$ of $S$, and a finite field $\Lambda$ of characteristic $\not=p$. 
A finite covering of $(S,\eta,s)$ stands for a trait $(S',\eta',s')$ equipped with a finite covering 
$S'\rightarrow S$. 

\subsection{}\label{app2}
Let $R$ be a complete discrete valuation ring, $L$ be the fraction field of $R$, $\fm$
be the maximal ideal of $R$, $L'$ be a finite separable extension of $L$, 
$R'$ be the integral closure of $R$ in $L'$. 
We say that $L'$ is a {\em stable} extension of $L$ if $\fm R'$ is the maximal ideal of $R'$.

\subsection{}\label{app3}
Let $R$ be a complete discrete valuation ring which is a $k$-algebra, 
$F$ be the residue field of $R$, $L$ be the fraction field of $R$, 
$\fm$ be the maximal ideal of $R$. We assume that $F$ is an extension of finite type of $k$.     
Then the $R$-module of absolute $1$-differential forms $\Omega^1_R$ is complete, separated,
and hence free of finite rank over $R$. We denote by $\Omega^1_R(\log)$ the sub-$R$-module
of $\Omega^1_L$ generated by $\Omega^1_R$ and $d\log(x)$ where $x$ a uniformizer of $R$ (cf. \cite{aml} 5.4). 
We put $\Omega^1_F(\log)=\Omega^1_R(\log)\otimes_{R}F$. We have a canonical exact sequence
\begin{equation}\label{app3a}
\xymatrix{
0\ar[r]&{\Omega^1_F}\ar[r]&{\Omega^1_F(\log)}\ar[r]^-(0.5)\res&{F}\ar[r]&0}.
\end{equation}

Let $\cF$ be a $\Lambda$-sheaf of rank $1$ over $\Spec(L)$. Kato \cite{kato1} associates to 
$\cF$ a Swan and a refined Swan conductors, that can also be defined using our ramification theory \cite{aml}. 
The Swan conductor $n=\sw(\cF)$ is an integer $\geq 0$, that vanishes if and only if $\cF$
is tamely ramified. The refined Swan conductor $\rsw(\cF)$ is an element of the $F$-vector space
\begin{equation}\label{app3b}
\Omega^1_F(\log)\otimes_F(\fm^{-n}/\fm^{-n+1}).
\end{equation} 

If $\cF$ is trivialized by a stable extension of $L$, then we have (\cite{kato1} remark after 6.8) 
\begin{equation}\label{app3c}
(\res\otimes 1)(\rsw(\cF))=0 \in \fm^{-n}/\fm^{-n+1}.
\end{equation}

\subsection{}\label{app4}
We denote by $\cC_S$ the following category. Objects of $\cC_S$ are normal affine $S$-schemes $H$
for which there exist an $S$-curve $X$ ({\em i.e.} a flat $S$-scheme of finite type and relative dimension $1$) 
and a closed point $x$ of $X_s$, such that $X$ is smooth over $S$ outside $x$, 
and $H$ is $S$-isomorphic to the henselization of $X$ at $x$. Let $H,H'$ be two objects of $\cC_S$.
A morphism $f\colon H'\rightarrow H$ of $\cC_S$ is a finite morphism
of $S$-schemes which is étale at the generic point of $H'$.

\subsection{}\label{app5}
Let $H$ be an object of $\cC_S$, $(S',\eta',s')$ be a finite covering of $(S,\eta,s)$.
Then $H\times_SS'$ is an object of $\cC_{S'}$ (\cite{kato3} 5.4).

\subsection{}\label{app6}
Let $H$ be an object of $\cC_S$. We denote by $H^\circ$ the set of height $1$ points of $H$,
$H^\circ_\eta=H_\eta\cap H^\circ$ and $H^\circ_s=H_s\cap H^\circ$.  Then,

(i) $H_\eta$ is geometrically regular over $\eta$, of dimension $1$, and for any $\fp\in H^\circ_\eta$,
the residue field $\kappa(\fp)$ of $H$ at $\fp$ is a finite extension of $K$. 

(ii) $H_s$ is a reduced henselian local scheme, of dimension $1$. 

Indeed, let $X$ be an $S$-curve and $x$ be a closed point of $X_s$ 
such that $X$ is smooth over $S$ outside $x$, and $H$ is $S$-isomorphic to the henselization of $X$ at $x$.  
The set of geometric points of $H$ is canonically isomorphic to the set of geometric points of $X$
which are generizations of $x$; 
moreover, the strict henselizations of $X$ and $H$ at associated geometric points are isomorphic (SGA 4 VIII 7.3).
We deduce that $H_\eta$ is regular of dimension $1$, and hence, geometrically regular over $\eta$ by \ref{app5}.
The second assertion of (i) is a consequence of (EGA IV 8.2.9 and its corollaries).
The scheme $H_s$ is the henselization of $X_s$ at $x$, which implies (ii). 

We denote by $\tH_s$ its normalization, 
which is a finite disjoint union of strictly local traits (indexed by $H^\circ_s$). We put 
\begin{equation}
\delta(H)=\dim_k(\co_{\tH_s}/\co_{H_s}).
\end{equation}

\subsection{}\label{app7}
We denote by $\fF_S$ the following category. Objects of $\fF_S$ are triples $(H,U,\cF)$, where 
$H$ is an object of $\cC_S$, $U$ is a non-empty open subscheme of $H_\eta$ and $\cF$ is a locally constant
constructible étale sheaf of $\Lambda$-modules over $U$. 
Let $(H,U,\cF), (H',U',\cF')$ be two objects of $\fF_S$. A morphism $(H',U',\cF')\rightarrow (H,U,\cF)$
of $\fF_S$ is a pair $(f,g)$ made of a morphism $f\colon H'\rightarrow H$ of $\cC_S$ such that $f(U')\subset U$, 
and a morphism $g\colon \cF'\rightarrow f_U^*\cF$, where $f_U\colon U'\rightarrow U$ is the restriction of $f$. 

Let $(S',\eta',s')$ be a finite covering of $(S,\eta,s)$. By \ref{app5}, 
the base change $S'\rightarrow S$ induces a natural functor $\fF_S\rightarrow \fF_{S'}$, 
that we denote by 
\begin{equation}\label{app7a}
(H,U,\cF)\mapsto (H,U,\cF)_{S'}.
\end{equation} 

\subsection{}\label{app75}
An object $(H,U,\cF)$ of $\fF_S$ is 
said to be {\em stable} if there exists a finite étale connected covering $U'$ of $U$ satisfying the following
conditions~:

(i) The pull-back of $\cF$ to $U'$ is constant.

(ii) The normalization $H'$ of $H$ in $U'$ belongs to $\cC_S$, and the residue fields of $H'$
at all points of $H'_\eta-U'_\eta$ are finite separable extensions of $K$. 

\begin{prop}[\cite{kato3} 6.3]\label{app8}
Let $(H,U,\cF)$ be an object of $\fF_S$. 

{\rm (i)} Let $(S',\eta',s')$ be a finite covering of $(S,\eta,s)$.
If $(H,U,\cF)$ is a stable object of $\fF_S$, then $(H,U,\cF)_{S'}$ is a stable object of $\fF_{S'}$.

{\rm (ii)} There exists $(S',\eta',s')$ a finite covering of $(S,\eta,s)$,
such that $(H,U,\cF)_{S'}$ is a stable object of $\fF_{S'}$.
\end{prop}

Proposition (i) follows from \ref{app5} and proposition (ii) follows from \cite{epp}. 

\subsection{}\label{app9}
Let $(H,U,\cF)$ be a stable object of $\fF_S$ such that $\cF$ has {\em rank} $1$ over $U$. 
For $\fp\in H^\circ$, we denote by $R_\fp$ the completion of the local ring of $H$ at $\fp$
(which is a discrete valuation ring), and by $\kappa(\fp)$ its residue field. 
Following (\cite{sga7-2} XVI), \cite{laumon2} and (\cite{kato3} 6.4), 
we call {\em total dimension} of $\cF$ at a point $\fp$, 
and denote by $\dimtot_{\fp}(\cF)$, the integer defined as follows. For $\fp\in H^\circ_\eta$,  we put 
\begin{equation}\label{app9a}
\dimtot_{\fp}(\cF)=[\kappa(\fp):K](\sw_\fp(\cF)+1),
\end{equation}
where $\sw_\fp(\cF)$ is the Swan conductor of the pull-back to $\cF$ over $\Spec(R_\fp)\times_HU$.  

For $\fp\in H^\circ_s$, we denote by $\tH_{s,\fp}$ the integral closure of $H_s$ in $\kappa(\fp)$
(which is a stricly local trait) and by $\ord_{s,\fp}$ the associated valuation of $\kappa(\fp)$, 
normalized by $\ord_{s,\fp}(\kappa(\fp)^\times)=\mZ$.
We denote also by $\ord_{s,\fp}\colon \Omega^1_{\kappa(\fp)}-\{0\}\rightarrow \mZ$ the valuation defined by 
$\ord_{s,\fp}(\alpha d\beta)=\ord_{s,\fp}(\alpha)$, if $\alpha,\beta\in \kappa(\fp)^\times$ and 
$\ord_{s,\fp}(\beta)=1$. We distinguish two cases~: 

(i)\ Assume that $\cF$ extends to a locally constant constructible sheaf of $\Lambda$-modules $\tcF$ 
over an open subscheme $\tU$ of $H$ that contains $\fp$. 
We denote by $\sw_{s,\fp}(\cF)$ the Swan conductor of the pull-back
of $\tcF$ to $\tH_{s,\fp}\times_H\tU$. We put 
\begin{equation}\label{app9b}
\dimtot_\fp(\cF)=\sw_{s,\fp}(\cF)+1.
\end{equation}
 
(ii)\ Assume that $\cF$ is ramified at $\fp$. We denote by $n=\sw_\fp(\cF)$ and $\rsw_\fp(\cF)$
the Swan and the refined Swan conductors of the pull-back of $\cF$ to $\Spec(R_\fp)\times_HU$ \eqref{app3}. 
Let $\pi$ be a uniformizer of $V$. Since $(H,U,\cF)$ is stable, we have 
$\res(\rsw_\fp(\cF)\otimes [\pi^n])=0$ \eqref{app3c}.
Hence, we can identify $\rsw_\fp(\cF)\otimes [\pi^n]$ with the image of 
an element $\omega \in \Omega^1_{\kappa(\fp)}$, which does not depend on the choice of $\pi$
up to a multiplication by an element of $k^\times$. We put 
\begin{equation}\label{app9c}
\dimtot_\fp(\cF)=-\ord_{s,\fp}(\omega).
\end{equation} 
We put 
\begin{eqnarray}
\varphi_\eta(H,U,\cF)&=&\sum_{\fp\in H_\eta-U}\dimtot_\fp(\cF),\label{app9d}\\
\varphi_s(H,U,\cF)&=&\sum_{\fp\in H^\circ_s}\dimtot_\fp(\cF).\label{app9e}
\end{eqnarray}

\begin{lem}[\cite{kato3} 6.5]\label{app10}
Let $(H,U,\cF)$ be a stable object of $\fF_S$ such that $\cF$ has rank $1$ over $U$, 
$(S',\eta',s')$ be a finite covering of $(S,\eta,s)$. We put $(H',U',\cF')=(H,U,\cF)_{S'}$. 

{\em (i)}\ For any $\fp'\in H'^\circ_s$ with image $\fp$ in $H^\circ_s$, 
we have $\dimtot_\fp(\cF)=\dimtot_{\fp'}(\cF')$. 

{\em (ii)}\ For any $\fp\in H_\eta-U$, we have 
\begin{equation}\label{app10a}
\dimtot_\fp(\cF)=\sum_{\fp'} \dimtot_{\fp'}(\cF),
\end{equation} 
where $\fp'$ runs over the points of $H'$ above $\fp$.
\end{lem}

\subsection{}\label{app11}
Let $(H,U,\cF)$ be an object of $\fF_S$  such that $\cF$ has rank $1$ over $U$. By \ref{app8}(ii), there 
exists a finite covering $(S',\eta',s')$ of $(S,\eta,s)$ such that $(H,U,\cF)_{S'}$ is a stable object of 
$\fF_{S'}$. We put  
\begin{eqnarray}
\varphi_\eta(H,U,\cF)&=&\varphi_{\eta'}((H,U,\cF)_{S'}),\label{app11a}\\
\varphi_s(H,U,\cF)&=&\varphi_{s'}((H,U,\cF)_{S'}).\label{app11b}
\end{eqnarray}
By \ref{app10}, these numbers do not depend on the choice of $(S',\eta',s')$.

\begin{teo}[Deligne, Kato]\label{app12}
Let $(H,U,\cF)$ be an object of $\fF_S$ such that $\cF$ has rank $1$ over $U$, 
$x$ be the closed point of $H$, $u\colon U\rightarrow H_\eta$ be the canonical injection. Then we have
\begin{equation}\label{app12a}
\dim(\Psi^0_x(u_!\cF))-\dim(\Psi^1_x(u_!\cF))=\varphi_s(H,U,\cF)-\varphi_\eta(H,U,\cF)-2\delta(H).
\end{equation}
\end{teo}  

If $\cF$ is unramified at every point of $H_s^\circ$, Deligne (\cite{laumon2} 5.1.1) proved 
the theorem for sheaves of any rank. In the general case, Kato (\cite{kato3} 6.7) proved 
the theorem for sheaves of any rank, with another definition of the invariant $\varphi_s(H,U,\cF)$.
One of the authors (T. Saito) \cite{saito} gave another proof for sheaves of any rank, 
with yet another definition of the invariant $\varphi_s(H,U,\cF)$. The latter corresponds 
to a second formula announced by Kato (\cite{kato2} 4.5). If $\cF$ has rank one, 
the invariant $\varphi_s(H,U,\cF)$ in Kato's latter formula coincides with 
our definition (\cite{kato1} remark after 6.8). 

Notice that formula \eqref{app12a} holds also in the case where $S$ has unequal characteristic.

\section{Dimension of the local Fourier transform}\label{apdlft}

\subsection{}\label{apdlft1}
We fix an algebraically closed field $k$ of characteristic $p>0$,
a finite field $\Lambda$ of characteristic $\not=p$
and a non-trivial additive character $\psi_0\colon \mF_p\rightarrow \Lambda^\times$.
We denote by $\cL_{\psi_0}$ the Artin-Schreier locally constant sheaf of $\Lambda$-modules 
of rank $1$ over $\mG_{a,k}$ associated to $\psi_0$. 
Apart from this change of conventions, we keep the same notation as in §\ref{results}.
In particular, we consider the sheaves $\cL_{\psi_0}(x\cx)$ and 
$\ocL_{\psi_0}(x\cx)$ over $\rA\times_k\crA$ and $\rP\times_k\crP$ respectively. 
For a scheme $W$ over $\rA\times_k\crA$ 
(resp. $\rP\times_k\crP$), we denote also by $\cL_{\psi_0}(x\cx)$ (resp. $\ocL_{\psi_0}(x\cx)$)
the pull-back of $\cL_{\psi_0}(x\cx)$ (resp. $\ocL_{\psi_0}(x\cx)$) to $W$.

\subsection{}\label{apdlft2}
Let $X$ be a smooth connected curve over $k$, $f\colon X\rightarrow \rP$ be a $k$-morphism,
étale over a dense open subscheme of $X$, $Y=f^{-1}(\rA)$, $s\in X(k)$, $z=f(s)$, $\cz\in \crP(k)$.
We denote by $\cT$ and $H$ the henselizations of $\crP$ and $X\times_k\crP$ at $\cz$ and $(s,\cz)$ 
respectively, by $\ctau$ the generic point of $\cT$ and (abusively) by $\cz$ the closed point of $\cT$. 
We consider $H$ as a $\cT$-scheme by the morphism $H\rightarrow \cT$ induced by 
the canonical projection $X\times_k\crP\rightarrow \crP$.
We denote by $U$ the inverse image of $Y\times_k\crA$ in $H_\ctau=H\times_{\cT}\ctau$, and let
\begin{equation}
\rho(s,\cz)=\varphi_{\cz}(H,U,\cL_{\psi_0}(x\cx)),
\end{equation}
where $\varphi_\cz$ is the invariant defined in \eqref{app11b}. 

\begin{lem}\label{apdlft3}
Under the assumptions of \eqref{apdlft2}, we have 
\begin{equation}
\rho(s,\cz)=\left\{
\begin{array}{clcr}
-\ord_s(f^*dx)& {\rm if}\ (z,\cz)\in \rP\times \cinfty,\\
1&{\rm if}\ (z,\cz)=(\infty,\czero),\\
\end{array}
\right.
\end{equation}
where $\ord_s(f^*dx)$ is the order of the non-zero meromorphic differential form $f^*dx$ over $X$.
\end{lem}
The case where $(z,\cz)=(\infty,\czero)$ follows directly from the definition.
Assume that $\cz=\cinfty$. 
We put $y=\cx^{-1}$ and consider the base change $\cT_1\rightarrow \cT$ given by $\cT[y_1]/(y_1^p-y)$.
We denote by $\fp$ the generic point of the special fiber of the canonical
projection $X\times_k\cT_1\rightarrow \cT_1$,
by $R$ the completion of the local ring of $X\times_k\cT_1$ at $\fp$,
by $K$ the fraction field of $R$, and by $b$ the image of $x$ in $R$ (which is a unit). 
Since $f^*dx\not=0$, $b$ is not a $p$-th power in $R$.
By (\cite{aml} §10), the Swan conductor of $\cL_{\psi_0}(x\cx)$ 
at $\fp$ is $p$, and its refined Swan conductor is the class of the differential form 
\[
db \otimes[y_1^{-p}]\in \Omega^1_{R}(\log)\otimes_R (y_1R)^{-p}/(y_1R)^{-p+1}.
\]
Moreover, $\cL_{\psi_0}(x\cx)$ is trivialized by a stable extension of $K$, 
namely the extension $L$ of $K$ defined by the equation $t^p-t=b/y_1^p$.
Indeed, the integral closure of $R$ in $L$ is generated over $R$ by $t_1=y_1t$ 
which satisfies the equation $t_1^p-y_1^{p-1}t_1=b$. The lemma follows.

\subsection{}\label{apdlft4}
We keep the assumptions of \ref{apdlft2}. Moreover, let $Y_0$ be a dense open subscheme
of $X$ contained in $Y$, $j\colon Y_0\rightarrow X$ be the canonical injection, 
$\cG$ be a locally constant sheaf of $\Lambda$-modules of rank $1$ over $Y_0$. 
We denote by $\pr_1\colon X\times_k\cT\rightarrow X$ the first projection,
and consider the complex of nearby cycles 
\[
\Psi(\pr_1^*(j_!\cG)\otimes\ocL_{\psi_0}(x\cx))
\]
in $\bD_c^b(X,\Lambda)$, relatively to the second projection $\pr_2\colon X\times_k\cT\rightarrow \cT$.  
We fix an algebraic closure of $k(\ctau)$ and denote by $\theta$ the associated geometric point of $\cT$. 
We consider the sheaf $\cG_{\theta}\otimes \cL_{\psi_0}(\cx f_{\theta})$
over $X_\theta=X\times_k\theta$ (cf. \ref{results1} for the notation).

\begin{prop}\label{apdlft5}
We keep the assumptions of \eqref{apdlft2} and \eqref{apdlft4}. 

{\rm (i)}\ If $s\in Y-Y_0$ and  $\cz=\cinfty$, then we have 
\[
\dim(\Psi^1_s(\pr_1^*(j_!\cG)\otimes\ocL_{\psi_0}(x\cx)))=\sw_s(\cG)+1+\ord_s(f^*dx).
\]

{\rm (ii)}\ If $(z,\cz)=(\infty,\cinfty)$, then we have 
\[
\dim(\Psi^1_s(\pr_1^*(j_!\cG)\otimes\ocL_{\psi_0}(x\cx)))=\sw_{s\times \theta}(\cG_{\theta}\otimes\cL_{\psi_0}(\cx f_{\theta}))+1+\ord_s(f^*dx).
\]

{\rm (iii)}\ If $(z,\cz)=(\infty,\czero)$, then we have 
\[
\dim(\Psi^1_s(\pr_1^*(j_!\cG)\otimes\ocL_{\psi_0}(x\cx)))=\sw_{s\times \theta}(\cG_{\theta}\otimes\cL_{\psi_0}(\cx f_{\theta}))-\sw_s(\cG).
\]
\end{prop}

This follows from \ref{app12} and \ref{apdlft3}.

\begin{prop}[\cite{laumon} 2.4.3]\label{apdlft6}
Let $z\in \rP(k)$, $\cz\in \crP(k)$, $T$ and $\cT$ be the henselizations of $\rP$
and $\crP$ at $z$ and $\cz$ respectively, $\tau$ and $\ctau$ be the generic points of $T$
and $\cT$ respectively. Let $\cF$ be a constructible sheaf
of $\Lambda$-modules over $\tau$, of rank $\rk(\cF)$ and Swan conductor $\sw(\cF)$, 
$\Theta(\cF)$ be the set of slopes of $\cF$. Then the rank of the local Fourier transform
of $\cF$ at $(z,\cz)$ \eqref{results2} is given by 
\begin{equation}\label{apdlft6a}
\rk(\fF^{(z,\cz)}(\cF))=\left\{
\begin{array}{clcr}
\sw(\cF)+\rk(\cF)& {\rm if}\ (z,\cz)\in \rA\times \cinfty,\\
\sw(\cF)-\rk(\cF)& {\rm if}\ (z,\cz)=(\infty,\cinfty)\ {\rm and} \ \Theta(\cF)\subset ]1,\infty[,\\
0& {\rm if}\ (z,\cz)=(\infty,\cinfty)\ {\rm and} \ \Theta(\cF)\subset [0,1],\\
\rk(\cF)-\sw(\cF)& {\rm if}\ (z,\cz)=(\infty,\czero)\ {\rm and} \ \Theta(\cF)\subset [0,1[,\\
0& {\rm if}\ (z,\cz)=(\infty,\czero)\ {\rm and} \ \Theta(\cF)\subset [1,\infty[.\\
\end{array}
\right.
\end{equation} 
\end{prop}

By Brauer induction, we may reduce to the case where $\cF=f_*(\cG)$, 
$f\colon S\rightarrow T$ is a finite morphism, étale above $\tau$, 
$S$ is the spectrum of a henselian discrete valuation ring, with generic point $\eta$,
and $\cG$ is a constructible sheaf of $\Lambda$-modules of rank $1$ over $\eta$. 
There exist a connected smooth curve $X$ over $k$, a $k$-morphism $\tf\colon X\rightarrow \rP$, 
a point $s\in X(k)$, a dense open subscheme $Y_0$ of $X$, and a locally constant constructible
sheaf of $\Lambda$-modules of rank $1$, $\tcG$, over $Y_0$, such that $S$ is isomorphic 
to the henselization of $X$ at $s$, $z=\tf(s)$, $f$ is induced by $\tf$, $\tf(Y_0)\subset \rA$,
and $\cG$ is isomorphic to the pull-back of $\tcG$ to $\eta$. 
We take again the notation of \eqref{apdlft2} and \eqref{apdlft4} (applied to $\tf$ and $\tcG$).
It follows from (\cite{sga7-2} XIII 2.1.7.1 and 2.1.7.2) that we have a canonical isomorphism
\begin{equation}\label{apdlft6b}
\fF^{(z,\cz)}(\cF)\simeq \Psi^1_s(\pr_1^*(j_!\tcG)\otimes\ocL(x\cx)).
\end{equation}

Let $R$ be the completion of the local ring of $S$, $K$ be its fraction field, $t$ be a uniformizer
of $R$, $b$ be the image of $x$ in $K$, $\ord$ be the valuation of $K$ normalized by $\ord(t)=1$. 
For $v\in K$, we put $v'=\frac{dv}{dt}\in K$. 
If $(z,\cz)\in A\times \cinfty$, then by \ref{apdlft5} and (\cite{serre1} VI §2), we have 
\begin{eqnarray}\label{apdlft6c}
\rk(\fF^{(z,\cz)}(\cF))&=&\sw(\cG)+1+\ord(b')\\
&=&\sw(\cG)+\ord(\frac{tb'}{b-x(z)})+\ord(b-x(z))\nonumber\\
&=&\sw(\cF)+\rk(\cF).\nonumber
\end{eqnarray}
We fix an algebraic closure of $k(\ctau)$ and let $\theta$ be the associated geometric point of $\cT$, 
and $T_{\{\theta\}}$ be the henselization of $T\times_k\theta$ at $z\times_k\theta$.
We denote by a subscript $\{\theta\}$ the objects deduced from objects over $T$ by the base change 
$T_{\{\theta\}}\rightarrow T$. Similarly as for \eqref{apdlft6c}, we have 
\begin{eqnarray}\label{apdlft6d}
\rk(\fF^{(\infty,\cinfty)}(\cF))&=&\sw(\cG_{\{\theta\}}\otimes f_{\{\theta\}}^*(\cL_{\psi_0}(x\cx)))+1+\ord(b')\\
&=&\sw(\cG_{\{\theta\}}\otimes f_{\{\theta\}}^*(\cL_{\psi_0}(x\cx)))+\ord(tb(b^{-1})')+\ord(b)\nonumber\\
&=&\sw(\cF_{\{\theta\}}\otimes \cL_{\psi_0}(x\cx))-\rk(\cF),\nonumber
\end{eqnarray}
\begin{eqnarray}\label{apdlft6e}
\rk(\fF^{(\infty,\czero)}(\cF))&=&\sw(\cG_{\{\theta\}}\otimes f_{\{\theta\}}^*(\cL_{\psi_0}(x\cx)))-\sw(\cG)\\
&=&\sw(\cF_{\{\theta\}}\otimes\cL_{\psi_0}(x\cx))-\sw(\cF).\nonumber
\end{eqnarray}
Let $\kappa$ be a geometric generic point of $T_{\{\theta\}}$, $I=\pi_1(\tau,\kappa)$,
$I_{\{\theta\}}=\pi_1(\tau_{\{\theta\}},\kappa)$, $I^{(a)}$ and  $I^{(a)}_{\{\theta\}}$
$(a\in \mQ_{\geq 0})$ be the classical logarithmic ramification filtrations of respectively $I$ and 
$I_{\{\theta\}}$ (\cite{serre1} IV, cf. \cite{as} for the notation).
For every $a\in \mQ_{\geq 0}$, the canonical surjective homomorphism 
$I_{\{\theta\}}\rightarrow I$ identifies $I^a$ with the image of $I_{\{\theta\}}^a$. 
We consider the slope decomposition of the representation $\cF_\kappa$ of $I$
\[
\cF_\kappa=\oplus_{\lambda\in \Theta(\cF)}\cF_{\kappa,\lambda}.
\]
By (\cite{laumon} 2.1.2.7), to conclude the proof of the proposition, it is enough to show that 
\begin{equation}\label{apdlft6f}
(\cF_{\kappa,1}\otimes \cL_{\psi_0}(x\cx)_{\kappa})^{I^{(1)}_{\{\theta\}}}=0.
\end{equation}
Recall that $I^{(1)}/I^{(1+)}$ is an $\mF_p$-vector space. By (\cite{aml} 14.3 and 14.4), 
we have an isomorphism 
\begin{equation}
{\rm Hom}_{\mZ}(I^{(1)}/I^{(1+)},\mF_p)\simeq k,
\end{equation}
and similarly for $I^{(1)}_{\{\theta\}}/I^{(1+)}_{\{\theta\}}$ (in {\em loc. cit.}, we trivialize the line
$N_{-1}$ by $x^{-1}$). Since we fixed a non-trivial character $\psi_0\colon \mF_p\rightarrow \Lambda^\times$, 
the action of $I^{(1)}/I^{(1+)}$ on $\cF_{\kappa,1}$ 
determines a finite set of characters $I^{(1)}/I^{(1+)}\rightarrow \mF_p$, and hence 
a finite set of points $\Sigma\subset k$. Similarly, the action of $I^{(1)}_{\{\theta\}}/I^{(1+)}_{\{\theta\}}$ 
on $\cL_{\psi_0}(x\cx)_{\kappa}$ determines the point $\cx \in k(\theta)$ (cf. \cite{aml} 9.13).  
Since $\cx\not\in \Sigma$, equation \eqref{apdlft6f} follows.

\end{document}